\documentclass[12pt,a4paper]{article}
\usepackage[english]{babel}
\usepackage{amsfonts,amssymb,amsmath, epsfig,cases}
\usepackage{float,color,graphicx,graphics,psfrag}
\usepackage{subcaption}
\usepackage{amstext,amscd}
\allowdisplaybreaks
\usepackage{hyperref}
\usepackage{amsthm}
\usepackage{dsfont}

\usepackage{mathtools}
\mathtoolsset{centercolon}

\usepackage{rotating}

\usepackage{booktabs}
\usepackage{threeparttable}

\usepackage{ifthen}

\usepackage{pgfplots,pgfplots}
\pgfplotsset{compat=1.14}
\usepgfplotslibrary{fillbetween}

\newboolean{usetikz}
\setboolean{usetikz}{false}
\newboolean{tikzextern}
\setboolean{tikzextern}{false}

\ifthenelse{\boolean{tikzextern}}{
    \usetikzlibrary{external}
    \tikzexternalize[prefix=imgs/]
}{}

\usepackage[margin=2cm]{geometry}

\usepackage{newpxtext,newpxmath}

\usepackage[draft,margin=false,inline=true,author=]{fixme}
\fxusetheme{colorsig}

\FXRegisterAuthor{cg}{acg}{Chiara}
\FXRegisterAuthor{cs}{acs}{Christian}
\FXRegisterAuthor{gj}{agj}{Gaspard}

\usepackage{enumerate}
\usepackage{xspace}
\usepackage{stackrel}

\def\paragraph#1{{\bf #1\ }}

\newcommand{\Eel}{E_{\mathrm{el}}}
\newcommand{\Ep}{E_{\mathrm{p}}}

\newcommand{\En}{E_{\mathrm{n}}}
\newcommand{\Ec}{E_{\mathrm{c}}}

\newcommand{\ffric}{F_\mathrm{fric}}
\newcommand{\Fcomp}{F_{\mathrm{comp}}}

\newcommand{\T}{\mathbb{S}^1}
\newcommand{\spint}{\smashoperator{\int_{\T}}}

\newcommand{\sgint}{\smashoperator{\int_{\Gamma_n}}}

\newcommand{\dd}{\mathrm{d}}
\newcommand{\e}{\mathrm{e}_\theta}

\definecolor{light_grey}{rgb}{0.90, 0.90, 0.90}
\definecolor{custom_yellow}{rgb}{0.90, 0.62, 0.00}
\definecolor{custom_blue}{rgb}{0.00, 0.45, 0.70}
\definecolor{custom_lightyellow}{rgb}{0.98, 0.92, 0.50}
\definecolor{custom_lightblue}{rgb}{0.80, 0.85, 1.00}
\newcommand{\plotlinewidth}{1.5pt}

\title{The influence of nucleus mechanics in modelling adhesion-independent cell migration in structured and confined environments.}
\author{C. Giverso$^{1}$, G. Jankowiak$^{2}$, L. Preziosi$^{1}, $C. Schmeiser$^{3}$}
\date{\today}

\begin{document}
\maketitle
    \begin{center}
    \small
1-Department of Mathematical Sciences, Politecnico di Torino, Corso Duca degli Abruzzi 24, 10129 Torino, Italy
\\
2-Department of Mathematics and Statistics, University of Konstanz, 78457 Konstanz, Germany
\\
3-Faculty of Mathematics, University of Vienna, Oskar-Morgenstern Platz 1, 1090 Wien, Austria
\end{center}

\begin{abstract}

Recent biological experiments \cite{lammermann, Reversat2020, Balzer} have shown that certain types of cells are able to move in structured and confined environment even without the activation of focal adhesion.
Focusing on this particular phenomenon and based on previous works \cite{Jankowiak_M3AS}, we derive a novel two-dimensional mechanical model, which relies on the following physical ingredients: the asymmetrical renewal of the actin cortex supporting the membrane, resulting in a backward flow of material; the mechanical description of the nucleus membrane and the inner nuclear material; the microtubule network guiding nucleus location; the contact interactions between the cell and the external environment.
The resulting fourth order system of partial differential equations is then solved numerically to conduct a study of the qualitative effects of the model parameters, mainly those governing the mechanical properties of the nucleus and the geometry of the confining structure. Coherently with biological observations, we find that cells characterized by a stiff nucleus are unable to migrate in channels that can be crossed by cells with a softer nucleus. Regarding the geometry, cell velocity and ability to migrate are influenced by the width of the channel and the wavelength of the external structure. Even though still preliminary, these results can be potentially useful in determining the physical limit of cell migration in confined environment and in designing scaffold for tissue engineering.
\end{abstract}

{\bf{keywords}: } cell migration, mathematical modelling, friction-based migration,
focal adhesion, cytoskeleton

\section{Introduction}
\label{sec:introduction}

Cell migration onto two dimensional (2D) substrates and inside three dimensional (3D) environment plays an essential role in many physiological and pathological processes, including embryonic development, wound healing, immune response, cancer progression and metastasis formation \cite{wolf_molecular_2006, bergert_force_2015, Trepat}.
The unconfined motion of cells onto 2D extracellular matrix (ECM) is a well-studied process and it is conventionally described by a continuous and highly coordinated cyclic processes; the elongation of protrusions at the leading edge driven by actin polymerization, the formation of integrin-mediated focal adhesions (FAs), myosin-mediated contraction and the detachment of the trailing edge \cite{Trepat, Balzer, abercrombie_locomotion_1970}.
This classical description requires that specific transmembrane adhesion proteins (integrins, among others) carry intracellular forces from the cytoskeleton to the substrate to propel the cell forward \cite{bergert_force_2015, rafelski_crawling_2004,vicente-manzanares_integrins_2009} and it is therefore called {\it adhesion-dependent migration} or {\it integrin-mediated migration}.

While this mechanism of motion is well understood, the physical challenges that cells have to face when moving in 3D environments are only now getting more attention, and recent researches indicate that in vivo cell migration can substantially deviate from migration on 2D unconfined substrates \cite{Lammerding_nucleus, Balzer}.
Indeed, during motion through tissues, ECM barriers, capillaries and lymph nodes, cells experience varying degrees of physical confinement and cell migration can thus be achieved with very different mechanisms \cite{lammermann, Balzer, bergert_force_2015, cell_migration_3D}.
In particular, it was biologically observed that cell migration inside 3D environments can occur even in the absence of focal adhesions, indicating that additional mechanisms for adhesion and migration are possible \cite{lammermann, Reversat2020}.

Such {\it adhesion-independent migration} has been observed in confining 3D ECM environments \cite{lammermann, friedl_amoeboid_2001, friedl_biology_2000, friedl_amoeboid_2001, Fraley, ONeill}, using different cell lines and technologies.
For cells of different types (dendritic cells in \cite{lammermann}; leukocytes in \cite{Reversat2020}; breast carcinoma, pancreatic carcinoma, and human osteosarcoma cells in \cite{Balzer}), it was observed that migration in 3D environments may not require myosin-mediated contraction and that inhibitors of integrins do not hamper migration through channels leading to cell confinement, although these treatments can hinder and even prevent motility in wider channels leading to unconfined migration.
Considering leukocyte migration \cite{Reversat2020}, it was shown, on one hand, that these do not migrate when confined between two parallel plates, in absence of adhesion. On the other hand, leukocyte adhesion-free motion is possible when supporting pilars or microfabricated structured channels are placed between these plates, under the conditions that the pilar size and spacing ---or the characteristic length of the side wall structure--- match roughly that of the cell length  \cite{Reversat2020}.
Even though the origin and transmission of propelling forces during focal adhesion-free migration are not fully understood, all these findings 
indicate that 2D but not 3D migration is integrin-dependent and that adhesion-free motility relies on a structured physical confinement that can induce cytoskeletal alterations reducing the dependence of cell motion on the adhesion-contraction force coupling.
In the absence of adhesions, non-specific transient interactions between transmembrane proteins and the substrate might generate friction, converting the protrusive actin cortex flow into cell movement \cite{bergert_force_2015, hawkins_spontaneous_2011, lammermann}.
Furthermore, it has been observed that confined migration depends largely on microtubule (MT) dynamics and might persists even when F-actin is disrupted \cite{Balzer, Stroka, Li}.

Although an increasing level of confinement can trigger transitions from integrin-based towards adhesion-independent migration modes in many cell types \cite{friedl_amoeboid_2001}, in the absence of matrix protease production, a too strong confinement either decreases or even prevents migration, due to cell stiffness \cite{wolf_2003, wolf_physical_2013}.
In particular, while the cytoplasm is very flexible and the cytoskeleton can actively remodel to undergo large deformations and penetrate small openings, the cell nucleus is normally 2 to 10 times stiffer than the surrounding cytoplasm and, with a typical diameter of  $3-10 \, \mu m$, occupies a large fraction of the cellular volume and is usually larger than many of the pores encountered in the extracellular environment \cite{Lammerding_nucleus, wolf_physical_2013, Cao}. Thus, the nucleus should undergo substantial deformations when the cell moves through 3D constrictions, and it may constitutes a rate-limiting factor during non-proteolytic migration of cells \cite{Lammerding_nucleus, wolf_physical_2013}.

To understand the bio-physical and mechanical factors involved in the process of cell migration, many mathematical models have been proposed in the past decades \cite{Jilkine, Holmes, Kruse, Mogilner}. Specifically, there have been
abundant works related to cell migration onto 2D substrates, either modeling the membrane mechanics and its signaling activity \cite{Elliott, Hecht} or describing in detail the cytosol dynamics \cite{Levine_2010, Recho_PRL, Recho_JMPS, manhart_extended_2015}.
However, coupled models, including the cytosolic machinery and membrane dynamics, have received little attention, even though they are
critical to understand cell migration \cite{Kruse, Mogilner,  Giverso_PhysRevE, Moure_Gomez2018}.
Furthermore, most of these models have focused on 2D adhesion-dependent cell motility, in which cells extend a stationary lamellipodium at the leading edge. Even in models accounting for amoeboid motion and which have been extended to model 3D confined migration \cite{Moure_Gomez, Moure2017, Moure_Gomez2018}, the cell motility substantially relies on adhesion, on acto-myosin protrusion-contraction, and on cell capability to sense an external field though membrane receptors.
On the contrary, adhesion-independent migration inside constrained 3D environments has received less attention and mathematical models have started to tackle this interesting mechanism only recently.
In particular, a simplified two-dimensional model for focal adhesion-independent cell swimming, based on the flow-friction driven force transmission, has been proposed by \cite{Wu_Othmer, Othmer}, while in \cite{kaoui_lateral_2008} the motion of closed phospholipid membranes suspended in a nonlinear shear gradient of a plane Poiseuille flow was investigated numerically in two dimensions.
A possible explanation to the chemical signalling activity regulating adhesion-independent migration was proposed in \cite{Elliott}, using a system of reaction-diffusion equations and assuming a Turing instability to model polymerization pattern on the cell surface, which drives the formation of pseudopod.
However, in all these cases \cite{kaoui_lateral_2008, Wu_Othmer, Othmer, Elliott, Moure2017, Moure_Gomez2018}, the presence of the nucleus as a limiting factor for cell migration and the effect of confinement are not taken into account.

The influence of nuclear deformations on the whole process of cell migration was included in some recent works.
Specifically, Cao et al. \cite{Cao} develop a chemo-mechanical model to study the nuclear strains and shapes, its plastic deformation and the threshold for the rupture
of its envelope, during migration through confined interstitial spaces.
In \cite{Lee2017}, a 2D model for cell migration through a dense network of host cells was proposed to reproduce glioma cell invasion. The moving cell is represented by two elastic closed curves, an inner curve corresponding to the nucleus of the cell and an outer curve corresponding to the cell basal membrane, whereas non-moving cells were represented by a single elastic curve.
In \cite{Vermolen}, the deformations of the cell and nucleus during invasion through a dense microenvironment were simulated incorporating stochastic processes and uncertainties in the input variables were evaluated using Monte Carlo uncertainty quantification simulations.
These models \cite{Lee2017, Vermolen, Cao} were able to reproduce correctly the hourglass cell and nucleus deformation observed in biological experiments, by relying on an external chemical factor.

In this work, we build on top of \cite{Jankowiak_M3AS} to develop a simplified framework to study whether adhesion-free migration could be driven by simple mechanical features.
We also test the influence of cell nucleus mechanical properties in the determination of the physical limit of cell migration.
Thus, we propose a model of force generation during adhesion-independent cell migration in confined environment, taking into account the flow-friction driven force transmission, the cell membrane polymerization and the nuclear deformations.
The cell is modelled by two membranes, an outer one representing the cell membrane and an inner one representing the nucleus. The two membranes are connected by microtubules, responsible of nucleus location inside the cell. The renewal of the actin network underneath the cell membrane is modelled by the evolution of the mass distribution along the membrane, with a source (resp. sink) term at the front (resp. back), while also taking into account the conservation of the center of mass. The model is used to simulate cell motion inside channels with structured walls with wave-lengths ranging in the order of magnitude of the cell and nucleus diameters.
We present the mathematical model in Section \ref{sec:model} and the numerical scheme we used to solve the systems of equations in Section \ref{sec:discretization}. Finally, in Section \ref{sec:resultsDiscussion}, we present and discuss the numerical results.

\section{The mathematical model}
\label{sec:model}
In this section, we present the continuum model ingredients for adhesion-free cell migration in domains containing rigid obstacles with a given geometry.
Motivated by the experimental setup of \cite{Reversat2020}, where the cell is confined between flat top and bottom surfaces and structured side walls, we choose a two-dimensional model.
Then, we consider the cell composed by two main compartments, the cytoplasm and the nucleus, both surrounded
by membranes. The nuclear and cellular membranes can be represented as closed curves of $\mathbb{R}^2$.
The cell wall (cortex) is schematically represented by a lipid bilayer and a complex underlying network of actin filaments. It is assumed to be elastic and is subjected to a pressure differential force acting in its exterior normal direction.
 The renewal of the actin network composing the cellular cortex is a key ingredient during many cellular behaviours and is here modelled by deposit and removal of material along the cortex.
 In rough biological terms, this corresponds to the polymerization and depolymerization of some parts of the actin filaments and gives the cell a preferential direction of movement. The actin is then transported in the cell to a new location where it polymerizes again (actin treadmilling). This transport mechanism needs to be taken into account to ensure conservation of momentum in the absence of external forces. In our model, this is done by including an additional reaction force, as detailed below.

 Finally, the nuclear membrane can be thought as a double phospholipidic bilayer with some resistance to bending and tension and it is eventually subjected to the differential pressure between the cytosol and the interior of the nucleus.

\subsection{Description of the model and notation}
\begin{figure}[ht]
    \def\svgwidth{0.75\columnwidth}
    \begin{center}
        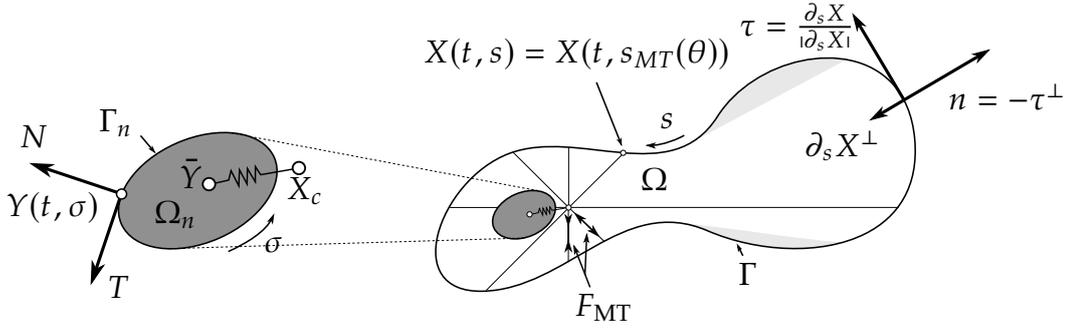
    \end{center}
    \caption{The parameterization and associated vector quantities, along with the representation of the microtubule network geometry, with the centrosome $X_c$ linked to the centroid of the nucleus $\bar{Y}$. The light shaded areas on the right are the ``blind'' areas, the boundary of which is not in $\Pi_{MT}(\T)$. This region would be empty if $\Gamma$ were star-shaped w.r.t $X_c$.
    The bottom arrow indicates the orientation of the curve. The elastic link between the centrosome and the nucleus is also represented, as well as the microtubular force $F_\mathrm{MT}$, both of which will be detailed later on.}
    \label{fig:cell parameterization}
\end{figure}

We set ourselves in $\mathbb{R}^2$ and consider the cell cortex and the nuclear membrane represented by the time-dependent closed curves $\Gamma(t)$ and $\Gamma_n(t)$ respectively,
\begin{align*}
    \Gamma(t) = \left\{X(t,s) : s\in\T\right\}, \text{ where } X(t, s) : [0, T] \times \T \rightarrow \mathbb{R}^2\,, \\
    \Gamma_n(t) = \left\{ Y(t,\sigma) : \sigma \in\T\right\}, \text{ where } Y(t, \sigma) : [0, T] \times \T \rightarrow \mathbb{R}^2\,,
\end{align*}
for some fixed maximal time $T$. For any fixed $t\in[0,T]$, $X(t,\cdot)$ and  $Y(t,\cdot)$ are defined on $\T$ (the unit circle in $\mathbb{R}^2$), which enforces closedness of the curves. The time derivatives are denoted by $\partial_t X$ and $\partial_t Y$ whereas the space derivatives (along the curves) are $\partial_s X$ and $\partial_{\sigma} Y$, respectively.
The variable $s$ corresponds to the density of actin along the cell membrane and therefore tracks Lagrangian particles. The total amount of actin on the cortex is normalized to~$1$ and for any non empty interval $[s_1, s_2]$, the normalized amount of actin on the corresponding piece of cortex is $|s_2 - s_1|$. Note that this can be generalized to a total amount of actin different from $1$.
On the other hand, the variable $\sigma$ should be regarded as a geometrical parameter.
Assuming that the curves are smooth enough, we denote by $\tau(t, s)$ and $n(t, s)$ the unit tangent and unit normal vectors to the cell membrane curve at $X(t, s)$ and with $T(t, \sigma)$ and $N(t, \sigma)$ the unit tangent and unit normal vectors to the nuclear membrane curve at $Y(t, \sigma)$.
Assuming positive orientation of the parametrization, we have
\begin{eqnarray*}
   \tau = \frac{\partial_s X}{|\partial_s X|} \,,\qquad n = -\tau^\bot \,, \qquad T = \frac{\partial_\sigma Y}{|\partial_\sigma Y|} \,,\qquad N = -T^\bot \,,
\end{eqnarray*}
with the convention $(a,b)^\bot = (-b,a)$.
A sketch of the notations employed in the paper is illustrated in Figure~\ref{fig:cell parameterization}. \\
Finally, let $\Omega(t) \subset \mathbb{R}^2$ be the set bounded by $\Gamma(t)$ and $\Omega_n(t) \subset \mathbb{R}^2$ the set bounded by $\Gamma_n(t)$, \emph{i.e.} $\Gamma(t) = \partial\Omega(t)$ and $\Gamma_n(t) = \partial\Omega_n(t)$: for biological consistency we have to guarantee that the nucleus is located inside the cell, i.e., $\Omega_n \subset \Omega$.

The nucleus and the cell cortex are coupled through the microtubule (MT) structure. Microtubules are composed by polymers of tubulin and constitute part of the cytoskeleton. They provide structure and shape to cells and are fundamental in a number of cellular processes. MTs are organized by the microtubule organizing center (MTOC), such as the centrosome found in many animal cells. The centrosome is connected to the cell nucleus and it plays a fundamental role in determining cell migration direction and persistence and in locating the cell nucleus \cite{fruleux_physical_2016, gundersen_nuclear_2013, beadle_role_2008,friedl_nuclear_2011}.
In particular, during migration, the nucleus is generally located at the back of the cell, but the centrosome (or microtubule organising center) is mostly seen at the front of the nucleus \cite{fruleux_physical_2016}, and the MTs are physically directed toward the migratory front in confined migration \cite{Balzer}.

In a continuous setting, the MT structure is fully described by setting the geometry of the microtubules (i.e., their shape and endpoints) and the distribution, or density, of the microtubules.
In this work, MTs are modelled as line segments connecting the centrosome, located in $X_c(t)$, to points on the cortex and we assume that they are homogeneously distributed around the centrosome.
In other words, the MTs form a rigid structure. 
We thus define, additionally to $X_c$, the angular velocity of the MT structure, $\omega$.

With that in mind, it is possible to define, for every time $t$, the microtubule anchoring points on the membrane cortex
as the first intersection of the half line starting at $X_c(t)$ with angle $\theta \in \T$ with the cell cortex. Formally, we define the map
\begin{equation*}
    \Pi_{MT}(t, \theta) = X_c(t) + \lambda_{MT}(t, \theta) \mathrm{e}_{\theta} \,,
\end{equation*}
where
\begin{equation*}
\lambda_{MT}(t, \theta) := \min \left\{\lambda \ge 0 : X_c(t) + \lambda \mathrm{e}_{\theta} \in \Gamma \right\} \text{ and }  \mathrm{e}_{\theta} = \begin{pmatrix} \cos(\theta) \\ \sin(\theta)\end{pmatrix}.
\end{equation*}
This map $\Pi_{MT}$ is well defined if $X_c(t) \in \Omega(t)$ and is surjective if $\Omega(t)$ is star-shaped with respect to $X_c(t)$.
Then, the cortex anchoring point
of the microtubule of angle $\theta$ is defined as
\begin{equation}
    X(t, s_{MT}(t, \theta)) := \Pi_{MT}(t, \theta)\,.
\end{equation}
Note that, with a slight abuse of notation, $\Pi_{MT}$ can also be seen as a map on $\T$ such that $\Pi_{MT}(t,\theta) = s_{MT}(t,\theta)$.
The segment $[X_c(t), X(t, s_{MT}(\theta))]$ then represents the microtubule at a given instant of time. The construction is illustrated in Figure~\ref{fig:cell parameterization}. In the following, the variables will be omitted whenever possible.

\subsection{Evolution of the cell membrane} \label{sec:evolution of the cortex}

For what concerns the evolution of the cell membrane, we refer to the model proposed in \cite{Jankowiak_M3AS}, properly modified in order to take into account the presence of the nucleus and the MT structure.
We thus introduce the arc length $\ell$, which is given by
\begin{equation}
\ell(t, s) = \int_0^s |\partial_s X(t, \sigma)|d\sigma \,.
\label{s2l}
\end{equation}
Equation \eqref{s2l} defines a map between the arc length $\ell$ and the Lagrangian variable $s$ which can be inverted in order to obtain the linear actin density $\rho(t, \ell)$:
\begin{equation*}
   s(t, \ell) = \int_0^\ell \rho(t, l)\;dl \,,
\end{equation*}
implying $\dfrac{d\ell}{ds} = |\partial_s X| = \rho^{-1}$.

The evolution of the actin density describes the active component of the model due to the heterogeneity of the polymerization rate across the cortex, is given by
\cite{Jankowiak_M3AS}:
\begin{equation*}
    \partial_t \rho(t, \ell) = f(t, \ell)\,,
\end{equation*}
where $f(t,\ell)$ is the rate of actin increase ($f(t,\ell)>0$) or decrease ($f(t,\ell)<0$).

As done in \cite{Jankowiak_M3AS}, we assume that the total amount of actin in the cortex is kept constant and that the cell polarization and subsequent actin polymerization manifests itself by a local imbalances producing a net increase of actin close to the cell front, balanced by a decrease close to the rear of the cell.
Since the total amount of actin in the cortex does not change in time, we require
\begin{equation*}
   \int_{\Gamma} f(t,\ell) \,d\ell = \spint f(t, \ell(t,s))|\partial_s X(t,s)|ds = 0 \,,\qquad t\ge 0\,.
\end{equation*}
The mass transfer rate, or polymerization rate, $r_\mathrm{pol} \ge 0$ is then defined as
\begin{equation}
    r_\mathrm{pol} \coloneq\int_\Gamma f^+\; d\ell = -\int_\Gamma f^- \; d\ell\,,
    \label{eq:def polymerization rate}
\end{equation}
where $\left(\,\cdot\,\right)^\pm$ denotes the positive and negative part.

\begin{figure}[ht]
    \def\svgwidth{0.5\columnwidth}
    \begin{center}
\begingroup%
  \makeatletter%
  \providecommand\color[2][]{%
    \errmessage{(Inkscape) Color is used for the text in Inkscape, but the package 'color.sty' is not loaded}%
    \renewcommand\color[2][]{}%
  }%
  \providecommand\transparent[1]{%
    \errmessage{(Inkscape) Transparency is used (non-zero) for the text in Inkscape, but the package 'transparent.sty' is not loaded}%
    \renewcommand\transparent[1]{}%
  }%
  \providecommand\rotatebox[2]{#2}%
  \newcommand*\fsize{\dimexpr\f@size pt\relax}%
  \newcommand*\lineheight[1]{\fontsize{\fsize}{#1\fsize}\selectfont}%
  \ifx\svgwidth\undefined%
    \setlength{\unitlength}{315.68867841bp}%
    \ifx\svgscale\undefined%
      \relax%
    \else%
      \setlength{\unitlength}{\unitlength * \real{\svgscale}}%
    \fi%
  \else%
    \setlength{\unitlength}{\svgwidth}%
  \fi%
  \global\let\svgwidth\undefined%
  \global\let\svgscale\undefined%
  \makeatother%
  \begin{picture}(1,0.38245813)%
    \lineheight{1}%
    \setlength\tabcolsep{0pt}%
    \put(0,0){\includegraphics[width=\unitlength,page=1]{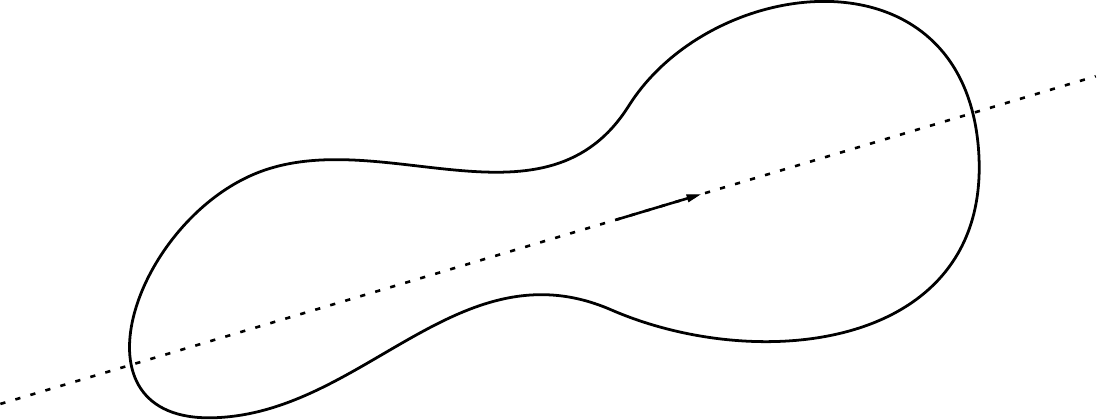}}%
    \put(0.5885713,0.15241657){\makebox(0,0)[lt]{\lineheight{1.25}\smash{\begin{tabular}[t]{l}$e_p$\end{tabular}}}}%
    \put(-0.18819682,0.09191039){\makebox(0,0)[lt]{\lineheight{1.25}\smash{\begin{tabular}[t]{l}$X(t,s_\mathrm{back}(t))$\end{tabular}}}}%
    \put(0.91490974,0.23848865){\makebox(0,0)[lt]{\lineheight{1.25}\smash{\begin{tabular}[t]{l}$X(t,s_\mathrm{front}(t))$\end{tabular}}}}%
    \put(0,0){\includegraphics[width=\unitlength,page=2]{imgs/polyregions.pdf}}%
  \end{picture}%
\endgroup%

    \end{center}
    \caption{The front and back of the cell, with the corresponding polymerization regions for a possible choice of $f$. The dark yellow and light yellow area highlight the support of $\tilde{f}_\mathrm{front}$ and $\tilde{f}_\mathrm{back}$, respectively.}
    \label{fig:polymerization regions}
\end{figure}

In practice, we assume that the cell is polarized in a given direction $e_p \in \mathbb{S}^1$, and that for each time $t$, there are unique $s_\mathrm{back}(t)$, $s_\mathrm{front}(t)$ so that
\[
    X(t, s_\mathrm{back}(t)) \cdot e_p(t) = \min_s X(t, s) \cdot e_p\,,
    \quad
    X(t, s_\mathrm{front}(t)) \cdot e_p(t) = \max_s X(t, s) \cdot e_p\,.
\]
This is illustrated in Fig.~\ref{fig:polymerization regions}.
A reasonable choice for $f$ is then a function with its (non negative) maximum at $X(t,s_\mathrm{front}(t))$ and (non positive) minimum at $X(t,s_\mathrm{back}(t))$.
For example
\begin{equation}
    f(t, \ell) = r_\mathrm{pol} \left( \tilde{f}_\mathrm{front}(t, d(X(s(\ell),t), s_\mathrm{front}(t)))) - \tilde{f}_\mathrm{back}(t, d(X(s(\ell), t), s_\mathrm{back}(t)))) \right)\,,
    \label{eq:poly source term}
\end{equation}
where $d$ denotes the distance along the cortex, and $\tilde{f}_\mathrm{front}$, $\tilde{f}_\mathrm{back}$ are non negative functions, normalized so that \eqref{eq:def polymerization rate} is fulfilled.
This explains their explicit dependency on time.

For what concerns the MT endpoints density on the cell cortex, we can define
\begin{equation}
    \label{eq:MT density}
    \rho_{MT}(s) =   \left| (\Pi_\text{MT}^{-1})'(s)\right|  =
    \begin{cases}
        \left| \frac{d \theta}{d s} \right| =   \left| \partial_s X \right| \, \dfrac{n\cdot (X(s) - X_c)}{|X(s) - X_c|^2} & \text{ if }X(s) \in \Pi_{MT}(\T)
        \\
        0 & \text{otherwise.}
    \end{cases}
\end{equation}
To derive eq. \eqref{eq:MT density}, we exploited the fact that
 $\dfrac{d\ell}{ds} = |\partial_s X| = \rho^{-1}$ and that
    \begin{equation*}
    \frac{dX(s)}{d \theta} \cdot \e^\perp = |X(s) - X_c|
    \Leftrightarrow
    \frac{dX(s)}{d \theta} \cdot (X(s)-X_c)^\perp = |X(s) - X_c|^2\,,
\end{equation*}
leading to
\begin{equation*}
    \frac{d \theta}{d\ell} = \frac{n\cdot (X(s) - X_c)}{|X(s) - X_c|^2}\,.
\end{equation*}

Then, the evolution of the cell cortex is determined by the Newton's Second Law, written in an overdamped regime
\begin{equation}
    \dfrac{D}{Dt} X- k_\tau v(\partial_t X_c, \omega, \theta) \rho_{MT} = - F_{MT} \rho_{MT}  - F_\mathrm{cont} + F_\mathrm{wall}+ \Fcomp + F_{c} \,, \label{eq:cell_membrane}
\end{equation}
where
the friction coefficient in front of the time derivative was set equal to unity, without loss of generality, by an appropriate choice of the time scale, as done in \cite{Jankowiak_M3AS}.
The second term on the l.h.s. and the first term in the r.h.s. of eq. \eqref{eq:cell_membrane} represent, respectively, the friction force generated by microtubules binding complexes sliding on the cortex  and the in-line force due to MTs elongation (as they will be detailed in Section \ref{sec:evolution of the centrosome}).
The second and third term on the r.h.s.  of eq. \eqref{eq:cell_membrane} take into account the contact between the cell cortex and the nucleus and the one between the cell and the channel wall (see Section \ref{sec:evolution of the nuclear membrane} for a detailed description of $F_\mathrm{cont}$). In the following, we assume that $F_\mathrm{wall}$ derives from a potential, a possible choice of which is proposed in Section~\ref{sec:resultsDiscussion}.
We observe that, in eq. \eqref{eq:cell_membrane}, the velocity of the cell cortex relative to the laboratory coordinates is given by the material derivative $DX(t,s)/Dt$, defined as
\begin{align}
  \dfrac{D}{Dt} X(t,s) &= \partial_t X(t,s) + \left(\int_0^{\ell(t, s)} f(t,\lambda)d\lambda\right) \partial_s X(t,s)\, \nonumber\\
  &= \partial_t X(t, s) + \left(\int_0^s f(t, \ell(t, \sigma)) |\partial_s X(t, \sigma)| d \sigma\right)  \partial_s X(t,s)\,,
  \nonumber\\
  &= \partial_t X(t, s) + F_T(t,s)\, .
  \label{eq:material derivative}
\end{align}
The compensating force $F_\mathrm{comp}$ must be chosen (for each time $t$) so that the center of mass is fixed if the effects of the confinement, nucleus and microtubule struture are removed, see \cite{Jankowiak_M3AS} for more details on the derivation.
In this case, we obtain after integration
\begin{equation}
    \spint \Fcomp(s) \; ds = \spint F_T(s) \; ds = -\spint X(s) f(\ell(s)) |\partial_s X(s)| \; d\sigma = -\int_\Gamma X(s(\ell)) f(\ell) \; d \ell\,.
    \label{eq:fcomp balance}
\end{equation}
Finally, the last term on the r.h.s. of eq. \eqref{eq:cell_membrane} takes into account the forces acting on the cell due to the cortex mechanical behaviour and the pressure difference in and out of the cell. It can be obtained from the cell membrane energy
$$F_{c}=- \dfrac{\delta \Ec} { \delta X} \, , \qquad
\Ec= \Eel +\Ep \, .$$
The membrane energy comprises an elastic term representing the cell and cell cortex elasticity $\Eel$ and a term related to the existence of a differential pressure $\Ep$.
More precisely, the associated elastic energy $\Eel$ is composed by two parts, the first one $\Eel^{(1)}$ is related to the response of the cell cortex to tension and the second one $\Eel^{(2)}$ is related to the response of the cell to deviations from its target area.
Formally $\Eel^{(1)}$ reads:
\begin{equation}\label{Espring}
\Eel^{(1)} = \frac{1}{2}  k_c \spint \left(|\partial_s X|-1\right)^2 ds\,,
\end{equation}
where $k_c$ is the mechanical parameter representing the stretchability of the cell membrane.
Furthermore, we include an elastic constraint on the cell area  $\left|\Omega\right|$ (which would correspond to the volume in three-dimensions). In particular, we assume that the elastic energy of the cell is minimized if the cell area is equal to a given target area $A_c^*$, i.e.,
\begin{equation*}
    \Eel^{(2)} = \mu_c \left( \left|\Omega\right|  - A_c^* \right)^2= \mu_c \left( -\frac{1}{2}\int_{\mathbb{T}^1} X\cdot \partial_s X^\perp \, \dd s - A_c^* \right)^2\,,
\end{equation*}
where $\mu_c$ is the mechanical parameter representing the elastic resistance of the cell bulk to variations of its area.\\
For what concerns the term $\Ep$, we suppose that the cell is subject to an internal cytoplasmic pressure, which results in a force in the direction of the normal to the curve. The force intensity per unit length is assumed to be uniform in space and constant in time, so that the associated energy $\Ep$ is:
\begin{equation*} \label{Ep}
    \Ep = -p \left|\Omega\right| = \frac{p}{2} \spint X\cdot \partial_s X^\perp\, \dd s\,,
\end{equation*}
where $p>0$ is the constant excess of pressure inside the cell, with respect to the extracellular pressure.

\subsection{Evolution of the microtubules' structure and the centrosome} \label{sec:evolution of the centrosome}

Microtubules (MTs) are one of the coupling mechanisms between the cortex and the nucleus. As stated before, each MT is modelled as a line segment connecting the centrosome, located in $X_c(t)$
to one point on the cortex $s_{MT}(t,\theta)= \Pi_{MT}(t,\theta) $.
The microtubule force $F_{MT} (t, \theta)$ is directed along the microtubule (see Fig.~\ref{fig:cell parameterization}) and its intensity is a function of the microtubule length $\left| X_c(t) - X(t, s_{MT}(\theta))\right| $, to be specified, i.e.,
\begin{equation}
    F_{MT} (t, \theta) = - f\left( \left| X(t, s_{MT}(t,\theta))- X_c(t) \right| \right) \mathrm{e}_{\theta} \,,
\end{equation}
were $\mathrm{e}_{\theta} $ is the outward unit vector representing the direction of a MT.
This force acts on the portion of the cell cortex where $\rho_{MT} \neq 0$ and on the centrosome, where MTs originate.
The total resultant of all microtubule forces acting on the centrosome is, then,
\begin{equation}
    \bar{F}_{MT}(t)=  \spint F_{MT} (t, \theta) \, \dd \theta \,.
\end{equation}

In addition to the in-line forces due to elongation, the cortex endpoints of the microtubules are also subject to a friction force, caused by the sliding of the binding complexes on the cortex. This force is directed opposite to
the velocity of the centrosome relative to the cortex. The friction also generates a torque, that makes the microtubule structure undergoing a rigid body rotation with angular velocity $\omega$.\\
Thus, let us define $v(\partial_t X_c, \omega, \theta)$ the speed of the centrosome relative to the cortex endpoint of the microtubule
of angle $\theta$. We have
\begin{equation*}
  v(t, \partial_t X_c, \omega, \theta) := \partial_t X_c(t) - \partial_t X(t, s_{MT}(t, \theta)) + |X_c(t) - X(t, s_{MT}(t, \theta))|\omega \mathrm{e}_{\theta}^\perp \,.
\end{equation*}
We remark that in the relative velocity, we consider $\partial_t X(t, s_{MT}(t, \theta))$ as opposed to the total derivative $\frac{d}{dt} X(t, s_{MT}(t, \theta))$. This corresponds to the speed of the point on the cortex where the microtubule is attached, and not the end of the microtubule itself.\\
We also decompose $v$ along  $\tau$ and $\e$:
\begin{equation*}
v = v_\tau \tau + v_{\e} \e = (\varPi_\tau + \varPi_{\e}) v\,,
\end{equation*}
where $\Pi_\tau$ and $\Pi_{e_\theta}$ are the corresponding \emph{oblique} projections:
\begin{equation*}
    \varPi_\tau = \frac{\tau \otimes \e^\perp}{\tau \cdot \e^\perp}, \text{ and }
    \varPi_{\e} = \frac{\e \otimes \tau^\perp}{\e \cdot \tau^\perp}\,.
\end{equation*}
Locally on the cortex, the friction of the microtubule against the cortex generates a force
\begin{equation*}
\ffric = -(k_\tau \varPi_\tau + k_{\e} \varPi_{\e}) v\,,
\end{equation*}
$k_\tau$, $k_{\mathrm{e}_\theta}$ are friction coefficients in the corresponding directions.
In the following, for the sake of simplicity, we will take $k_{\e} = k_\tau$, but the model can be easily generalized to the case $k_{\e} \neq k_\tau$.

Calling $F_\mathrm{int}$ the forces acting on the microtuble structure, which include both the forces directed along the microtubules, $\bar{F}_{MT}$, and the sum of all forces acting on the centrosome itself, $\bar{F}_{X_c}$, we have the following force balance in $\mathbb{R}^2$
\begin{equation}
    \label{eq:centrosome force balance}
    -k_\tau \left[2\pi\partial_t X_c + \int  \left(\omega|X_c - X_\theta|\e^\perp - \partial_t X_\theta \right) \dd \theta\right] + F_\mathrm{int} = 0\,,
\end{equation}
where we used the short-hand $X_\theta = X(s_{MT}(\theta))$ and $F_\mathrm{int}= \bar{F}_{X_c}+\bar{F}_{MT}$. The expression for the force $\bar{F}_{X_c}$ will be detailed in section \ref{sec:evolution of the centrosome}.\\
Since the microtubule structure has zero moment of inertia w.r.t $X_c$, we also have balance of torques:
\begin{align*}
    \int (X_\theta - X_c) \times F_\mathrm{fric} \;\dd \theta = 0
    &\Leftrightarrow
    -k_\tau \int (X_\theta - X_c)^\perp \cdot v\; \dd \theta = 0
    \\
    & \Leftrightarrow
    -k_\tau \int |X_\theta - X_c| \e^\perp \cdot v\; \dd \theta = 0\,,
\end{align*}
which can be rewritten by decomposing $v$ on $\e$ and $\e^\perp$:
\begin{align}
    \label{eq:centrosome torque balance}
    0 &= -k_\tau \int |X_\theta - X_c| \left[\partial_t (X_c - X_\theta)\cdot \e^\perp + \omega |X_\theta - X_c|\e^\perp \right] \dd \theta\nonumber
    \\
    &= -k_\tau \left[\omega \int|X_c - X_\theta|^2 + \partial_t X_c \cdot \int |X_c - X_\theta|\e^\perp - \int |X_c - X_\theta| \partial_t X_\theta \cdot \e^\perp\right] \dd \theta \,.
\end{align}

Gathering \eqref{eq:centrosome force balance} and \eqref{eq:centrosome torque balance} we get the system
\begin{equation}
    \label{eq:friction first order condition}
    A \begin{pmatrix} \partial_t X_c \\ \omega \end{pmatrix} = B\,,
\end{equation}
where
\begin{align*}
    A &= k_\tau \begin{pmatrix}
        2 \pi  I_2 & \int |X_c - X_\theta|\e^\perp
        \\
        \left(\int |X_c - X_\theta|\e^\perp\right)^T & \int |X_c - X_\theta|^2
    \end{pmatrix}\,, \quad
    B = \begin{pmatrix}
        k_\tau \int \partial_t X_\theta + F_\mathrm{int} \\
        k_\tau \int |X_c - X_\theta|\e^\perp \cdot \partial_t X_\theta
    \end{pmatrix}\,,
\end{align*}
$I_2$ being the identity matrix in 2D.
We remark that, by Jensen's inequality, we have
\begin{equation*}
    \det A =  2\pi\, k_\tau^3 \left(2\pi\int |X_c - X_\theta|^2 - \left|\int |X_c - X_\theta|\e^\perp\right|^2\right) = 2\pi\, k_\tau^3\, \Delta \ge 0\,,
\end{equation*}
The equality occurs only if $X_\theta$ reduces to a single point, so we can assume that the determinant is always positive.
Then, system \eqref{eq:friction first order condition}

 can be solved to get the angular velocity $\omega$ of the MT structure and the velocity of the centrosome:
\begin{align}
    \omega &= \Delta^{-1}
    \left[-\left( k_\tau^{-1} F_\mathrm{int} + \int \partial_t X_\theta\right) \cdot \int |X_c - X_\theta|\e^\perp
    + 2\pi\int |X_c - X_\theta|\e^\perp \cdot \partial_t X_\theta\right]
    \\
    2\pi \partial_t X_c &= - \omega \int |X_c - X_\theta|\e^\perp + \int \partial_t X_\theta + k_\tau^{-1} F_\mathrm{int}\,.
    \label{eq:evo centrosome}
\end{align}

\subsection{Evolution of the nuclear membrane} \label{sec:evolution of the nuclear membrane}

We consider the nucleus material to be harder to deform than the rest of the cell and with a negligible relaxation time.
The nuclear membrane has a certain mechanical behaviour and it is connected to the centrosome, which contributes to the positioning of the nucleus inside the cell.
Furthermore, the nucleus cannot cross the cellular membrane and therefore interacts with it.
Then, calling $F_{n}$ the forces internal to the nuclear membrane, $F_{X_cN}$ the force acting on each point of the nuclear membrane due the interaction between the nucleus and the centrosome,
$F_\mathrm{cont,n}(\sigma)$ the contact force with the cortex, the Newton's Second Law for the nuclear membrane reads
\begin{align} \label{eq:nucleus_balance}
    F_{n} - F_{X_cN} +F_{\mathrm{cont},n} = 0 \,,
\end{align}
where
the meaning and derivation of the different terms in eq. \eqref{eq:nucleus_balance} is depicted below.\\
The force $F_{n}$ related to the mechanical behaviour of the membrane can be derived from the nuclear membrane energy, $\En$, through the relation
\[F_{n}=- |\partial_\sigma Y|^{-1} \delta \En / \delta Y\,,\]
where the energy $\En$ is
\begin{equation} \label{eq:nucleus_energy}
    \En = \dfrac{k_b}{2} \sgint (K-K_0)^2 \,\dd \ell + \lambda \sgint  \dd \ell + \Delta p_n \smashoperator{\int_{\Omega_n}}\,\dd A \, + \,  \mu_n \left( \smashoperator{\int_{\Omega_n}}\,\dd A - A_n^* \right)^2 \, ,
\end{equation}
where $\dd \ell= |\partial_\sigma Y| \dd \sigma$ is the arc length element and $\dd A $ is the area element.
The first term in eq. \eqref{eq:nucleus_energy} represents the energy associated to the membrane bending, $k_b$ being the bending modulus of the nuclear membrane, $K$ the local curvature and $K_0$ the characteristic (or spontaneous) curvature,
which is assumed to be zero in what follows, according to \cite{kaoui_lateral_2008}. The second term represents the tensile stress acting on the membrane and $\lambda$ can be thought as the nucleus surface tension.
$\Delta p_n$ is the difference between the pressure in the cytosol and the pressure inside the nucleus and the last term represents the volumetric elastic constraint associated to changes in nucleus area $A_n$ from a defined target area $ A_{n}^*$.

Taking all these contributions into account, the computation of the variation $\delta\En$ of the nuclear membrane energy function (see Section \ref{appendix:variation}) leads to the following nuclear membrane mechanical force
\begin{align*}
F_n =  k_b  \left(  \dfrac{\partial^2 K}{\partial \ell^2} + \dfrac{1}{2} K^3  \right) N  - \lambda  K \, N - \Delta p_n N  - \mu_n  \left( \spint Y \cdot \partial_{\sigma} Y ^\perp \, \dd \sigma  -A_{n}^* \right)  N \, .
\end{align*}
Then, we assume that the centrosome is linked to the nuclear membrane points, thanks to cytoskeletal filaments acting as a spring. This coupling is uniform along the length of the nuclear membrane.
Therefore, we assume a linear elastic link between the centrosome and the nuclear membrane, i.e.,
$$F_{X_cN}  = - k_e (X_c - Y ) \, ,$$
being $k_e$ the elasticity of these links.
Once a proper constitutive form is assumed for $F_{X_cN}$, it is also possible to define the total force $\bar{F}_{X_c}$ acting on the centrosome due to the centrosome-nuclear interaction, which is thus part of $F_\mathrm{int}$ in eq.  \eqref{eq:centrosome force balance}:
$$\bar{F}_{X_c}  =\spint F_{X_cN}  \left|\partial_{\sigma} Y(\sigma)\right|\, \dd \sigma  =- \spint k_e (X_c - Y ) \left|\partial_{\sigma} Y(\sigma)\right|\, \dd \sigma = - L_n k_e (X_c -\bar{Y}) \, ,$$
where $\bar{Y}$ is the centroid of the nucleus (see Fig.~\ref{fig:cell parameterization}) defined as
\begin{equation}
    \bar{Y}(t)= \dfrac{1}{L_n} \spint Y(t, \sigma) \left|\partial_{\sigma} Y(t, \sigma)\right|\, \dd \sigma \,,
\end{equation}
and $L_n$ is the length of the nuclear membrane at time $t$:
\begin{equation}
    L_n(t) = \spint \left|\partial_\sigma Y(t, \sigma)\right|\, \dd \sigma \, .
\end{equation}

Finally, as previously mentioned, the nucleus is constrained to remain within the cell, since it cannot cross the cortex. This is guaranteed by including a penalization.
The contact force is given by
\begin{equation}
\begin{aligned}
    F_{\mathrm{cont},n}(\sigma)
& = -\nabla_x\left[ \int_{\mathbb{T}^1} W_\mathrm{cont}(|x - X(s)|) \;ds\right]\bigg|_{x=Y(\sigma)}
\\
& = -\int \frac{Y(\sigma) - X(s)}{|Y(\sigma) - X(s)|}W'_\mathrm{cont}(|Y(\sigma) - X(s)|) \;ds\,.
\end{aligned}
\label{eq:fcont}
\end{equation}
where $W_\mathrm{cont}$ is a decreasing function with compact support.

By symmetry, we specify the corresponding force acting on the cortex, $F_\mathrm{cont}$, appearing in \eqref{eq:cell_membrane}.
\begin{equation}
    F_\mathrm{cont}(s) = \int \frac{Y(\sigma) - X(s)}{|Y(\sigma) - X(s)|}W'_\mathrm{cont}(|Y(\sigma) - X(s)|) \;d\sigma\,,
    \label{eq:fcont cortex}
\end{equation}
so that $\int F_\mathrm{cont}(s) \; ds = - \int F_{\mathrm{cont},n}(\sigma) \; d\sigma\,.$

\section{Numerical discretization} \label{sec:discretization}

The whole system composed by eqs. \eqref{eq:cell_membrane}-\eqref{eq:friction first order condition}-\eqref{eq:nucleus_balance} is discretized in space using finite differences,
the resulting system of ODEs is then solved using split step time stepping scheme of order one.
At each time step, the new position of the nuclear membrane is computed using the explicit scheme described below. Then the position of the cortex and centrosome
are updated using a semi-implicit scheme.

\subsection{Cortex, centrosome and microtubules structure}
We assume the total mass of actin on the cortex to be normalized to 1 and we introduce $N_1\in\mathbb{N}$ grid-points for the discretization of $s$ (corresponding to the cortex) such that $s_i = i\Delta s$,
$\Delta s =\frac{1}{N_1}$ for $i$ considered on a discrete torus, meaning that $i$ is identified with $i+N_1$. Given a time step $\Delta t>0$, we denote $t^j = j\Delta t$.
In what follows, subscripts (resp. supscripts) correspond to space (resp. time), and we define $X_i^j:=X(t^j, s_i)$. For the sake of legibility, indices or exponents are omitted if possible.
For integrals, we use the midpoint rule, meaning that
\[ \spint f(X(s)) \; ds \simeq \Delta s \sum_i f(X_i)\,, \text{ and } \int f(X(s(\ell))) \; d \ell \simeq \sum_i f(X_i) |X_{i+1} - X_{i-1}|/2\,. \]

We define $\tau_i = \frac{X_{i+1} - X_{i}}{|X_{i+1} - X_{i}|}$ and $n_i = \tau_i^\perp$.

The cortex elastic force, the bulk elasticity and pressure forces are discretized as
\[
    F_{c,i}^j = \frac{k_c}{\Delta s} \left( \left(\frac{|X_{i+1}^j - X_i^j|}{\Delta s} - 1\right) \tau_i^j - \left(\frac{|X_{i}^j - X_{i-1}^j|}{\Delta s} - 1\right) \tau_{i-1}^j \right)
    + \left(p - \mu_c (A^j - A_c^*)\right)\frac{X_{i+1}^j - X_{i-1}^j}{2\Delta s}\,,
\]
where the (polygonal) area $A^j$ is computed as
\[
    A^j = \sum_i \frac{1}{4} X_i^j \cdot (|X_{i+1}^j - X_i^j| n_i^j + |X_{i}^j - X_{i-1}^j| n_{i-1}^j)\,.
\]
The interaction force with the wall simply becomes
\[
    F_{\mathrm{wall},i} = -\nabla W_\mathrm{wall}(X_i)\,,
\]
whereas the transport term in~\eqref{eq:material derivative} reads
\[
    F_{T,i}^j = \tau_i^j \Delta s\sum_{k\le i} f_i^k\,.
\]
Following \eqref{eq:fcomp balance}, the compensating force $\Fcomp$ must satisfy
\[
    \sum_i F_{\mathrm{comp},i} = -\sum_i X_i f_i \frac{|X_{i+1} - X_{i-1}|}{2}\,.
\]
In the present work, we are mostly interested in the impact of the nucleus on the dynamics,
so we can choose $F_{\mathrm{comp},i}$ independent of $i$ for simplicity.

The computation of $\Pi_{MT}$ (or equivalently $s_{MT}(\theta)$) and $\rho_{MT,i}^j$ requires the construction of the visibility polygon of the cortex from the centrosome.
We do not discuss the construction here and refer to \cite{Lee_1983,JoeSimpson_1987} for details. The quadrature formulae and corresponding expression for $\Pi_{MT}$ are
detailed in Appendix~\ref{sec:discretization and quadrature formulae}.

It remains to deal with the contact force between the cell cortex and the nuclear membrane, which we simply take
as the discretization of \eqref{eq:fcont cortex}:
\[
    F_{\mathrm{cont},i} = \sum_k \frac{Y_k - X_i}{|Y_k - X_i|} W_\mathrm{cont}'(|Y_k - X_i|)\,.
\]

The equation \eqref{eq:friction first order condition} for the angular velocity $\omega$ of the MT structure and the velocity of the centrosome is discretizated similarly in time and space, using the quadrature formulae reported in the Appendix (see section \ref{sec:discretization and quadrature formulae}) to compute the integral of the different quantities.

For the time iteration, we use an implicit Euler scheme, so that for both the cortex and the centrosome we obtain a system of the following form at each time step:
\[
    \left[
    \begin{pmatrix}
        I_{2N_1} & &
        \\
                 & I_{2} &
        \\
        & & 0
    \end{pmatrix}
    + \Delta t
    \begin{pmatrix}
        M_{X,X} & M_{X_c,X}^T & M_{\omega,X}^T
        \\
        M_{X_c,X} & M_{X_c,X_c} & M_{X_c,\omega}
        \\
        M_{\omega,X} & M_{\omega,X_c} & M_{\omega,\omega}
    \end{pmatrix}
    \right]
    \begin{pmatrix}
        \vdots
        \\
        X_{i,1}^{j+1} - X_{i,1}^j
        \\
        \vdots
        \\
        X_{c}^{j+1} - X_{c}^j
        \\
        \omega^{j+1} - \omega^j
    \end{pmatrix}
    =
    \Delta t
    \begin{pmatrix}
        \vdots
        \\
        r_{i}^{j}
        \\
        \vdots
        \\
        r_{c}^{j}
        \\
        r_\omega^j
    \end{pmatrix}\,,
\]
where the $r_i, r_c, r_\omega$ on the r.h.s. correspond to the discretized terms detailed in this section, and the matrices $M$ are the corresponding jacobian matrices.
The left-most matrix above corresponds to the discretization of the time derivatives. Since the problem is formulated in term of the angular velocity $\omega$,
whose governing equation do not involve derivatives in time, the last element of the diagonal is zero.

\subsection{Nuclear membrane}
The evolution of the nuclear membrane is derived rewriting the proposed equation \eqref{eq:nucleus_balance} with the formulation and the efficient discretization scheme proposed in \cite{mikula_computational_2004, benes_comparison_2009}
appropriately adapted to our setting (see the Appendix, section \ref{appendix:Mikula} and \ref{App:Mikula_discretized}).
The evolution of the nucleus is then given by the following system
\begin{equation}
    \begin{cases}
    \partial_t K &= - k_b \left(\partial_\ell^4 K + \frac{1}{2} \partial_\ell^2(K^3)\right) + \partial_\ell (\alpha K) - K(K\beta + \partial_\ell \alpha) + \partial_\ell^2 (\nabla W(Y(\ell)) \cdot N) + \partial_\ell^2 (W(Y(\ell)) K)
    \\
     \partial_t \nu &= -k_b \partial_\ell^4 \nu - \frac{k_b}{2}\partial_\ell(\partial_\ell \nu)^3 + \partial_\ell(\nabla W(Y(\ell))\cdot N) + \partial_\ell (W \partial_\ell \nu ) + \alpha \partial_\ell \nu\
    \\
    \partial_t \eta &= K \beta + \partial_\ell \alpha
    \\
       \partial_t Y &= \alpha \partial_\ell Y - k_b \left(\partial_\ell^4 Y + \frac{3}{2} \partial_\ell(K^2 \partial_\ell Y)\right) + W(Y(\ell)) \partial_\ell^2 Y
                  \\& - \Delta p_n \, N  - \mu_n (A_n -A_{n}^*) N  - \left( \nabla W(Y(\ell)) \cdot N \right)  N  \,,
\label{eq:nucleus specialized_main}
\end{cases}
\end{equation}
where $\ell$ is the arc length, $K$, $\nu$ and $\eta = \log|\partial_\sigma Y|$ are, respectively, the curvature, the tangential angle and the logarithm of the local length of the nucleus envelope $\Gamma_n$ at a a point $Y \in \Gamma_n $, whereas $W(Y(\ell))$ is the potential combining the contact interaction with the cell cortex, the elastic constraint with the centrosome and the nuclear membrane surface tension.
In eq. \eqref{eq:nucleus specialized_main} we take
\begin{align}
   \partial_\ell \alpha &= -K\beta + \langle K \beta \rangle + (L_n/g-1) \zeta
    \label{eq:nucleus tangent velocity_main}
    \\
    \beta &= k_b \left(\partial_\ell^2 K + \frac{1}{2} K^3\right) - \Delta p_n  - \mu_n (A_n -A_{n}^*)  - \nabla W \cdot N - W K \,,    \label{eq:nucleus normal velocity_main}
\end{align}
where $\langle.\rangle = L_n^{-1} \int_{\Gamma_n} . \; \dd \ell$  is an averaging operator over the whole curve $\Gamma_n$, which makes \eqref{eq:nucleus tangent velocity_main} a non-local equation, whereas $\zeta >0$ is a given positive constant, included in order to avoid nodes to concentrates on points, leading to poor approximation and eventually inversion of ill-conditioned matrices.\\
To write the corresponding discretization of \eqref{eq:nucleus specialized_main}, we uniformly discretized the fixed parametrization interval $[0, 1]$
in $N_2$ subintervals, each of equal length $h=1/N_2$ and indexed by $i \in [0, N_2-1]$. For time, we use the same discretization introduced for the cell membrane, so that the point $Y(ih, j\Delta t)$ is written $Y_i^j$.
The measure of the finite element $[Y_{i-1}^j, Y_i^j]$ at time $t_j$, is given by $r_i = |Y_{i-1} - Y_i|$.
Then, the system of equations  \eqref{eq:nucleus specialized_main}-\eqref{eq:nucleus tangent velocity_main}- \eqref{eq:nucleus normal velocity_main} is solved for the discrete quantities $\alpha_i^j$, $\beta_i^j$, $K_i^j$, $\nu_i^j$, $\eta_i^j$, $Y_i^j$. In particular $\alpha_i^j$ denotes the tangential velocity of the node $Y_i^j$, whereas $\beta_i^j$, $K_i^j$, $\nu_i^j$, $\eta_i^j$ are piecewise constant approximations of the corresponding quantities on the finite element $[Y_{i-1}^j, Y_i^j]$.

The algebraic system determining the evolution of the nuclear envelope is reported in the Appendix \ref{App:Mikula_discretized}, along with some comments on the derivation of the discretized equations.
Finally, the discretized system of equations has been solved implementing a numerical code with Julia\footnote{\texttt{julialang.org}} \cite{Julia}.

\section{Numerical experiments} \label{sec:resultsDiscussion}
\subsection{Setup}
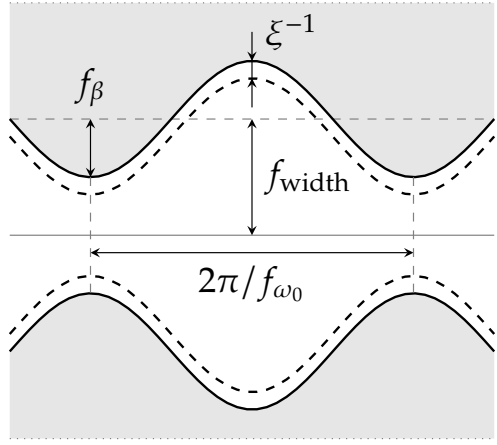
\begin{figure}[t]
 \centering
 \ifthenelse{\boolean{usetikz}}{
     \tikzsetfigurename{tikz_channel_}
\begin{tikzpicture}
    \begin{axis}[
            hide axis,
            ymax = 4,
        ]

        \addplot[dashed, color=gray, domain=-180:180] { 2 };
        \addplot[name path=top, dotted, color=gray, domain=-180:180] { 4 };
        \addplot[name path=bottom, dotted, color=gray, domain=-180:180] { -3.5 };
        \addplot[color=gray, domain=-180:180] { 0 };

        \addplot[name path=wall, thick, color=black, samples=50, domain=-180:180] {  2 + cos(1.5*x) };
        \addplot[thick, color=black, samples=50, domain=-180:180, dashed] {  1.7 + cos(1.5*x) };

        \addplot[thick, color=black, samples=50, domain=-180:180, dashed] { -1.7 - cos(1.5*x) };
        \addplot[thick, name path=dwall, color=black, samples=50, domain=-180:180] { -2 - cos(1.5*x) };

        \addplot[color=light_grey] fill between[
            of = wall and top,
          ];
        \addplot[color=light_grey] fill between[
            of = dwall and bottom,
          ];

        \draw[-stealth] (0,2.2) -- (0, 2.7);
        \draw (0,2.6) -- (0, 3.1);
        \draw[stealth-] (0,3) -- (0, 3.5) node[right] {$\xi^{-1}$};

        \draw[stealth-stealth] (0, 0) -- node[right]{$f_\mathrm{width}$} (0, 2);

        \draw[stealth-stealth] (-120, -0.3) -- node[below]{$2\pi/f_{\omega_0}$} (120, -0.3);
        \draw[dashed, color=gray] (-120, 1) -- (-120, -1);
        \draw[dashed, color=gray] (120, 1) -- (120, -1);

        \draw[stealth-stealth] (-120, 1) -- (-120, 2) node[above]{$f_\beta$};

    \end{axis}
\end{tikzpicture}}{
\includegraphics[width=0.45\textwidth]{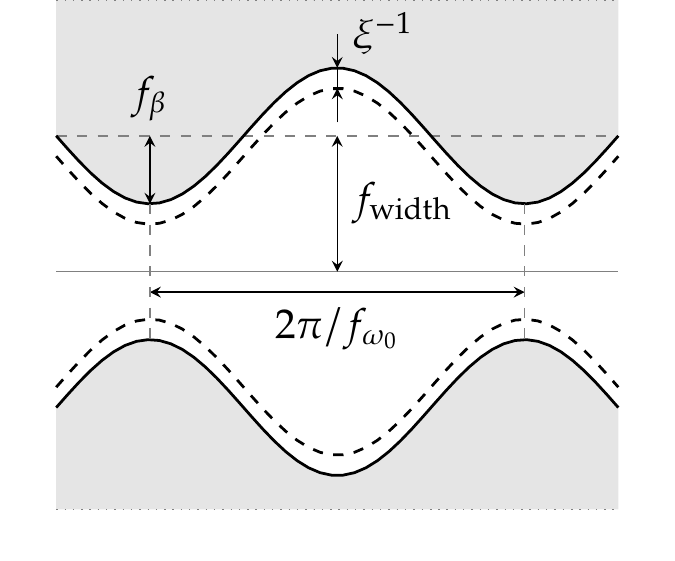}
}
{\footnotesize \caption{The sinusoidal channels used in the numerical simulation are described by the parameters: (i) $f_\mathrm{width}$, representing half of the mean width of the channel, (ii) $f_{\beta}$, which is the amplitude of the oscillation of the wall, (iii) $f_{\omega_0}$, which is the pulsation of the sinusoidal channel.  \label{Fig:channel}}}
\end{figure}

In this section, we present the numerical results obtained solving the discretized system of equations presented in Section \ref{sec:discretization}.
In particular, we consider channels with structured side walls of the following form:

\[
    \Omega_\mathrm{wall} = \left\{ (x, y) \in \mathbb{R}^2 : |y| \le f(x) =  f_{\beta}\, \sin ( f_{\omega_0} x) + \tfrac{1}{2} f_\mathrm{width} \right\}\,,
\]
where $f_\mathrm{width} $ represents half of the mean width of the channel, $f_{\beta} < f_\mathrm{width}$ is the amplitude of the oscillation of the wall, and $f_{\omega_0}$ is the pulsation of the sinusoidal channel (see Fig. \ref{Fig:channel} for an illustration of these quantities).
We remark that for $ f_{\beta}=0$, one recovers the flat walls geometry. We then choose
$$F_\mathrm{wall}(x, y) = \nabla \left[ g_{ \xi}(f(x) + y) + g_{ \xi}(f(x) - y) \right]\, ,$$
where $g_{ \xi}(x) = \min({ \xi} x - 1, 0)^2 \log({ \xi} x)$ is a smooth barrier function with
        $\lim_{x \rightarrow 0^+} g_{ \xi}(x) = +\infty$ and $g_{ \xi}(x) = 0$ for $x > { \xi}$.

Concerning the polymerization, we follow Section~\ref{sec:evolution of the cortex} and choose $\tilde{f}$ as a super-Gaussian: $\tilde{f} = \exp(-(\frac{x^2}{2w})^P)/C(t)$, where $C$ is a normalization factor. We choose $w = 0.5$ and $P = 3$.

Finally, we will consider the case in which the cell is represented
\begin{itemize}
\item only by the cell envelope, as done in \cite{Jankowiak_M3AS}, for comparison;
\item by the cell envelope and the cell nucleus, linked together by the microtubule  structure, as explained in the Section \ref{sec:model}.
\end{itemize}

 The initial condition for the cell cortex is chosen as the evenly spaced discretization of a closed curve $C_0$ which matches the side walls of the channel ---albeit with a smaller width--- in order to have an initial cell area equal to the target area $A_c^*$.
 More precisely, it is the union of the following 4 curves:
 \begin{itemize}
     \item $\{(x_0^\mathrm{min}, y) \in \mathbb{R}^2 : x_0^\mathrm{min} = -\frac{\pi}{2 f_{\omega_0}}, -f(x_0^\mathrm{min}) \le y \le f(x_0^\mathrm{min})\}$
     \item $\{(x, y) \in \mathbb{R}^2 : x_0^\mathrm{min} \le x \le x_0^\mathrm{max}, y = - f(x) + {\xi}\}$
     \item $\{(x, y) \in \mathbb{R}^2 : x_0^\mathrm{min} \le x \le x_0^\mathrm{max}, y = f(x) - { \xi}\}$
     \item $\{(x_0^\mathrm{max}, y) \in \mathbb{R}^2 : -f(x_0^\mathrm{max}) \le y \le f(x_0^\mathrm{max})\}$
 \end{itemize}
 where $x_0^\mathrm{max}$ is chosen so that the area enclosed by $C_0$ is $A_c^*$. The initial condition for the nucleus is the circle $\Gamma_{n,0}$ centered on $(\pi/2 f_{\omega_0}), 0)$ and such that $\Gamma_{n,0} + B_{{ \xi}^{-1}}(0) \subset \Omega_\mathrm{wall}$,
 where $B_{{\xi}^{-1}}(0)$ is the open ball of radius ${ \xi}^{-1}$ and $+$ denotes the Minkowski sum.

The spatio-temporal evolution of the cell and nuclear envelope for some benchmark simulations are reported in Section \ref{res1}, whereas the effect of the different parameters of the model on cell ability to move and its velocity inside sinusoidal channels is investigated in Section \ref{res2}.

\subsection{Results: to move or not to move}  \label{res1}
\begin{figure}[h]
 \centering
 \includegraphics[height=0.35\textwidth]{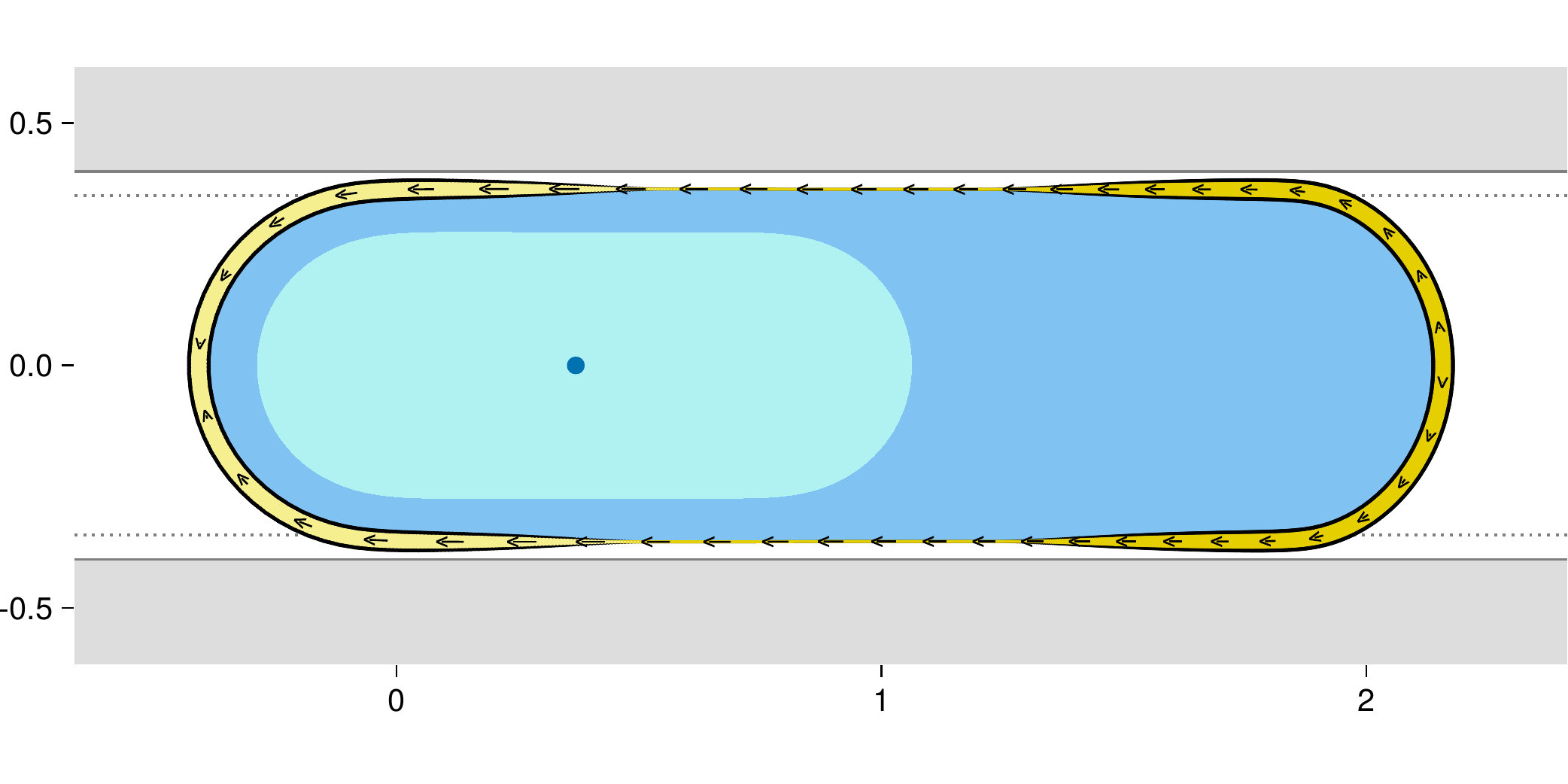}
 {\footnotesize \caption{Equilibrium configuration for $f_\beta = 0$, the values for the remaining parameters are presented in Table~\ref{table:numerical parameters}. The cell is in blue and the nucleus in light blue. The dark blue dot is the centrosome.
 The thickness of the yellow (resp. light yellow) region on the cell membrane indicate the strength of the polymerization (resp. depolymerization). The arrows indicate the flow of the cortex relative to the cell.
 \label{Fig:flatchannel}}}
\end{figure}

In this section, we investigate the ability of the model to reproduce cell migration inside a channel and the spatio-temporal evolution of the cell and the nucleus shapes. \\
We first consider the case of a cell positioned inside a channel with flat walls ($f_\beta = 0$): in this case, independently on the widths of the channel, no migration is observed in the simulations, in spite of polarization and corresponding cortex flow
(see Fig. \ref{Fig:flatchannel})

This result confirms that adhesion-free motility relies on structured confinement and it is in agreement with the experiments performed on leukocytes \cite{Reversat2020} and with the model described in \cite{Jankowiak_M3AS}, where the nucleus is not considered.

We then describe the motion of the cell inside a channel with sinusoidal walls.
In this configuration, the cell takes an hourglass configuration with nucleus deformation which is often observed in vivo,
so that channels of this type can be argued to mimic the inter- and extra- cellular space in which cells naturally migrate.

In this case, as it is known from biology, the capability of the cell to migrate inside the channel and its speed of migration are highly influenced by the presence of the nucleus, which is the stiffest part in the cell and can therefore remain stuck in the rear of the cell, preventing cell motion.
In the present model, the resistance of the nucleus to deformations is highly controlled by the bending modulus $k_b$ and by the elastic area-change constraint $\mu_n$.
Therefore, in Fig.~\ref{Fig:motion} (and in the videos in the Supplementary materials) we report the evolution of the cell and the nucleus shapes at the same instant of time, for some typical simulations, obtained varying the parameters $k_b$ and $\mu_n$. Namely, we consider
\begin{enumerate}[(a)]
\item the migrating cell without the cell nucleus, as done in  \cite{Jankowiak_M3AS} (Fig.~\ref{Fig:motion}a);
\item the migrating cell with a low bending modulus of the nucleus  (Fig.~\ref{Fig:motion}b, with ${k_b} = 10^{-2.5}, \mu_n = 50$);
\item the migrating cell with an intermediate bending modulus of the nucleus (Fig.~\ref{Fig:motion}c, with ${k_b} = 10^{-1.5}, \mu_n = 50$);
\item[(c')] the migrating cell with a higher value of the relaxation parameter $\mu_n$, and the same bending modulus as in (c) (Fig.~\ref{Fig:motion}c', with ${k_b} = 10^{-1.5}, {\mu_n} = 100$);
\item the non-migrating cell with a high bending modulus of the nucleus (Fig.~\ref{Fig:motion}d, with $k_b = 10^{-0.5}, \mu_n = 50$).
\end{enumerate}

\begin{sidewaysfigure}
 \centering
 \begin{tabular}{ccccc}
     \raisebox{1cm}{(a)}
     &\includegraphics[width=0.22\textwidth]{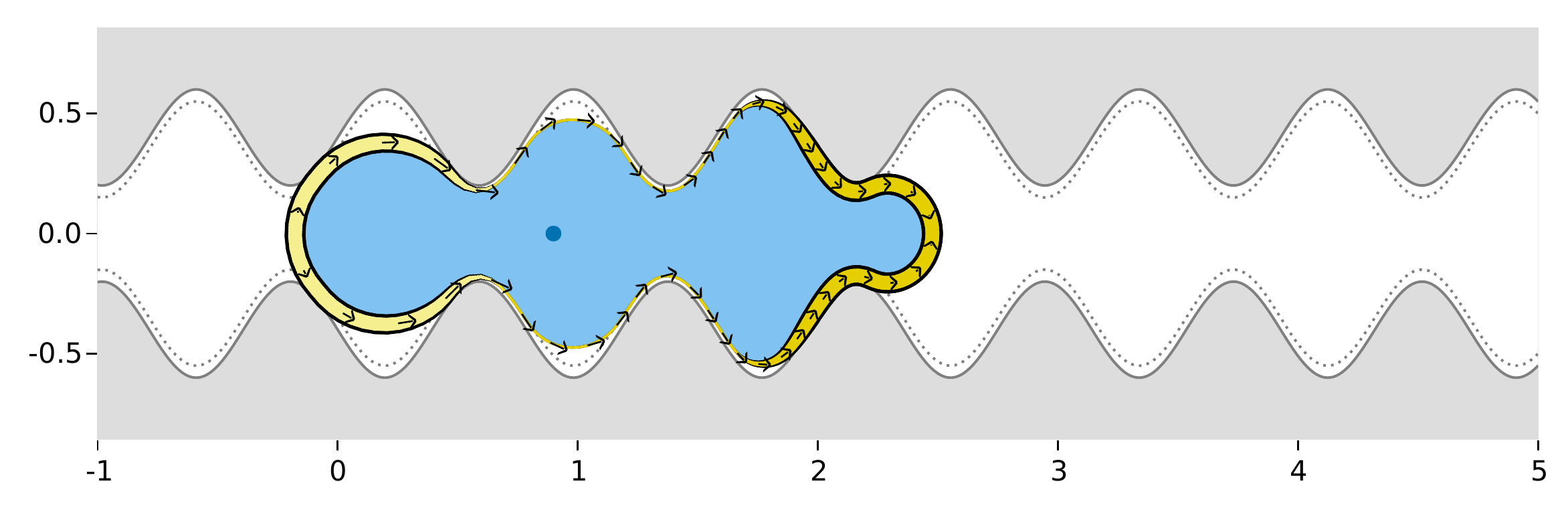}
        &\includegraphics[width=0.22\textwidth]{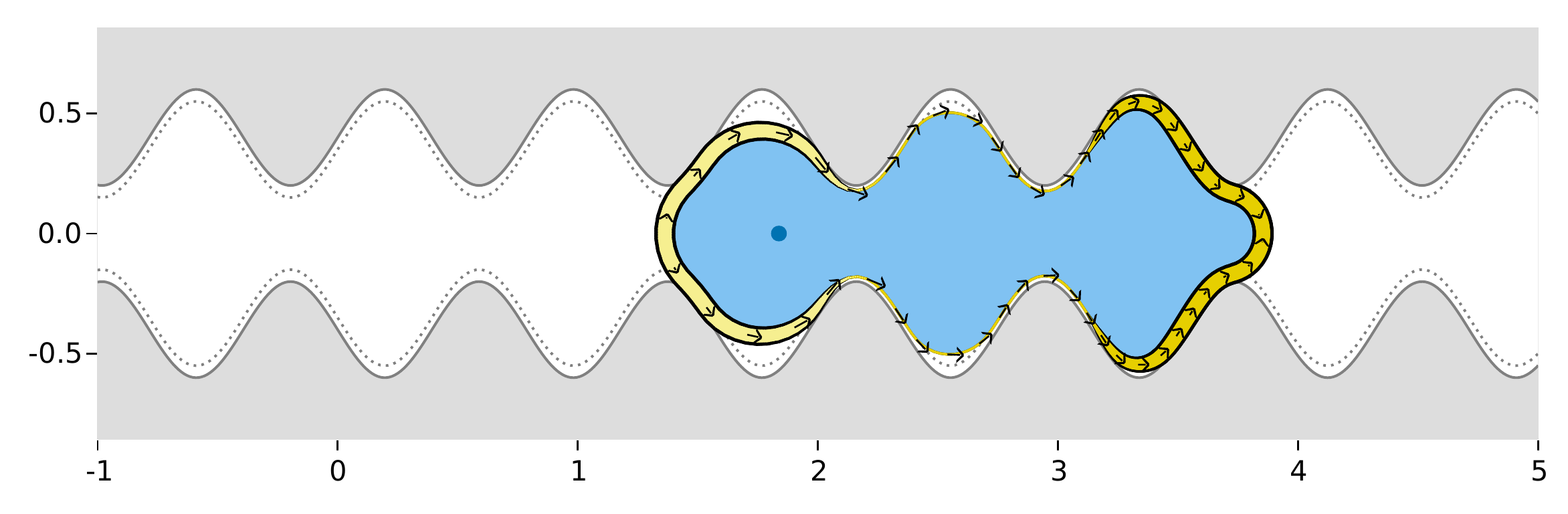}
        &\includegraphics[width=0.22\textwidth]{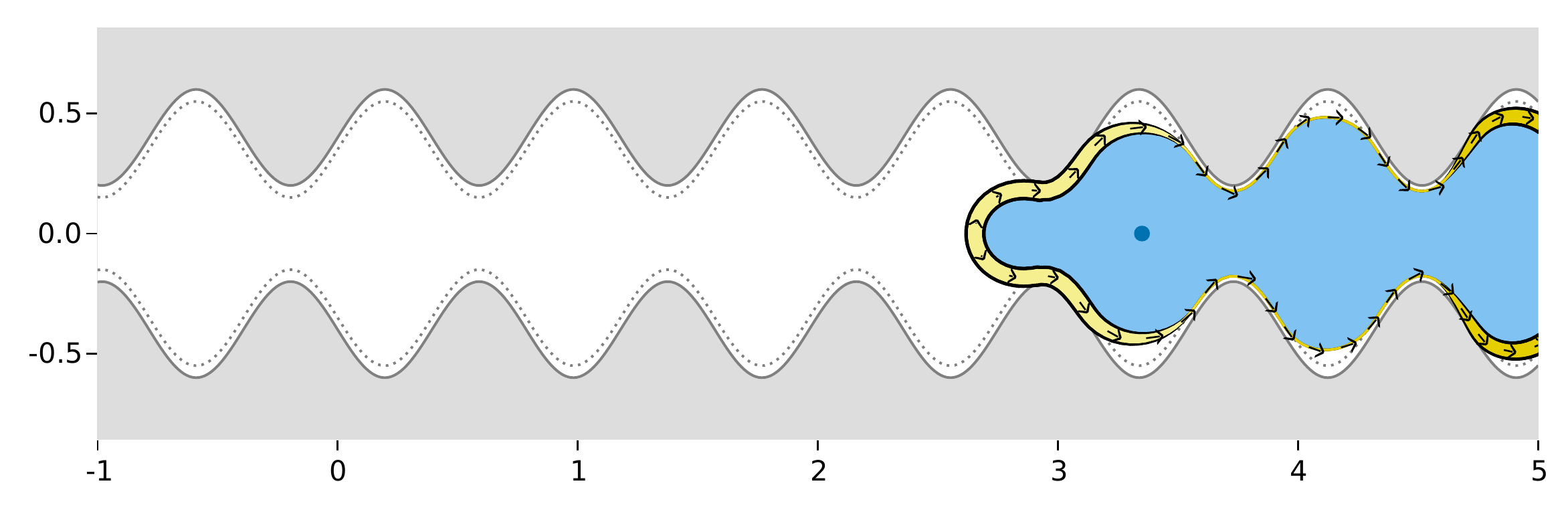}
        &\includegraphics[width=0.22\textwidth]{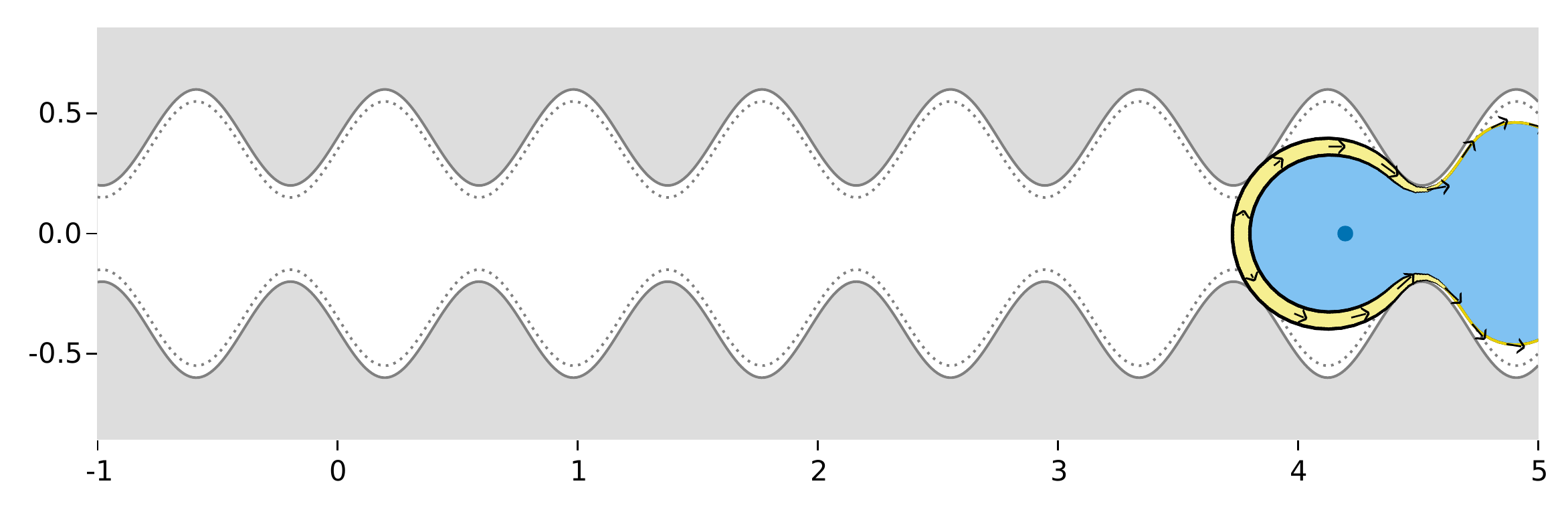}
\\
     \raisebox{1cm}{(b)}
        &\includegraphics[width=0.22\textwidth]{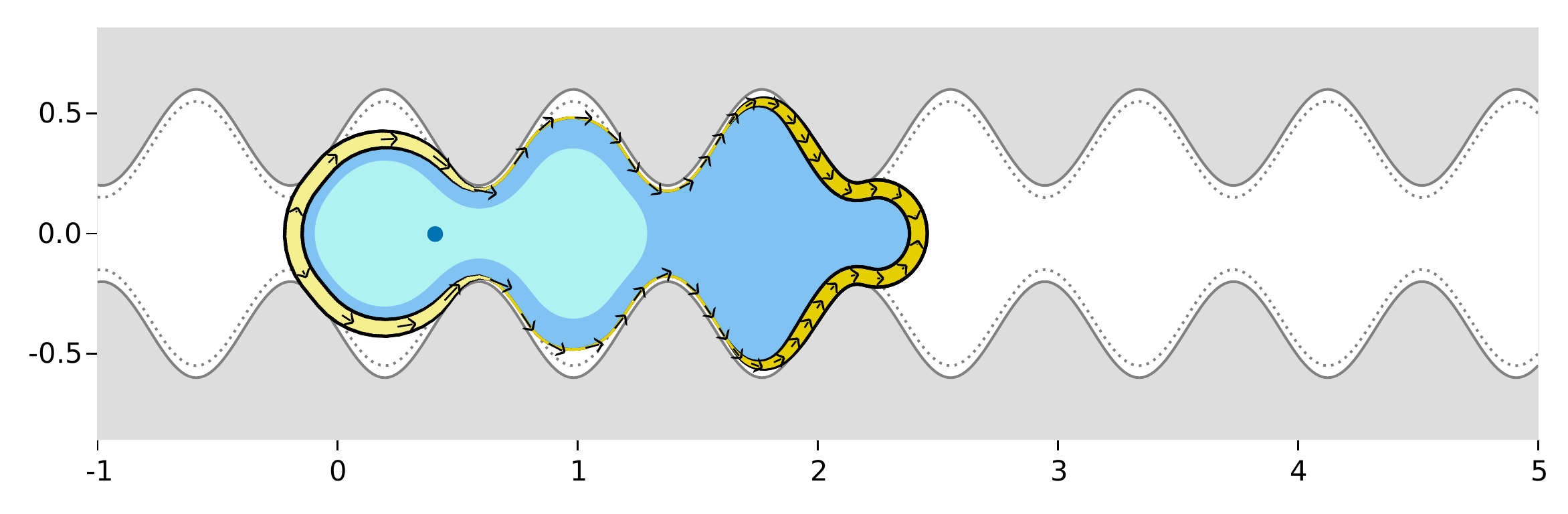}
        &\includegraphics[width=0.22\textwidth]{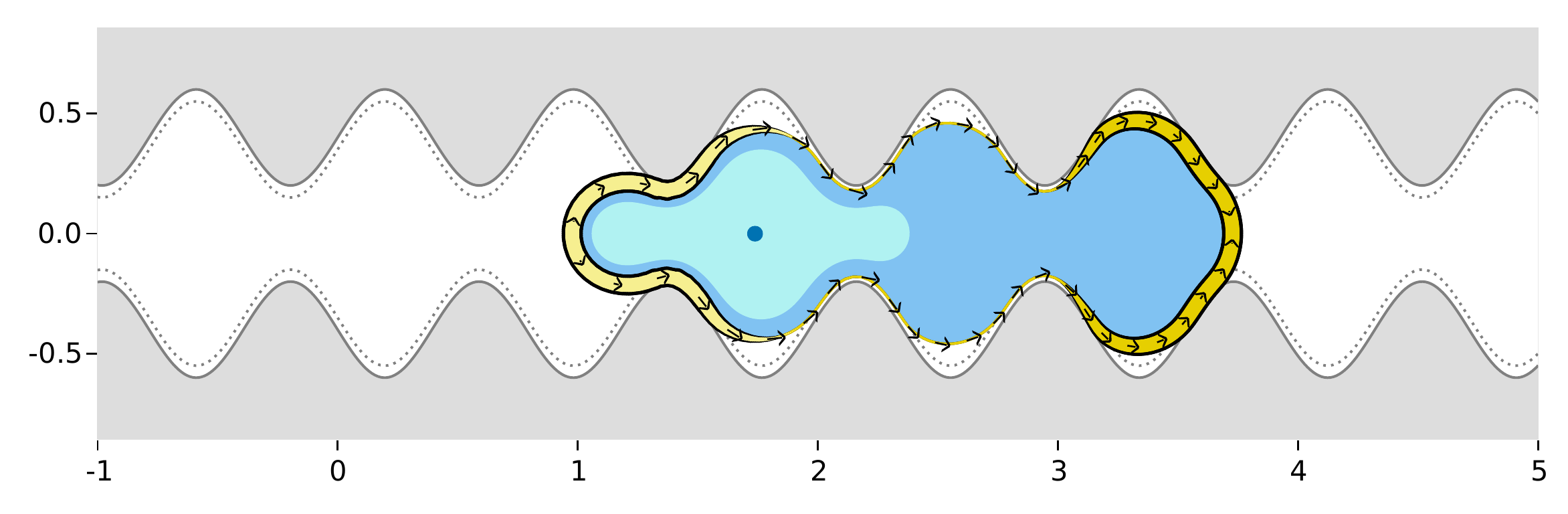}
        &\includegraphics[width=0.22\textwidth]{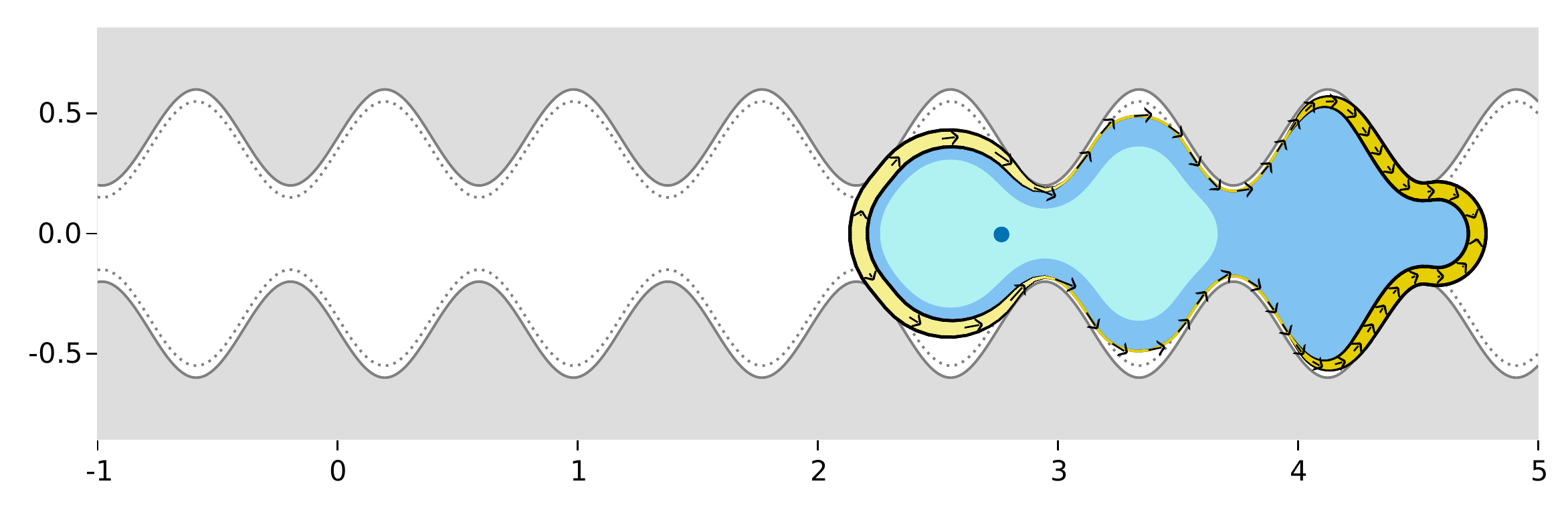}
        &\includegraphics[width=0.22\textwidth]{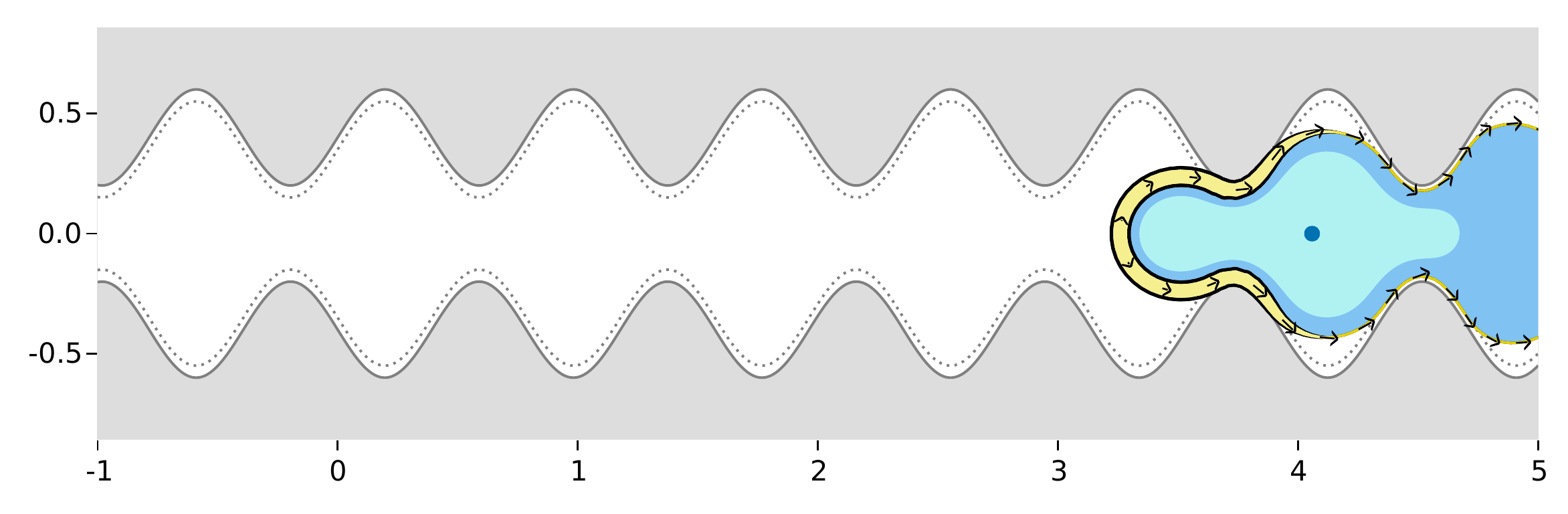}
\\
     \raisebox{1cm}{(c)}
        &\includegraphics[width=0.22\textwidth]{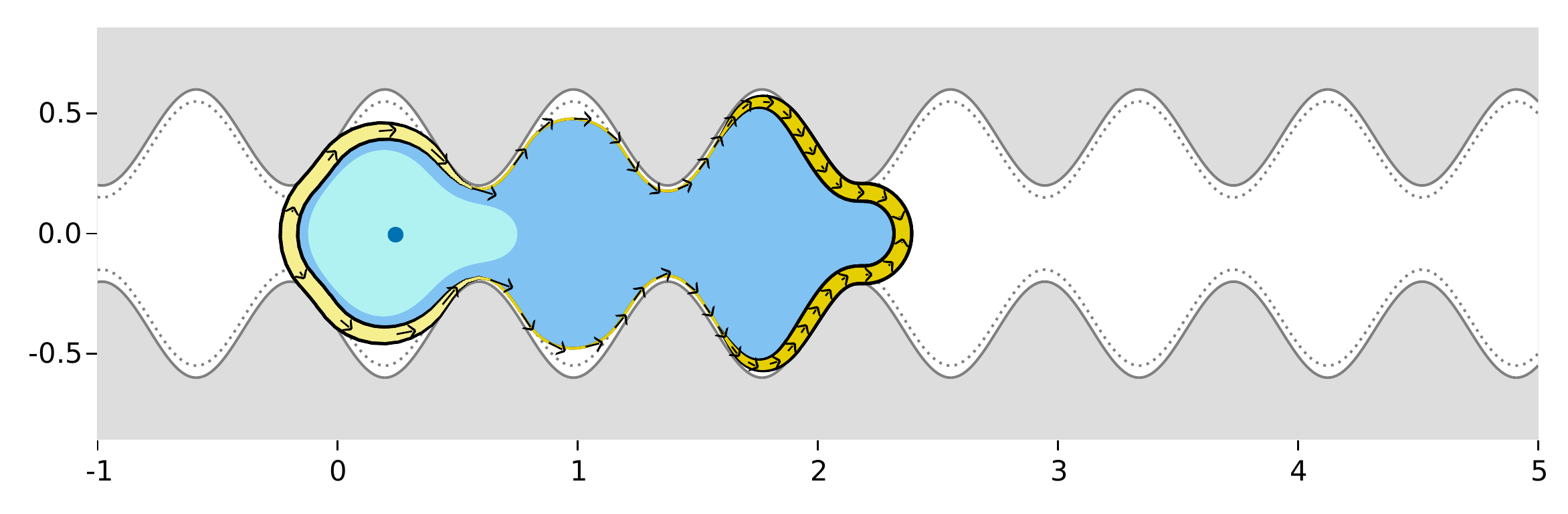}
        &\includegraphics[width=0.22\textwidth]{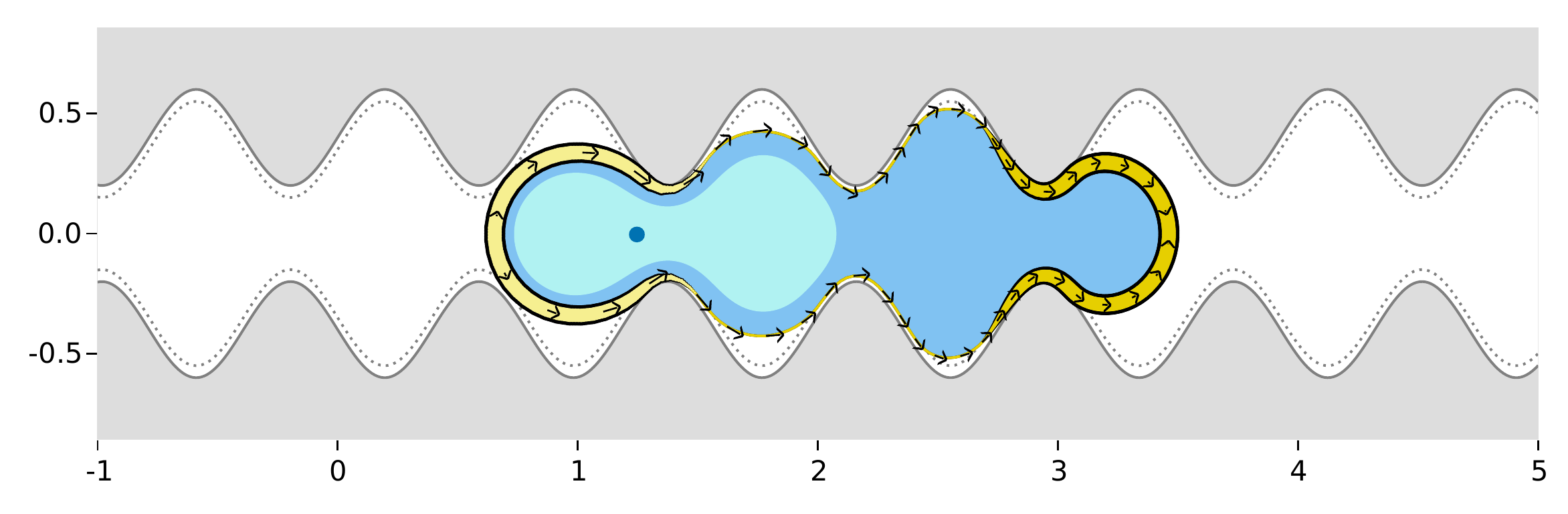}
        &\includegraphics[width=0.22\textwidth]{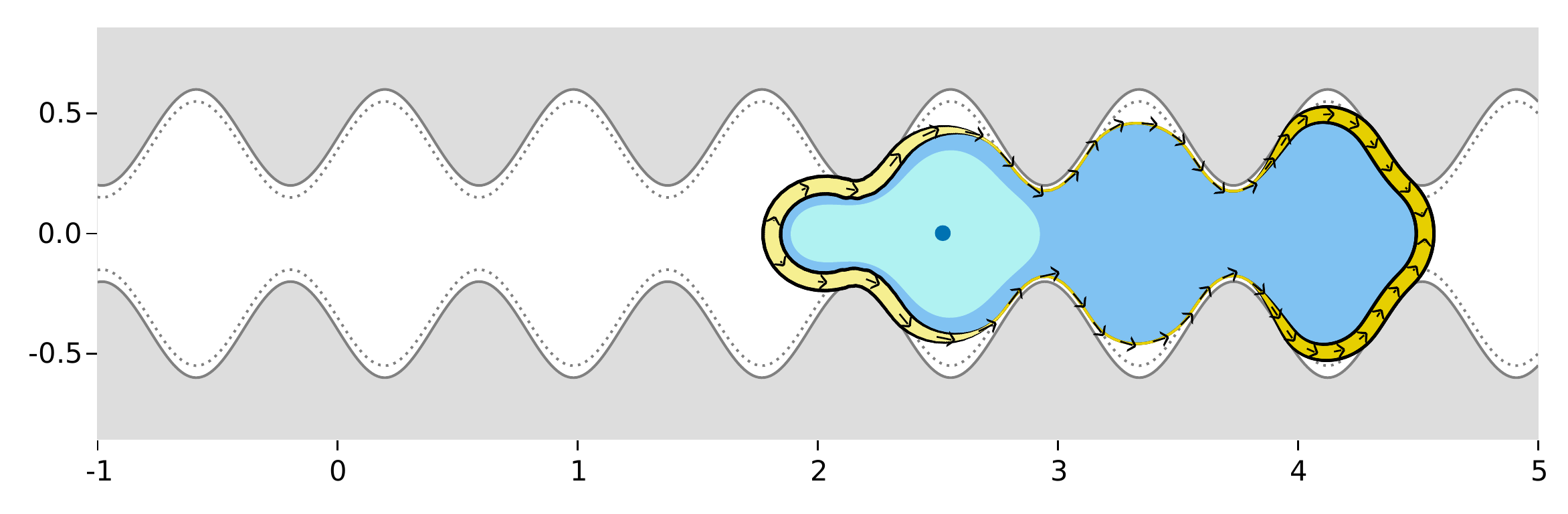}
        &\includegraphics[width=0.22\textwidth]{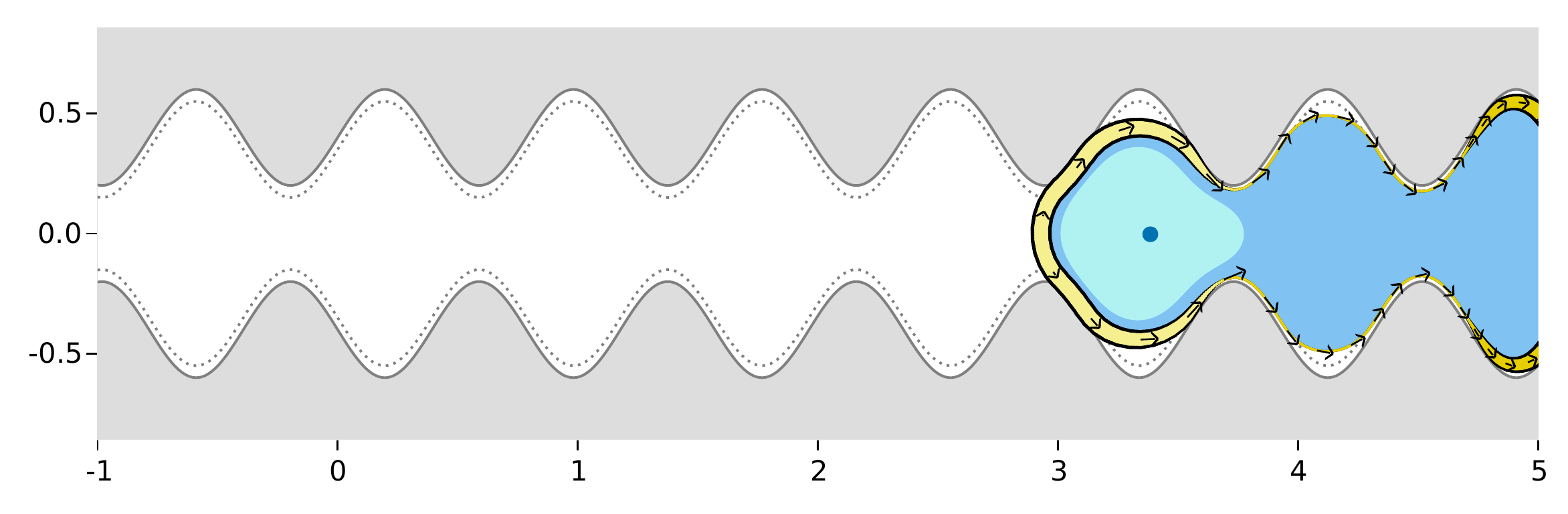}
\\
     \raisebox{1cm}{(c')}
        &\includegraphics[width=0.22\textwidth]{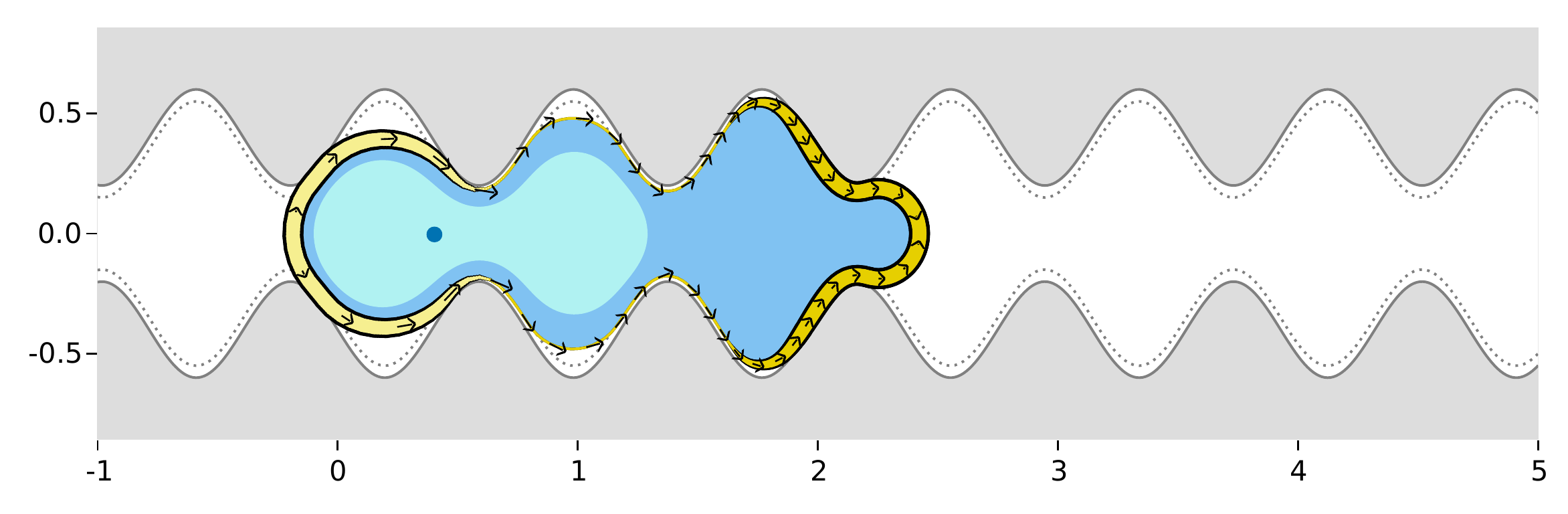}
        &\includegraphics[width=0.22\textwidth]{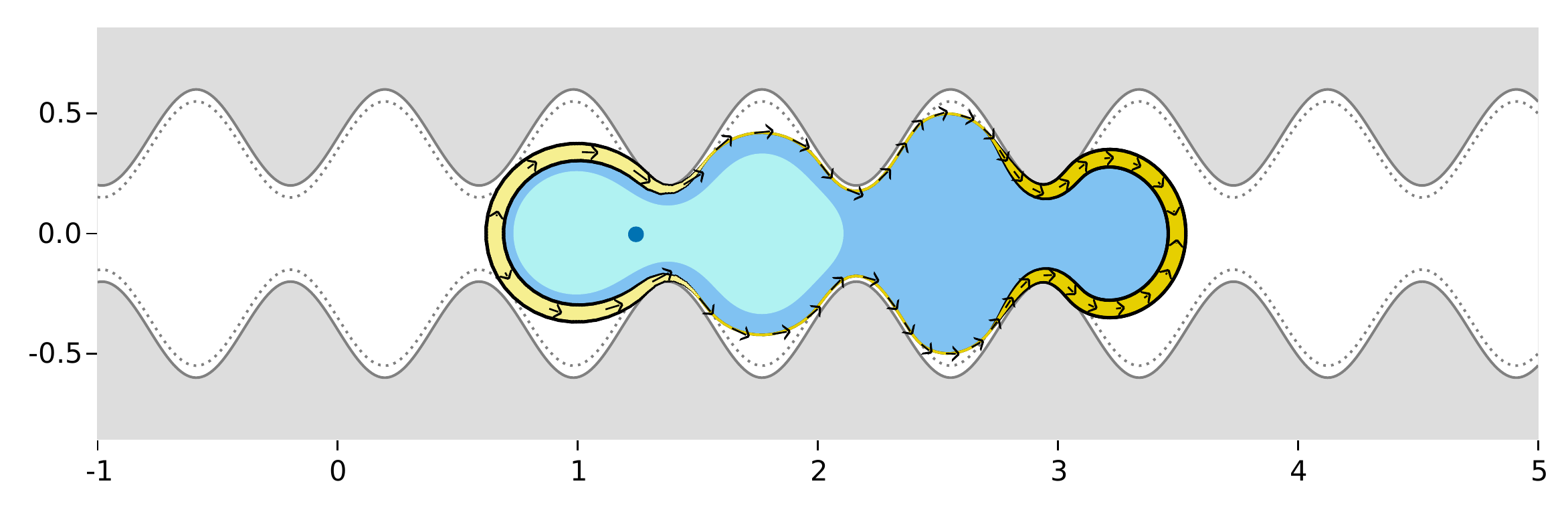}
        &\includegraphics[width=0.22\textwidth]{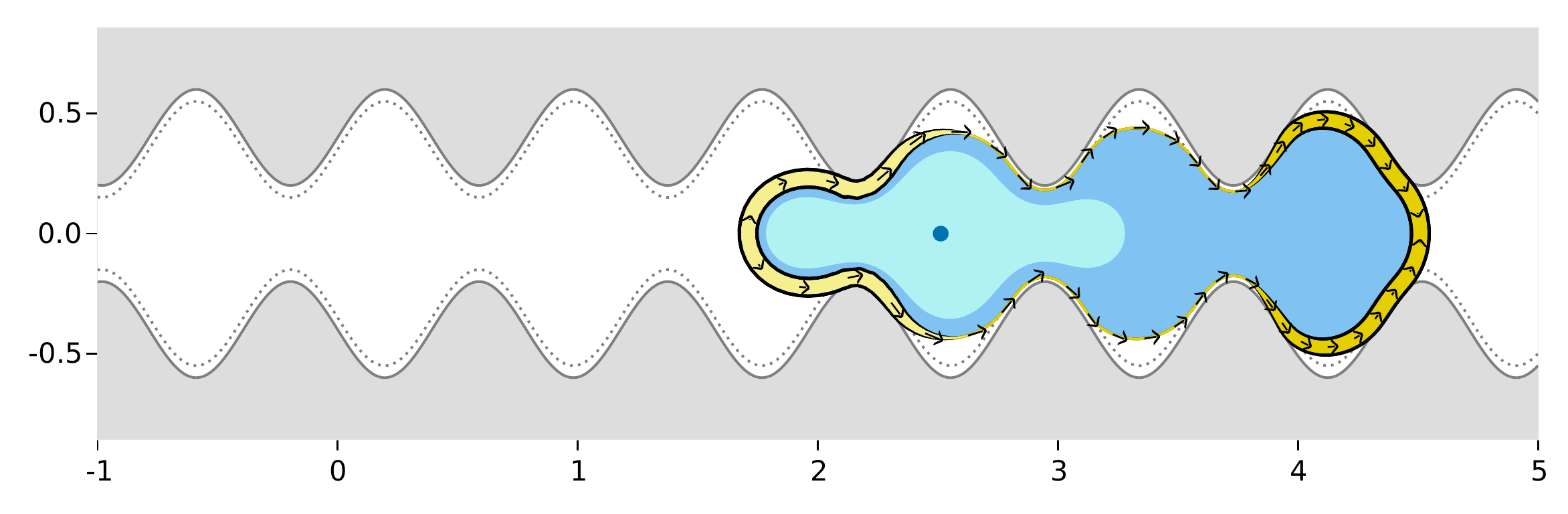}
        &\includegraphics[width=0.22\textwidth]{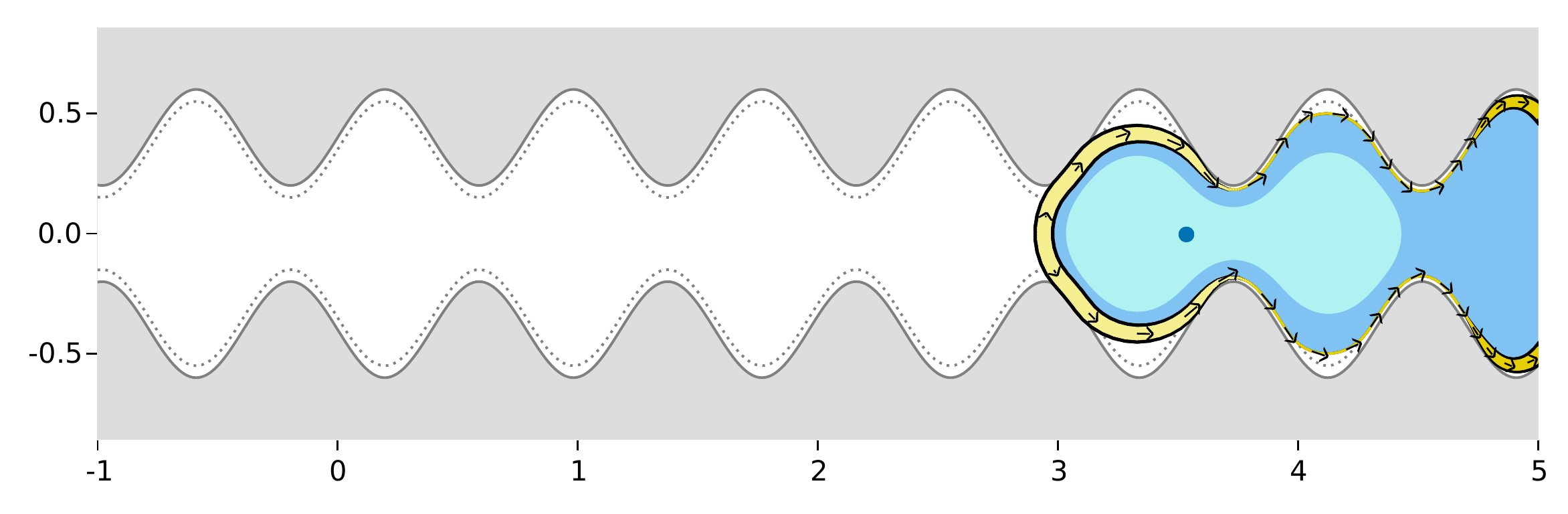}
\\
     \raisebox{1cm}{(d)}
        &\includegraphics[width=0.22\textwidth]{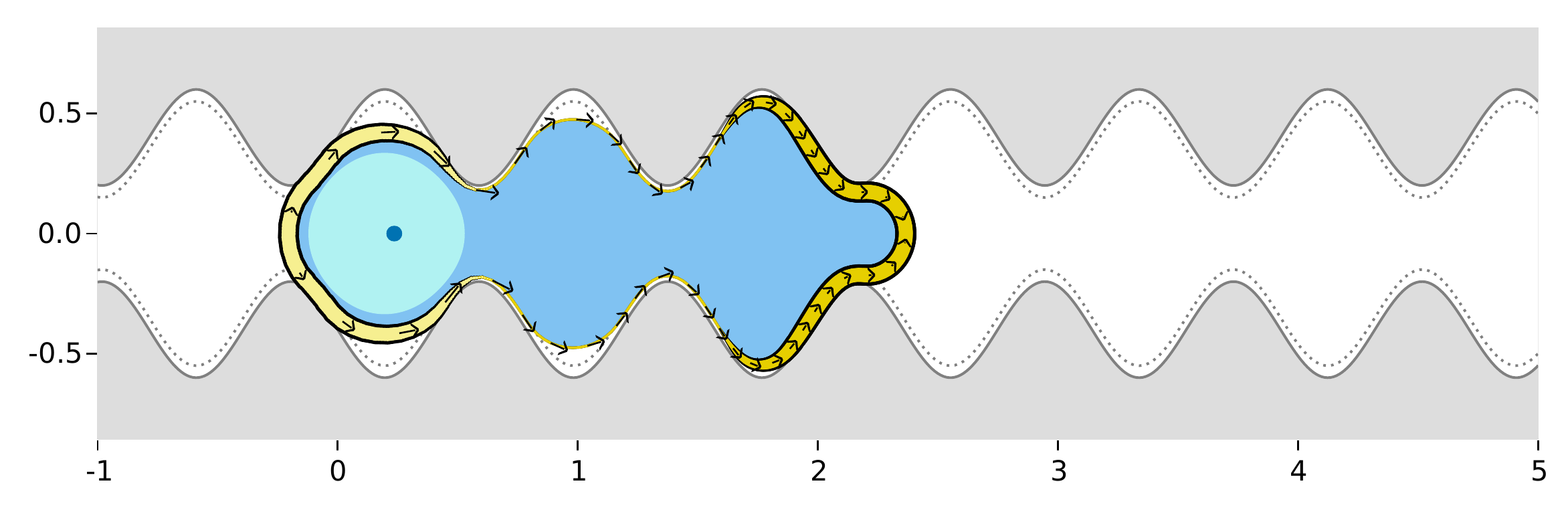}
        &\includegraphics[width=0.22\textwidth]{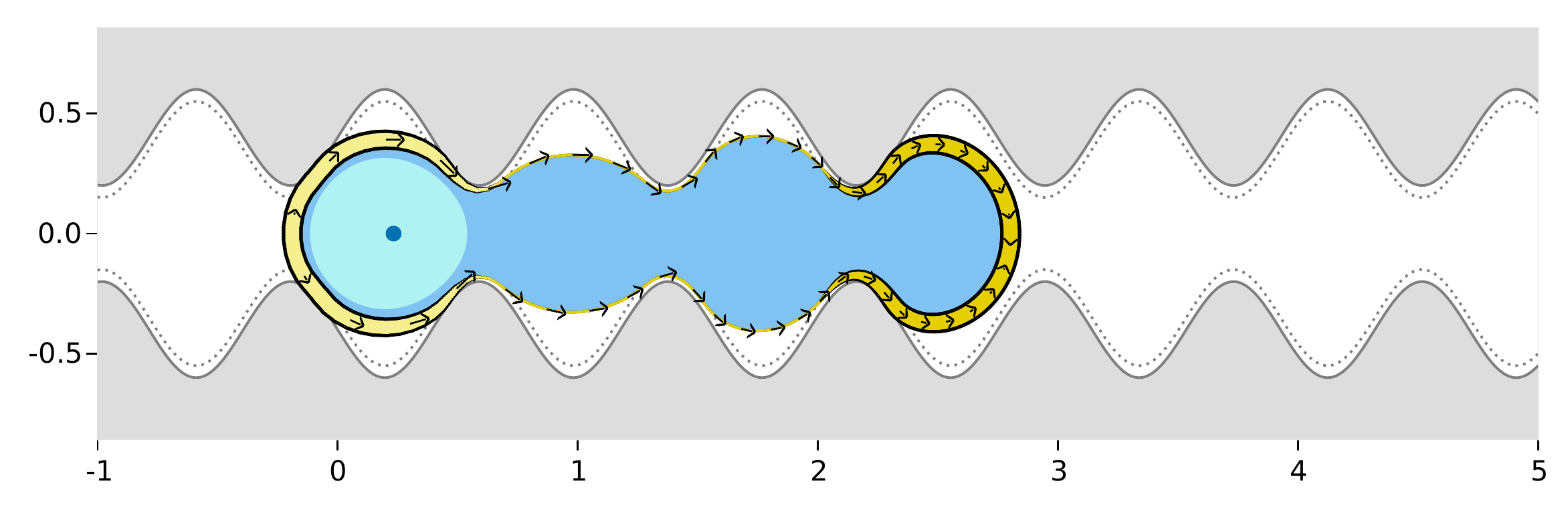}
        &\includegraphics[width=0.22\textwidth]{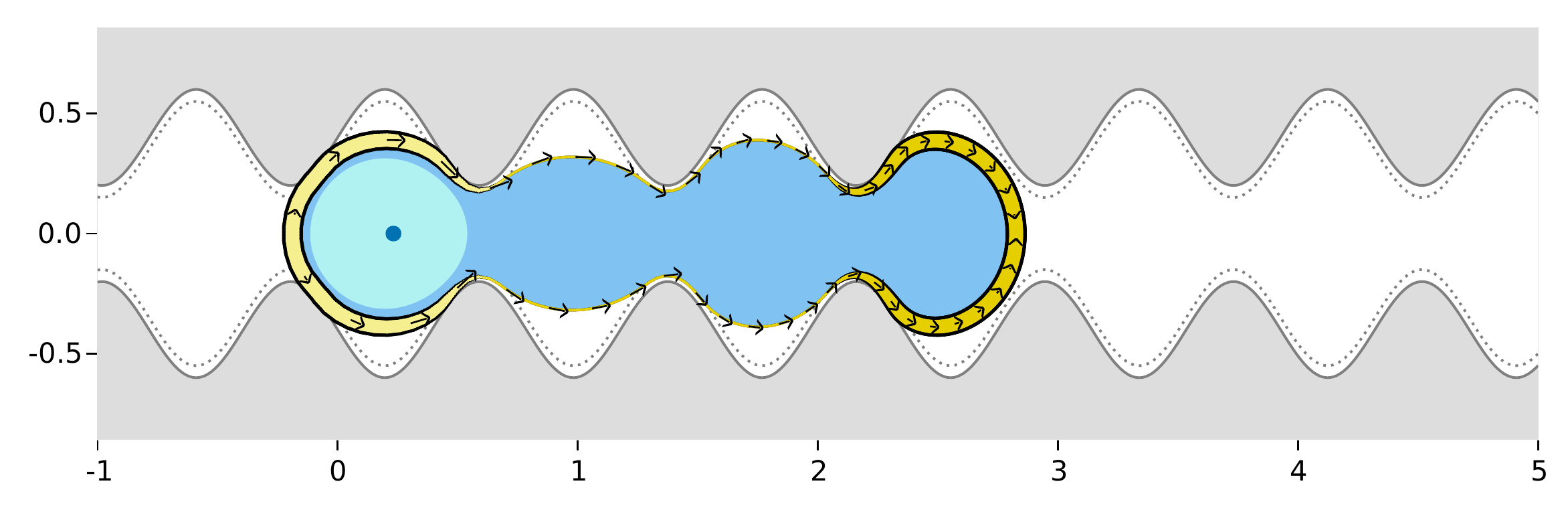}
        &\includegraphics[width=0.22\textwidth]{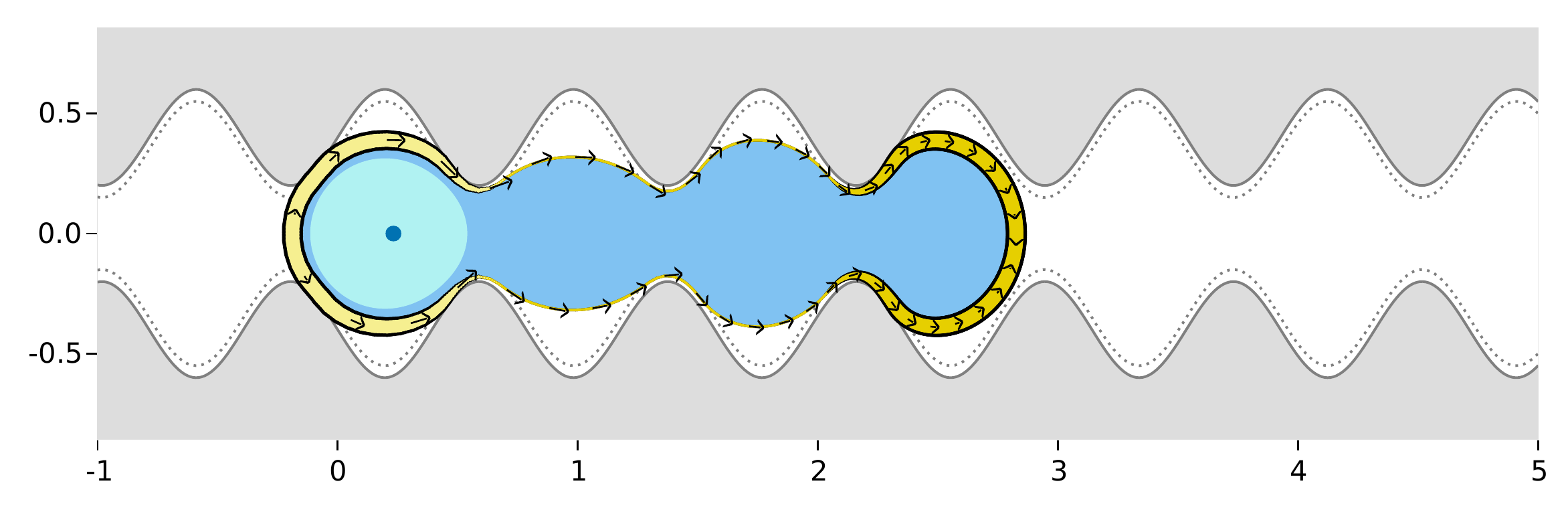}
\\
        & $t = 0.1$ & $t = 0.5$ & $t = 0.9$ & $t = 1.3$
 \end{tabular}
{\footnotesize \caption{Snapshots of the cortex/nucleus/centrosome system at different times, without nucleus (a) and with increasing nucleus stiffness $k_b$.
        The average velocity decreases in the presence of the nucleus, and decreases further are $k_b$ increases ($k_b = 10^{-2.5}$ (b), $k_b= 10^{-1.5}$ (c)), eventually reaching $0$ ($k_b = 10^{-0.5}$ (d)).
        The row (c') illustrates the situation for a larger value of the nucleus area constraint relaxation parameter $\mu_n = 100$. Other rows correspond to $\mu_n = 50$.
The remaining parameters are presented in Table~\ref{table:numerical parameters}. \label{Fig:motion}}}
\end{sidewaysfigure}

The other dimensionless parameters chosen in the simulation are summarized in Table \ref{table:numerical parameters}.
We observe that, for the particular choice of parameters set in the simulations reported in Fig. \ref{Fig:motion}a to \ref{Fig:motion}c', the cell is able to migrate inside the channel even when the nucleus is explicitly modelled. The cell and the nucleus size and shape change during the motion as a balance of the cell area and nucleus area penalizations, their mechanical properties, microtubules and polymerizing forces, and the contact with the channel wall.
 In particular, it is possible to see that the cell first protrudes the cytoplasm to fill the maximum of sinusoidal spaces ahead the cell nucleus. The maximum cytoplasm extension is controlled by the cell target area, the cell membrane target surface, and by the mechanical parameters of the cell membrane. Since these parameters are kept fixed in the simulations in Fig.~\ref{Fig:motion}, the cytoplasm extension inside the channel is comparable for the different cases. When the cytoplasm completely fills the new space, the nucleus is pulled by the microtutubule structure through the constriction in the channel and it forms a bleb in the front until the nucleus is pushed inside the second constriction in the sinusoidal channel. When the nucleus has a low bending stiffness ${k_b}$ (see Fig. \ref{Fig:motion}b), it acquires an hourglass shape by passing through constrictions in the channel and the formation of nuclear protrusions is evident both at the cell front and at the cell rear. However, for higher values of bending stiffness (see Fig. \ref{Fig:motion}c) the nucleus shrinks to pass through constrictions and the formation of blebs and the development of an hourglass shape are less pronounced. Moreover, if we increase the nucleus area-change constraint, $\mu_n$, by keeping the value of $k_b$ fixed, the nucleus cannot shrink to pass through the constrictions and intense nuclear deformations, with blebs both in the front and in the rear of the nucleus, are observed (see Fig. \ref{Fig:motion}c'). In all these cases, when the nucleus fills the new sinusoidal space ahead the cell, the cytoplasm protrudes in the next available sinusoidal space and the process is repeated cyclically, allowing the cell to move forward inside the channel.

 We remark that, although the cell motion occurs inside a simplified geometry, the nucleus hourglass deformation and the formation of blebs in correspondence of small openings in the extracellular space is a characteristic observed during cell motion inside intricate ECM \cite{Wolf_2007, wolf2013physical, Beadle}.

 Furthermore, the importance of including the description of the nucleus becomes clear from simulations reported in Fig. \ref{Fig:motion}d, in order to account for those situations in which the cell can remain trapped inside the channel, because the nucleus cannot deform and squeeze through the small opening of the channel, thus preventing the cell motion.
 Indeed, in this latter case even though the cytoplasm protrudes inside the channel and the nucleus is pulled by microtubules, the energy required to deform and bend the nucleus is too high (due to the high value of $k_b$), so that the nucleus gets stuck in the rear of the cell (see Fig. \ref{Fig:motion}d). This finding is in qualitative agreement with a number of experimental works, such as \cite{Wolf_2007, wolf2013physical, Rolli, Beadle}, where the cell migratory capability is associated with nuclear deformations, and the existence of a critical ECM gap size below which cell migration is entirely hampered in the absence during non-proteolytic has been observed. Such a critical size was termed ``the physical limit of cell migration'' ~\cite{wolf2013physical}.

Moreover, these results are in agreement with other mathematical models, dealing with the influence of the mechanical properties of the nucleus onto cell migratory ability into channels composed by extracellular matrix \cite{giverso2014influence, giverso2018nucleus, scianna2013cellular, scianna2013modeling, scianna2014cellular} and through a dense network of static cells \cite{Lee2017}.
In particular, with respect to previous mechanical models \cite{giverso2014influence, giverso2018nucleus}, this work provides better insights into the phenomenon, since the whole dynamics of the process is investigated, along with the influence of the cell membrane. Furthermore, the motion of the cell in this case does not require the presence of an external flux as in \cite{Lee2017}, but only relies on cell deformation, cortex polymerization and the microtubule activity.
Finally, compared with previous models derived using an extended version of a Cellular Potts Model \cite{scianna2013cellular, scianna2013modeling, scianna2014cellular}, this work allows, in principle, to obtain quantitative results of the whole migratory process, by including in the model identifiable mechanical parameters.\\
Of course, the dynamics of cell and nucleus deformation and motion that can be reproduced by our model is wider than the one captured by these benchmark simulations, therefore we will investigate in Section \ref{res2} the dependence of the mean cell speed and the nucleus shape on the parameters of the model.

\subsection{Influence of the model parameters} \label{res2}

\begin{figure}[H]
    \captionsetup{position=top}
    \centering
    \begin{subfigure}[c]{0.04\textwidth}
        \subcaption{}
        \label{Fig:v_kb}
    \end{subfigure}
    \begin{minipage}[c]{0.94\textwidth}
        \centering
        \includegraphics[width=\textwidth]{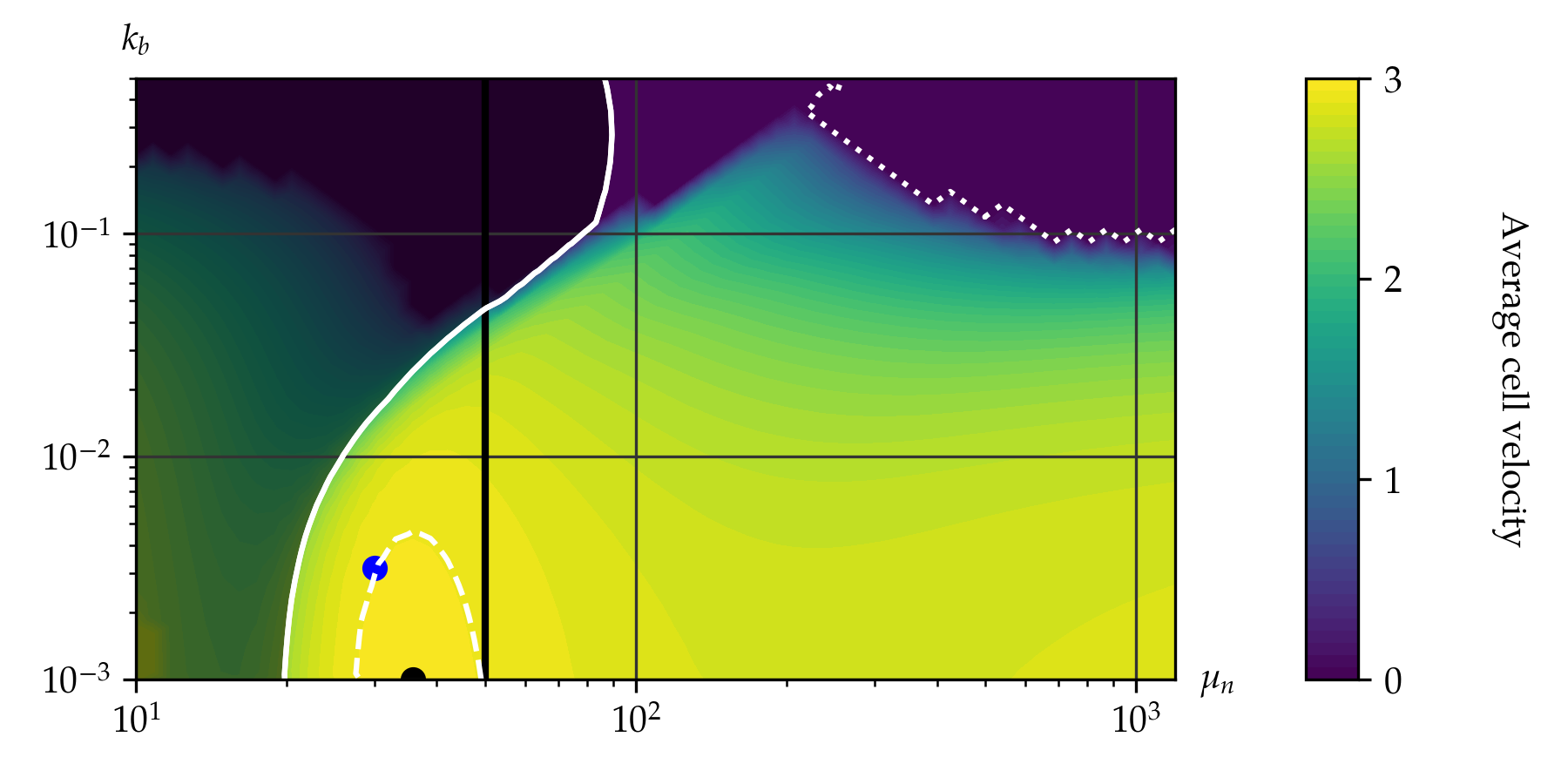}
    \end{minipage}
    \\
    \begin{subfigure}[c]{0.04\textwidth}
        \subcaption{}
        \label{Fig:motion_min_area}
    \end{subfigure}
    \begin{minipage}[c]{0.94\textwidth}
        \centering
        \includegraphics[width=\textwidth]{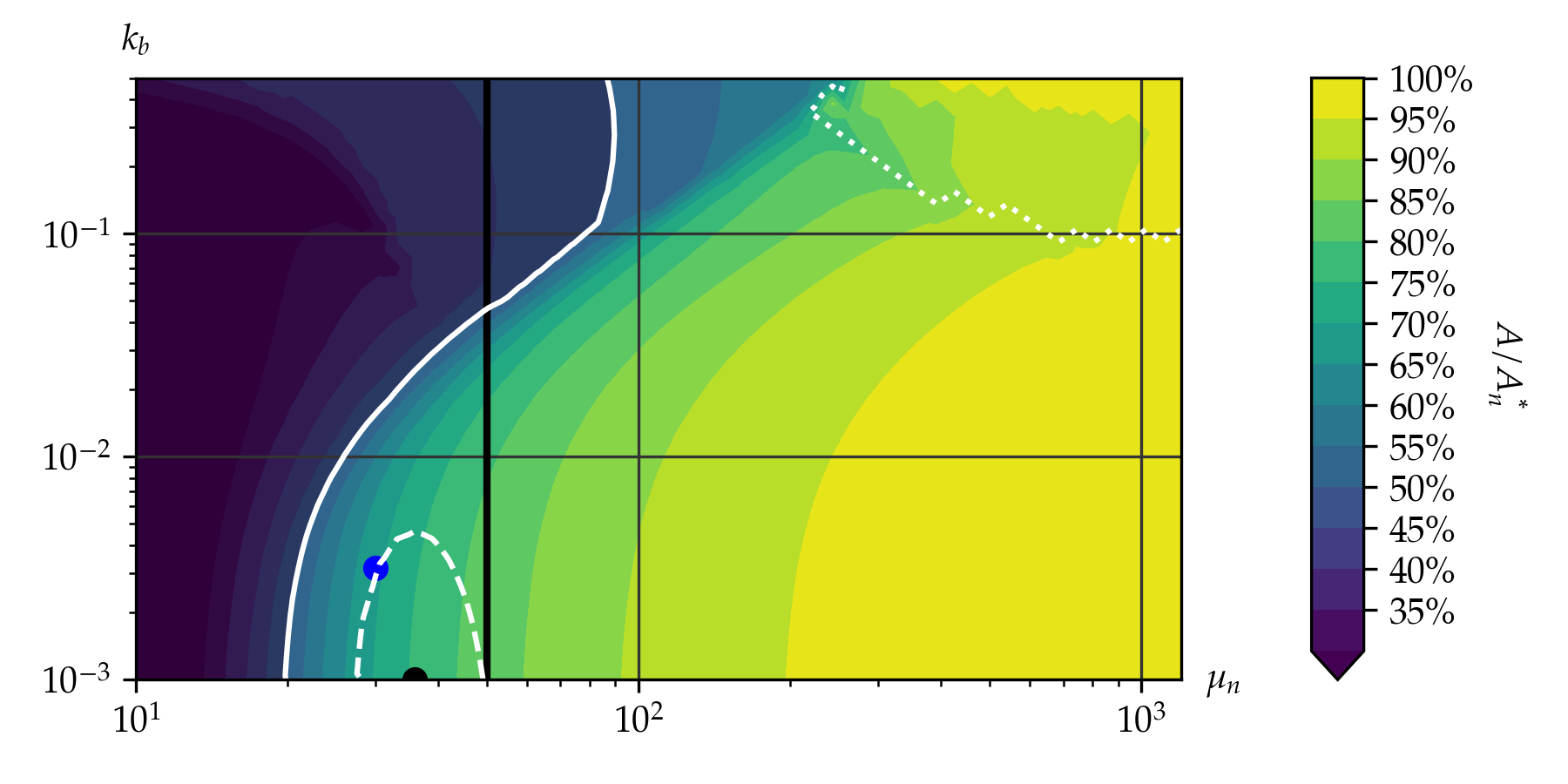}
    \end{minipage}
\caption{
    Average cell speed (a) and ratio between the average and the target nucleus areas (b) while moving inside a sinusoidal channel.
    The shaded region on the left corresponds to a ratio $A/A_n^* < 50\%$.
    The region bounded by the dotted line corresponds to $A/A_n^* > 70\%$ and zero speed.
    The black dot denotes the point of maximum speed, while region for which the velocity is above
    $99\%$ of the maximum velocity is marked with the dashed white line.
The blue dot corresponds to the values of $\mu_n$ and $k_{b}$ taken for Figure \ref{fig:beta_omega_width}.
 The parameters used are given in Table \ref{table:numerical parameters}.
}
    \label{Fig:v_mu_kb}
\end{figure}

In order to understand how the different parameters of the model affect the capability of the cell to migrate inside sinusoidal channels and its mean speed, we perform a set of numerical experiments, changing one parameter at a time.
To compare the results, we defined the mean cell speed as the average speed of the tip of the cell over one period of the cyclic motion, and the mean nucleus area as the average area enclosed by the nuclear membrane over the same period.
The value of the mean cell speed for varying nuclear stiffness is shown in Fig.~\ref{Fig:v_kb}, with the area-change constraint ($\mu_n$) and the bending modulus of the nucleus ($k_b$) along the horizontal and vertical axes, respectively.
From Fig.~\ref{Fig:v_kb}, one sees that the capability of the cell to move inside the channel and its speed highly depend on the nucleus mechanical properties.

Indeed, cell migration is hampered when the energy needed to bend the nuclear envelope is too high (i.e., high values of the parameter $k_b$), since in this case the nucleus is not able to intensively deform and acquire the hourglass shape required to pass the constrictions in the sinusoidal channel. Therefore, for a given value of $\mu_n$, the speed of the cell decreases for increasing values of the bending modulus $k_b$, until the cell is stuck inside the channel (zero speed). The threshold for $k_b$  allowing cell motion depends on the value of $\mu_n$.

Specifically, the constraint imposed by the bending stiffness of the nucleus is more restrictive when the chromatin inside the nucleus is little compactible, i.e. for high values of the parameter $\mu_n$. Indeed, when $\mu_n$ is sufficiently high the nuclear deformations occurs maintaining the nuclear area close to the target one $A_n^*$. This is evident looking at the plot in Fig. \ref{Fig:motion_min_area}, where we report the average nucleus area (still computed over one period of the cyclic motion) with respect to the nucleus target area, for the same values of the nucleus mechanical parameters used in Fig.~\ref{Fig:v_kb}: for high values of the parameter $\mu_n$ the nucleus moves preserving in average more than the $90\%$ of its target area.

The plot in Fig. \ref{Fig:motion_min_area} also allows to comment on the admissibility of the cell velocities predicted by the model.
Indeed, from Fig. \ref{Fig:motion_min_area} it is clear that for small values of the parameters $\mu_n$ the nucleus can decrease its area (which would correspond to the volume in a three-dimensional simulation) below a physiological threshold (shaded region on the left of Fig. \ref{Fig:motion_min_area}).
Even though there is biological evidence \cite{friedl_nuclear_2011, Rowat_2006, Versaevel2012}
that the volume of the nucleus is not preserved during large elongations ---which suggests that the nuclear envelope is permeable to aqueous material and that the chromatin structure can compact itself (chromatin condensation)--- the nucleus volume cannot shrink under a minimum threshold.
In particular, in \cite{Rowat_2006}, it was shown that although nuclei experienced a marked loss of total volume under aspiration, it stabilized above $30-40\%$ of the initial nuclear volume.
Therefore, mechanical parameters allowing the nuclear area going below the $30\%-40\%$ of the target area should be disregarded (see the white isoline in Fig.~\ref{Fig:v_kb}, corresponding to the threshold $A/A_n^* = 50\%$).
We observe that, relaxing the constraint on the area change (i.e., lowering the value of $\mu_n$), so that the nucleus can shrink, the cell can move inside the channel even for high values of bending modulus $k_b$, albeit with lower velocity (green region in  the top of Fig.~\ref{Fig:v_kb}). \\
Furthermore, the proposed model predicts the existence of an optimal region in the space of the mechanical parameters $k_b$-$\mu_n$ for which the cell speed is maximal (dashed region in  Fig.~\ref{Fig:v_kb}) and the nuclear area is in the range $70\%-85\%$ of the target area.

\begin{table}
  \begin{threeparttable}[b]
    \caption{Values and ranges of the dimensionless parameters used for the numerical experiments}
    \label{table:numerical parameters}
    \begin{tabular}{lccc}
        & Symbol & Value & Range
        \\
        \toprule
        \textbf{Cell cortex related parameters} &&&
        \\
        \midrule
        Cell/environment pressure differential & $\Delta p_c$ & $2.56$ &
        \\
        Membrane elasticity & $k_c$ & $0.3$ & 
        \\
        Cell target area & $A_c^*$ & $1.8$ &
        \\
        Cell area constraint relaxation constant & $\mu_c$ & $50$ &
        \\
        Cortex polymerization rate & $r_\mathrm{pol}$ & $10$ &
        \\
        \midrule
        \textbf{Nucleus related parameters} &&& 
        \\
        \midrule
        Nucleus/cell pressure differential & $\Delta p_n$ & 1 & \\
        Nucleus target area & $A_n^*$ & $0.7$ & \\
        Nucleus area constraint relaxation constant & $\mu_n$ & 30 & $[10,70]$\tnote{1} \\
        Nucleus bending stiffness & $k_{b}$ & $3.16 \times 10^{-3}$ & $[5\times 10^{-4}, 5\times 10^{-1}]$\tnote{1} \\
                                  & & & $\{10^{-2.5}, 10^{-1.5}, 10^{-0.5} \}$\tnote{2} \\
        Nucleus/Cortex interaction charact. length & $\xi_\mathrm{cont}$ & 10 & \\
        Nucleus/Cortex interaction coefficient & $k_\mathrm{cont}$ & 5 & \\
        \midrule
        \textbf{Centrosome related parameters} &&&
        \\
        \midrule
        Microtubules friction coefficient & $k_\tau$ & $10^{-4}$ &
        \\
        Centrosome link stiffness & $k_{e}$ & $10^{-3}$ & \\
        \midrule
        \textbf{Channel geometry related parameters} &&&
        \\
        \midrule
        Sharpness\tnote{*} & $\xi$ & 20 & \\
        Depth & $f_\beta$ & $0.2$ & $[0.03, 0.35]$\tnote{3} \\ 
        Pulsation & $f_{\omega_0}$ & $8$ & $[5, 11.7]$\tnote{3} \\
        Mean width & $f_\mathrm{width}$ & $0.4$ & $[0.27, 0.8]$\tnote{3} \\ 
        \midrule
        \textbf{Numerical parameters} &&&
        \\
        \midrule
        Size of the cortex discretization & $N_c$ & $250$ &
        \\
        Initial time step & $\Delta t$ & $2 \times 10^{-4}$ &
        \\
        Size of the nucleus discretization & $N_n$ & 200 &
        \\
        \bottomrule
    \end{tabular}
    \begin{tablenotes}
        \footnotesize
        \item[*]{Inverse of the width of the approximating potential, see Fig.~\ref{Fig:channel}}
        \item[1]{Figures \ref{Fig:v_mu_kb}}
        \item[2]{Figure \ref{Fig:motion}}
        \item[3]{Figure \ref{fig:beta_omega_width}}
    \end{tablenotes}
    \end{threeparttable}
\end{table}

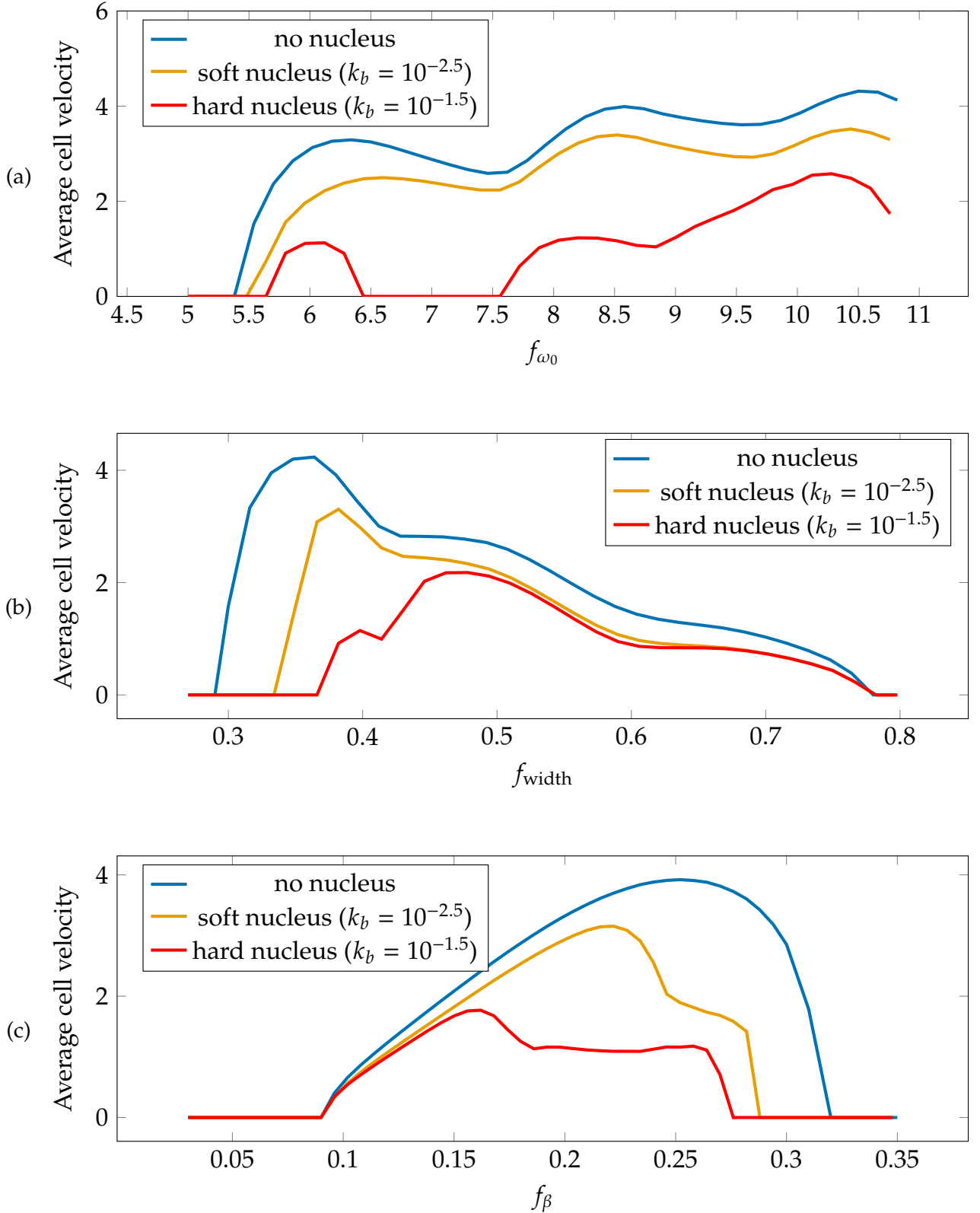
\begin{figure}[H]
    \captionsetup{position=top}
    \centering
    \begin{subfigure}[c]{0.04\textwidth}
        \subcaption{}
        \label{fig:omega}
    \end{subfigure}
    \begin{minipage}[c]{0.94\textwidth}
        \centering

        \ifthenelse{\boolean{usetikz}}{
            \tikzsetfigurename{tikz_owb_}
            \begin{tikzpicture}
                \pgfplotsset{
                    width=0.99\textwidth,
                    height=0.4\textwidth,
                }

                \pgfplotstableread[col sep=comma]{data/bow/no_nuc_omega_3.txt} \omeganonuc;
                \pgfplotstableread[col sep=comma]{data/bow/nuc_omega_1.txt} \omeganuc;
                \pgfplotstableread[col sep=comma]{data/bow/nuc_omega_1_hard.txt} \omeganuchard;
                \begin{axis}[
                    ymin=0, ymax=6,
                    xlabel=$f_{\omega_0}$,
                    ylabel=Average cell velocity,
                    legend pos=north west
                    ]
                    \addplot[color=custom_blue, line width=\plotlinewidth] table[x index=0, y index=1] \omeganonuc;
                    \addlegendentry{no nucleus}
                    \addplot[color=custom_yellow, line width=\plotlinewidth] table[x index=0, y index=1] \omeganuc;
                    \addlegendentry{soft nucleus ($k_b = 10^{-2.5}$)}
                    \addplot[color=red, line width=\plotlinewidth] table[x index=0, y index=1] \omeganuchard;
                    \addlegendentry{hard nucleus ($k_b = 10^{-1.5}$)}
                \end{axis}
            \end{tikzpicture}
        }{\includegraphics[width=\textwidth]{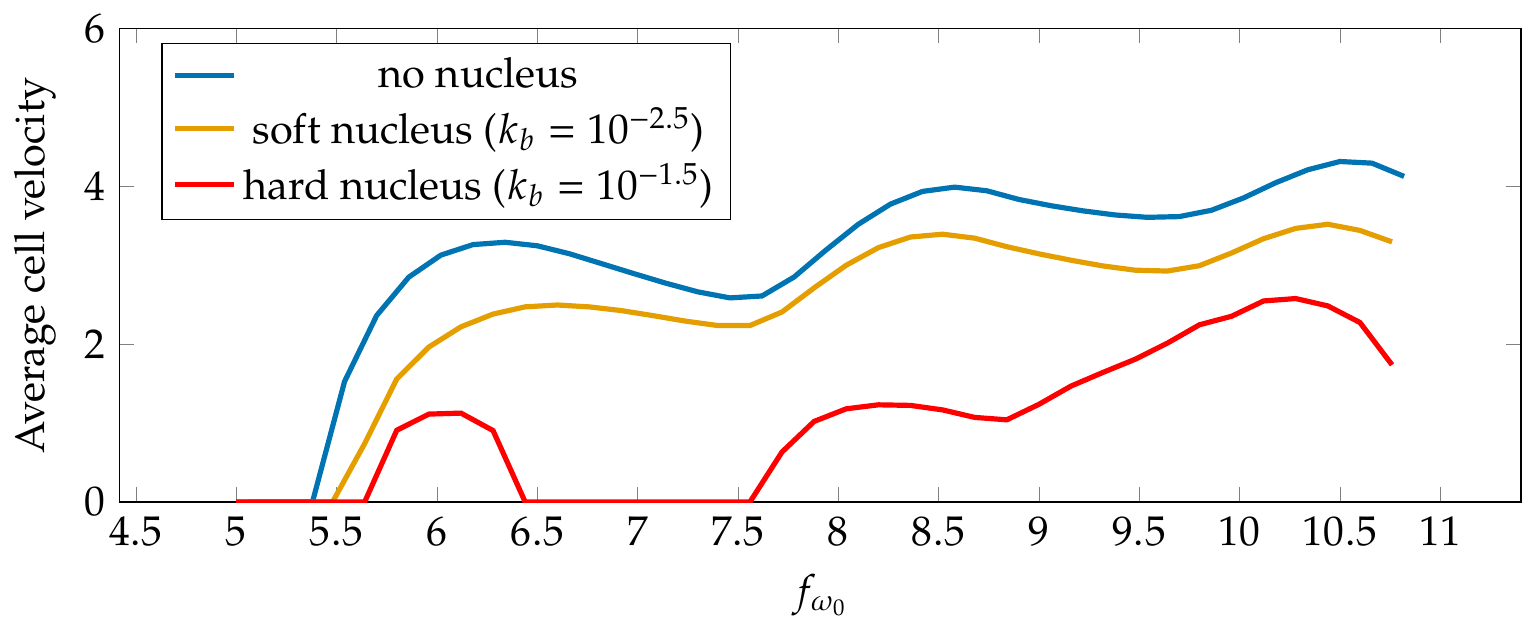}}
    \end{minipage}

    \vspace{1cm}

    \begin{subfigure}[c]{0.04\textwidth}
        \subcaption{}
        \label{fig:width}
    \end{subfigure}
    \begin{minipage}[c]{0.94\textwidth}
        \centering
        \ifthenelse{\boolean{usetikz}}{
            \tikzsetfigurename{tikz_owb_}
    \begin{tikzpicture}
        \pgfplotsset{
            width=0.99\textwidth,
            height=0.4\textwidth,
            xtick = {0.3,0.4,0.5,0.6,0.7,0.8},
            xticklabels = {0.3,0.4,0.5,0.6,0.7,0.8},
        }
        \pgfplotstableread[col sep=comma]{data/bow/no_nuc_width_3.txt} \widthnonuc;
        \pgfplotstableread[col sep=comma]{data/bow/nuc_width_1.txt} \widthnuc;
        \pgfplotstableread[col sep=comma]{data/bow/nuc_width_1_hard.txt} \widthnuchard;

        \begin{axis}[
            xlabel=$f_\textrm{width}$,
            ylabel=Average cell velocity
            ]
            \addplot[color=custom_blue, line width=\plotlinewidth] table[x index=0, y index=1] \widthnonuc;
            \addlegendentry{no nucleus}
            \addplot[color=custom_yellow, line width=\plotlinewidth] table[x index=0, y index=1] \widthnuc;
            \addlegendentry{soft nucleus ($k_b = 10^{-2.5}$)}
            \addplot[color=red, line width=\plotlinewidth] table[x index=0, y index=1] \widthnuchard;
            \addlegendentry{hard nucleus ($k_b = 10^{-1.5}$)}
        \end{axis}

    \end{tikzpicture}
    }{\includegraphics[width=\textwidth]{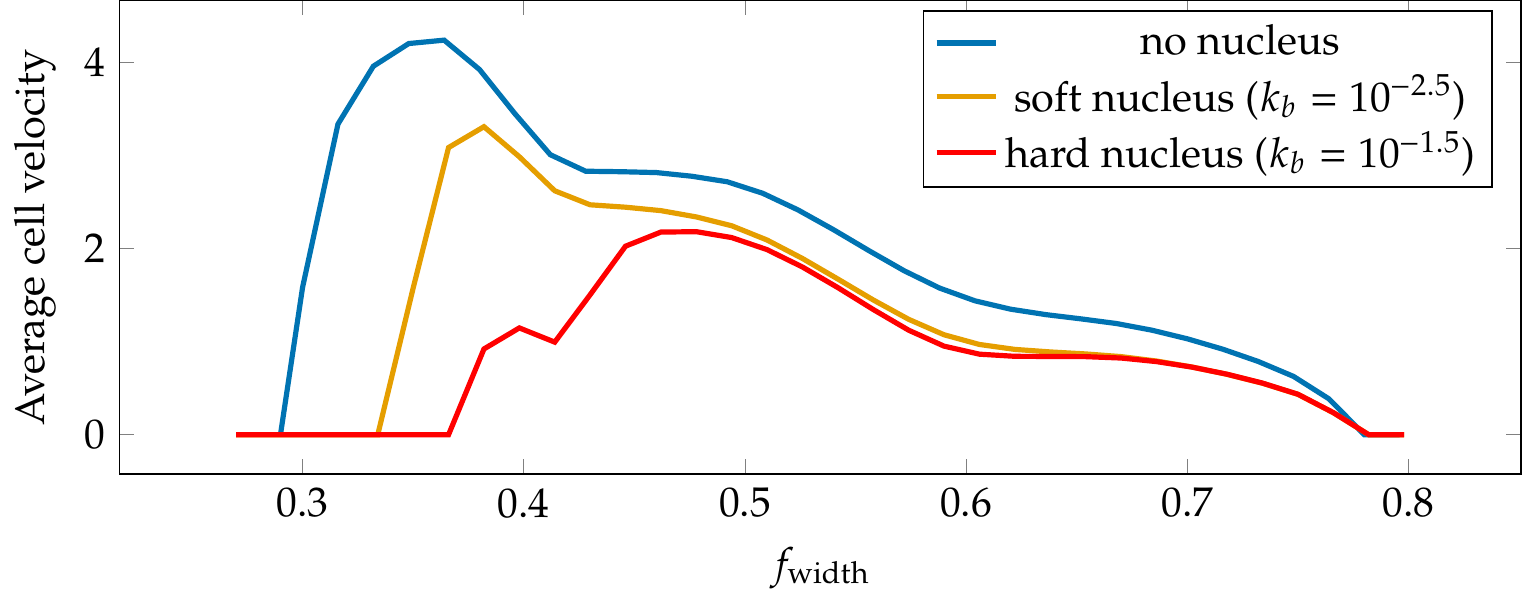}}
    \end{minipage}

    \vspace{1cm}

    \begin{subfigure}[c]{0.04\textwidth}
        \subcaption{}
        \label{fig:beta}
    \end{subfigure}
    \begin{minipage}[c]{0.94\textwidth}
        \centering
        \ifthenelse{\boolean{usetikz}}{
            \tikzsetfigurename{tikz_owb_}
    \begin{tikzpicture}
        \pgfplotsset{
            width=0.99\textwidth,
            height=0.4\textwidth,
            xtick = {0.05,0.1,0.15,0.2,0.25,0.3,0.35,0.4},
            xticklabels = {0.05,0.1,0.15,0.2,0.25,0.3,0.35,0.4},
        }

        \pgfplotstableread[col sep=comma]{data/bow/no_nuc_beta_3.txt} \betanonuc;
        \pgfplotstableread[col sep=comma]{data/bow/nuc_beta_1.txt} \betanuc;
        \pgfplotstableread[col sep=comma]{data/bow/nuc_beta_1_hard.txt} \betanuchard;

        \begin{axis}[
            xlabel=$f_\beta$,
            ylabel=Average cell velocity,
            legend pos=north west
            ]
            \addplot[color=custom_blue, line width=\plotlinewidth] table[x index=0, y index=1] \betanonuc;
            \addlegendentry{no nucleus}
            \addplot[color=custom_yellow, line width=\plotlinewidth] table[x index=0, y index=1] \betanuc;
            \addlegendentry{soft nucleus ($k_b = 10^{-2.5}$)}
            \addplot[color=red, line width=\plotlinewidth] table[x index=0, y index=1] \betanuchard;
            \addlegendentry{hard nucleus ($k_b = 10^{-1.5}$)}
        \end{axis}
    \end{tikzpicture}
    }{\includegraphics[width=\textwidth]{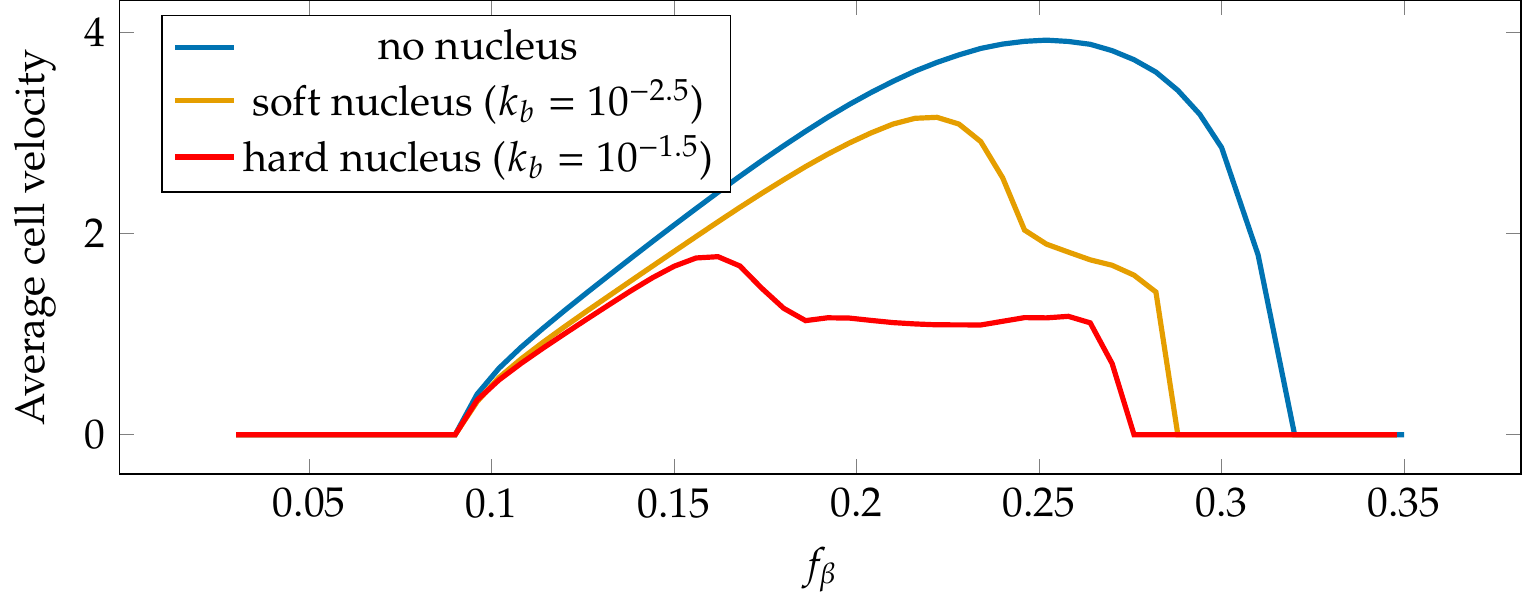}}
    \end{minipage}

    \caption{Comparison of the motilities associated with cell with and without nucleus, for a range of channel parameters.
    Top: Wave number of the channel pattern, middle: depth of the pattern, bottom: mean width of the channel.
    Model parameters are summed up in Table~\ref{table:numerical parameters}.}
    \label{fig:beta_omega_width}
\end{figure}

Finally, in Fig.~\ref{fig:beta_omega_width}, we consider the cell average velocity for varying values of the parameters describing the channel geometry and we compare the results obtained either with or without the description of the nucleus.
In particular, we first consider the pulsation of the channel in determining the cell capability to move.
As observed for the cell without the nucleus (blue line in Fig.~\ref{fig:omega}) the cell is able to migrate only if the channel is sufficiently structured, i.e. for channel pulsation above a minimum threshold, $f_{\omega_{min}}$.

The value of the threshold $f_{\omega_{min}}$ is only slightly influenced by the presence of the nucleus, that mainly affects the cell average speed, which is slowed down by the presence of the nucleus. Specifically, the cell speed decreases when increasing the value of $k_b$ (see the yellow and red curves).
However, for high values of $k_b$ the cell gets stuck (see red curve in Fig.~\ref{fig:omega}) even for intermediate values of channel pulsation,
since the hourglass nuclear deformation is impeded by the high bending modulus and the space between two subsequent constriction is enough to host the cell nucleus. Increasing further the channel pulsation, the nucleus characterized by an high bending modulus is again able to migrate since the space between two constriction becomes restrictive for the preservation of the nuclear area and the nucleus can deform acquiring an ellipsoidal shape, without following the structure of the channel walls.
We remark that, since the discretization of the cell membrane is finite, the channel structure will not be well resolved as $f_ {\omega_0}$ is taken larger and larger. This explains why the first plot in Fig.~\ref{fig:omega} is cut on the right hand side. Note that this will occur for any discretization size.

The influence of the channel width on the cell velocity is analysed Fig.~\ref{fig:width}: without considering the nucleus (blue curve), the motion of the cell is allowed for smaller channel width and this ``physical limit of cell migration'' is related to the capability of the cell envelope and the cytosol to deform and enter even small gaps. Indeed, the size of the smallest constriction inside the channel is $2(f_\mathrm{width} - f_{\beta})$. When we include the nucleus, the cell is no longer able to move inside channels with a too narrow necking (i.e., small values of $f_\mathrm{width}$), since the nucleus hinders cell capability to deform (yellow curve). Rising the bending modulus (red curve), the cell can move inside the channel for large values of $f_\mathrm{width}$ only. On the other hand, the upper limit is the same in all three cases, since it is related to the capability of the whole cell to maintain the contact with the channel wall and it is not affected by the presence and the mechanical properties of the nucleus.
Furthermore, the plots in Fig.~\ref{fig:width} show the well-known bimodal behaviour of cell velocity for varying channel width: cells cannot migrate both inside very small channels (since deformations would be too large) and inside very large channels (since the cell does not touch the boundary of the channel.
The velocity of cell motion is also in this case slowed down by the presence of the nucleus. Specifically the speed decreases for increasing values of $k_b$, as previously observed.

The bimodal behaviour in the cell velocity can be observed also changing the value of $f_{\beta}$ and keeping the channel width fixed (see Fig.~\ref{fig:beta}). In this case, the upper limit of $f_{\beta}$ is related to size of the small constriction that the cell could enter in order to move inside the channel and thus it is highly influenced by the presence of the nucleus and by its mechanical properties.
On the other hand, when $f_\beta$ is small the channel is almost flat and the cell cannot migrate, independently on the presence of the nucleus and its mechanical properties. Therefore, the lower limit for $f_{\beta}$ is the same both for the cell with and without the nucleus, since it is related to the mechanism of motion, which requires a sufficient structure of the lateral walls.

\section{Conclusions}
In this work we have proposed a continuous mechanical model for describing single cell adhesion-independent migration inside restrictive 3D environments.
The cell migration is driven by a local imbalance of polymerization and depolymerization of the actin network underneath the cell membrane, which induces a cortex flow. The cell shape is determined by a balance of cytoplasmic pressure, elastic behaviour of the cell cortex and interactions with the subcellular elements and the channel walls.
Differently from previous works \cite{Jankowiak_M3AS}, the influence of the nucleus in the process of cell migration is explicitly taken into account by including an inner nuclear membrane, connected to the cell membrane by the microtubule structure, responsible of nucleus location inside the cell. The nucleus shows a resistance to bending, to stretching and to changes in the area enclosed by the nuclear membrane and it becomes the limiting factor in determining the ability of the cell to migrate across small neckings of the channel.
Therefore, the proposed model represents an advancement with respect to the state of the art, since it provides a purely mechanical description of adhesion-free migration, without requiring external chemical stimuli as done in \cite{Lee2017, Vermolen, Cao}. It also allows testing the influence of nucleus mechanical properties in the determination of the physical limit of cell migration, which has been neglected in previous models \cite{Wu_Othmer, Othmer, kaoui_lateral_2008, Moure2017, Moure_Gomez2018, Jankowiak_M3AS}.

The model equations have been discretized and solved numerically in order to simulate the process in a 2D geometry, corresponding to a section of the 3D channel with lateral structured walls and top and bottom flat walls.
To formulate an efficient and stable (under appropriate time step restrictions) discretization scheme, we adapt to our model the setting proposed in \cite{mikula_computational_2004, benes_comparison_2009}.

The numerical experimente reproduce qualitatively the behaviour observed in the biological experiments \cite{Reversat2020}, where it is pointed out that adhesion-independent migration needs both confinement and sufficiently structured channel walls. This behaviour is purely related to the cell membrane behaviour and it is not influenced by the presence of the nucleus, since the threshold for the channel structure depth (i.e. lowest value of $f_\beta$ that allows cell motion) predicted by our model is the same of the one predicted by the model developed in \cite{Jankowiak_M3AS}.
However, the cell speed inside the channel and the physical limit of cell migration (i.e. the size of the smallest opening in the channel) substantially depends on the presence of the nucleus and on its mechanical properties, as demonstrated by the parametric study showing the dependency on geometric properties of the channel and on the mechanical properties of the nucleus.
In particular, in agreement with biological experiments  \cite{wolf_2003, wolf_physical_2013}, we observe that a little-deformable nucleus is a limiting factor for cell migration inside restrictive environments. Indeed, by keeping both the geometry of the channel and the cell membrane parameters fixed, we can show the transition from a migrating cell to a non-moving cell, by changing the mechanical properties of the nucleus only (i.e., increasing the bending modulus $k_b$ and the elastic area constraint $\mu_n$).

Furthermore, the sensitivity analysis of the mechanical parameters of the nucleus predicts the optimal nuclear deformability, for which the cell reaches its maximal speed, for low values of $k_b$ and intermediate values of $\mu_n$.

Without claiming to provide quantitative numerical measurements, the results presented in this work are a proof of concept for a comprehensive model of single cell adhesion-independent migration taking into account cell and nucleus mechanics.
Nonetheless, these results need to be validated quantitatively by comparing the predicted evolution with the actual spatio-temporal evolution of the moving cell inside a structured channel, by also keeping track of nuclear deformations. Therefore, additionnal work is needed to inverstigate more realistic in vitro and in vivo conditions, in order to
quantitatively validate the model.

From a modelling point of view, this study has to be seen as a first step and several components would benefit from a more detailed treatment.
In particular, the role of the compensating force, that can be attributed to the transport of actin monomers and their polymerization, certainly requires further studies.
For example, one could include an explicit description of the acto-myosin machinery and its influence in cell polymerization.
Another interesting direction would be to take into account the resistance to deformation of the cortex, which can then be considered as visco-elastic and not only elastic.
Finally, the mechanical description of the nucleus is kept rather simple in the present model and a detailed description of both the deformations occurring inside the nucleus and the exchange of liquid between the nucleus and the surrounding cytosol during the whole process of cell migration, is still missing.
Taking all these effects into account could lead to a more comprehensive understanding of the multiple factors involved in determining the physical limit of cell migration during non-proteolytic migration of cells.

\section*{Acknowledgments}
G.J. acknowledges the support of the Austrian Academy of Sciences, project NST 0001 and of the Austrian Science Fund, grant no. F65.
The work of C.S. has been supported by the Austrian Science Fund, grants no. F65 and W1245.
The work of C.G and L.P. was supported by the MUR,  through the grant Dipartimenti di Eccellenza 2018--2022 Project n. E11G18000350001 and through the project PRIN2017 n. 2017KL4EF3 and PRIN2020 n.
2020F3NCPX.

\section*{Statements and Declarations}
The authors declare no conflict of interest.

\section*{Appendix}
\appendix

\section{Nuclear energy variation} \label{appendix:variation}
We here compute the variation $\delta\En$ of the nuclear membrane energy function, in order to obtain the force acting on the nuclear membrane.
We named
\begin{equation} \label{eq:nucleus_energy expanded}
\En^{(1)} = \dfrac{k_b}{2} \sgint K^2 \,\dd \ell \, ;  \quad \En^{(2)} =\lambda \sgint \, \dd \ell \, ;  \quad \En^{(3)} =\Delta p_n \smashoperator{\int_{\Omega_n}} \,\dd A \,  ; \quad \En^{(4)} =\mu_n \left( \smashoperator{\int_{\Omega_n}}\,\dd A -A_{n}^* \right)^2,
\end{equation}
so that $\En= \En^{(1)} + \En^{(2)} + \En^{(3)} + \En^{(4)}$ and $\delta \En= \delta \En^{(1)} + \delta \En^{(2)} + \delta \En^{(3)} +\delta \En^{(4)} $.
In order to compute $\delta \En^{(j)}$, with $j=1,2,3,4$, it is convenient to use the parametrization of the curve with respect to $\sigma$, since $\dd \ell$ will undergo variations as well for displacement of $Y$.

The first energy term with respect to $\sigma$ (assuming $K_0=0$)  reads
\begin{align*}
\En^{(1)} = \dfrac{k_b}{2} \sgint K^2 \,\dd \ell =  \dfrac{k_b}{2} \spint \left( (\partial_{\sigma \sigma} Y)^2 - \left( \dfrac{\partial_{\sigma \sigma} Y \cdot \partial_{\sigma} Y}{|\partial_{\sigma} Y|}\right)^2\right)  |\partial_{\sigma} Y|^{-3}\,\dd \sigma
\end{align*}
and its variation is \cite{kaoui_lateral_2008}
\begin{align*}
\dfrac{\delta \En^{(1)}}{\delta Y} = - k_b \left|\partial_{\sigma} Y \right|  \left(  \dfrac{\partial^2 K}{\partial \ell^2} + \dfrac{1}{2} K^3  \right) N \, .
\end{align*}
The variation of the second energy term can be easily computed
\begin{align*}
\delta \En^{(2)} &= \lambda \sgint  \delta(\dd \ell )=
\lambda \spint  \delta(\left| \partial_{\sigma} Y \right| ) \, \dd \sigma =
\lambda \spint \dfrac{\partial_{\sigma} Y }{\left| \partial_{\sigma} Y \right|} \cdot \delta(\partial_{\sigma} Y ) \, \dd \sigma =
-\lambda \spint \partial_{\sigma} \left( \dfrac{\partial_{\sigma} Y }{\left| \partial_{\sigma} Y \right|} \right) \cdot \delta Y  \, \dd  \sigma = \\
&= -\lambda \spint \partial_{\sigma} T \cdot \delta Y  \, \dd  \sigma = \lambda \sgint K \, N \cdot \delta Y  \left|\partial_{\sigma} Y \right| \, \dd  \sigma \, ,
\end{align*}
and thus
\begin{align*}
\dfrac{\delta \En^{(2)}}{\delta Y}  &= \lambda \left|\partial_{\sigma} Y \right|  K \, N .
\end{align*}
In order to calculate $\delta \En^{(3)} $ and  $\delta \En^{(4)} $, we use the identity
\begin{align*}
\En^{(3)} &= \smashoperator{\int_{\Omega_n}}  \, \dd A= \dfrac{1}{2}  \sgint  Y \cdot N \, \dd \ell=-\dfrac{1}{2}  \sgint  Y \cdot \partial_{\ell} Y ^\perp \, \dd \ell=- \dfrac{1}{2}  \spint Y \cdot \partial_{\sigma} Y ^\perp \, \dd \sigma
\end{align*}
so that
\begin{align*}
\delta \En^{(3)} &=  - \Delta p_n \spint \partial_{\sigma} Y^\perp \cdot \delta Y  \dd \sigma \, ,\\
\delta \En^{(4)} &=  - \mu_n \left( \spint Y \cdot \partial_{\sigma} Y ^\perp \, \dd \sigma  -A_{n}^* \right) \spint \partial_{\sigma} Y^\perp \cdot \delta Y  \dd \sigma
\end{align*}
Therefore, we have
\begin{align*}
\dfrac{\delta \En^{(3)}}{\delta Y}  &=  - \Delta p_n  \partial_{\sigma} Y^\perp = \Delta p_n  \left|\partial_{\sigma} Y \right|  N \, ,\\
\dfrac{\delta \En^{(4)}}{\delta Y}  &=  - \mu_n  \left( \spint Y \cdot \partial_{\sigma} Y ^\perp \, \dd \sigma  -A_{n}^* \right) \partial_{\sigma} Y^\perp = \mu_n  \left( \spint Y \cdot \partial_{\sigma} Y ^\perp \, \dd \sigma  -A_{n}^* \right) \left|\partial_{\sigma} Y \right|  N
\end{align*}
Finally, the nuclear membrane mechanical force reads
\begin{align*}
F_n =  k_b  \left(  \dfrac{\partial^2 K}{\partial \ell^2} + \dfrac{1}{2} K^3  \right) N  - \lambda  K \, N - \Delta p_n N - \mu_n  \left( \spint Y \cdot \partial_{\sigma} Y ^\perp \, \dd \sigma  -A_{n}^* \right)  N
\end{align*}

\section{Discretization and quadrature formulae for the centrosome and the microtubule structure}
\label{sec:discretization and quadrature formulae}

Assume that $X$ is a polygon with nodes $X_i = X(\theta_i)$ such that
\begin{equation*}
    X(\theta_i) = X_c + r_i
    \begin{pmatrix} \cos \theta_i \\ \sin \theta_i \end{pmatrix}, \text{ so that } |X(\theta_i) - X_c| = r_i\,,
\end{equation*}
where $X_c$ is the centrosome. Then the line segment between $X_i$ and $X_{i+1}$ has the following equation in polar coordinates
\begin{equation*}
    r_i(\theta)(a_i \cos\theta + b_i \sin\theta) = 1
\end{equation*}
where
\begin{equation*}
    a_i = \frac{r_{i+1} \sin \theta_{i+1} - r_{i} \sin \theta_{i}}{r_{i} r_{i+1} \sin(\theta_{i+1} - \theta_{i})},
    \text{ and }
    b_i = -\frac{r_{i+1} \cos \theta_{i+1} - r_{i} \cos \theta_{i}}{r_{i} r_{i+1} \sin(\theta_{i+1} - \theta_{i})},
\end{equation*}
Then we have the following identities:
\begin{align*}
    \int_0^{2\pi} |X_c - X_\theta|^2 \dd \theta &= \sum_i \int_{\theta_i}^{\theta_{i+1}}
    \frac{\dd \theta}{(a_i \cos \theta + b_i \sin\theta)^2}
    = \sum_i (a_i^2 \cot \theta + a_i b_i)^{-1} \Big|_{\theta_i}^{\theta_{i+1}}
    \\
    \int_0^{2\pi} |X_c - X_\theta|\e^\perp \dd \theta
    &= \sum_i \int_{\theta_i}^{\theta_{i+1}}
    \begin{pmatrix}-\sin\theta\\\cos\theta\end{pmatrix}
    \frac{\dd \theta}{(a_i \cos \theta + b_i \sin\theta)}
    \\
    &= \sum_i (a_i^2+b_i^2)^{-1}
    \begin{pmatrix}
        - a_i \log(r_i(\theta)) - b_i \theta
        \\
        - b_i \log(r_i(\theta)) + a_i \theta
    \end{pmatrix}\Bigg|_{\theta_i}^{\theta_{i+1}}
\end{align*}

In order to integrate function over the cortex $f(\Pi_{MT}(\theta)) = f(s_{MT}(\theta))$ w.r.t. $\theta$, we need to know $\Pi_{MT}(\theta)$ and ${\Pi_{MT}^{-1}}'(s)$
in order to perform the change of variables
\begin{equation}
    \label{eq:integration cortex theta}
    \int_{\theta_i}^{\theta_{i+1}} f(\Pi_{MT}(\theta)) \dd \theta = \int_{s_i}^{s_{i+1}} f(s) |{\Pi_{MT}^{-1}}'(s)| \dd s
\end{equation}
The expression of $|{\Pi_{MT}^{-1}}'(s)|$ is derived from \eqref{eq:MT density}:
\begin{align*}
    |{\Pi_{MT}^{-1}}'(s)| &= \frac{|n(s)\cdot(X(s) - X_c)|}{|X(s)-X_c|^2}
= \left|(X_{i+1} - X_i)^\perp\cdot \frac{X_i - X_c + \lambda(X_{i+1} - X_i)}{|X_i - X_c + \lambda(X_{i+1} - X_i)|^2}\right|
\\
&= \frac{|(X_{i+1} - X_i)^\perp\cdot (X_i - X_c)|}{|X_i - X_c + \lambda(X_{i+1} - X_i)|^2}
\end{align*}

If in \eqref{eq:integration cortex theta} we take a piecewise linear approximation for $f$, we
can compute the corresponding integral if we can compute
\begin{equation*}
    \int_0^1 \frac{\lambda^\alpha \dd \lambda}{|X_i-X_c + \lambda(X_{i+1}-X_i)|^2}\,, \quad \alpha\in \left\{0,1\right\}\,,
\end{equation*}
which can be put in the form
\begin{equation*}
    \int_0^1 \frac{\lambda^\alpha \dd \lambda}{c_0 + c_1 \lambda + c_2 \lambda^2}
\end{equation*}
where
\begin{align*}
    c_0 = |X_i - X_c|^2\,,\quad
    c_1 = 2(X_{i+1}-X_i)\cdot(X_i - X_c)\,,\quad
    c_2 = |X_{i+1} - X_i|^2
\end{align*}

\begin{align*}
    \alpha = 0&: \frac{2\tan^{-1}\frac{2c_2\lambda + c_1}{\sqrt{4c_2c_0-c_1^2}}}{\sqrt{4c_2c_0-c_1^2}}
    \\
    \alpha = 1&: \frac{1}{2c_2}\left[\log\left(\lambda(c_2\lambda+c_1)+c_0\right)-\frac{2c_1 \tan^{-1}\frac{2c_2\lambda +c_1}{\sqrt{4c_2c_0-c_1^2}}}{\sqrt{4c_2c_0-c_1^2}}\right]
    \\
    \alpha = 2&: \frac{1}{2c_2^2} \left[2c_2\lambda - c_1\log\left(\lambda(c_2\lambda+c_1)+c_0)\right) + \frac{2\left(c_1^2-2c_2c_0\right)\tan^{-1}\frac{2c_2\lambda+c_1}{\sqrt{4c_2c_0-c_1^2}}}{\sqrt{4c_2c_0-c_1^2}} \right]
\end{align*}

The last expression we need to compute is of the kind $\int |X_\theta - X_c| \e^\perp f(s_{MT}(\theta)) \dd \theta$.
To do so, we notice that
\begin{equation*}
    |\Pi_{MT}^{-1}(s)||X_c - X_\theta|\e^\perp =
    |(X_{i+1} - X_i)^\perp\cdot (X_i - X_c)| \frac{(X_i - X_c + \lambda \left(X_{i+1} - X_i\right))^\perp}{|X_i - X_c + \lambda(X_{i+1} - X_i)|^2}\,,
\end{equation*}
and the corresponding integral can be computed directly, using the case $\alpha = 2$ above.

\section{The nuclear membrane evolution equations}

\subsection{The formulation of Mikula-\v{S}ev\v{c}ovi\v{c}} \label{appendix:Mikula}

The nuclear membrane $\Gamma_n(t)$ can be either parameterized by arc length $\ell$ or by $\sigma$ which always span the same interval $[0, 1]$
\begin{equation*}
    \Gamma_n(t) := \left\{Y(t, \ell) : \ell \in [0, L_n(t)] \right\} = \left\{ Y(t, \sigma) : \sigma \in [0, 1] \right\} \, .
\end{equation*}
Therefore we can define the local length element $g := |\partial_\sigma Y(\sigma)|$, i.e., $\dd \ell= g \dd \sigma$
The curvature is denoted by $K$ and $\nu$ is the tangential angle, so that $T = (\cos \nu, \sin \nu)$ and $N = (\sin \nu, -\cos \nu)$. The region enclosed by $\Gamma_n$
id called $\Omega_n$ and $\dd A$ denotes its surface measure.

The curve's total energy $E_n$ is given by:
\begin{equation*}
    E_n^{TOT} = \En^{(1)}+ \En^{(2)}+  \En^{(3)} +  \int_{\Gamma_n} W(Y(\ell)) \dd \ell\,, \label{eq:nucleus_energy_Mikula}
\end{equation*}
where the first three terms are given by \eqref{eq:nucleus_energy expanded}, while $W$ denotes a potential, which combines the contribution related to the nuclear membrane surface tension (i.e. $\En^{(4)}$ given by \eqref{eq:nucleus_energy expanded}), the contact interaction with the cell cortex and the elastic constraint linking the nucleus with the centrosome.
The evolution of the planar curve $\Gamma_n$ is, then, given by \cite{mikula_computational_2004}:
\begin{equation}
\begin{cases}
    \partial_t K &= -\partial_\ell^2 \beta + \alpha \partial_\ell K - K^2 \beta
    \\
    \partial_t \nu &= -\partial_\ell \beta + \alpha K
    \\
    \partial_t g &= g K \beta + g\partial_\ell \alpha
    \\
    \partial_t Y &= \alpha T + \beta N,,
\end{cases}
    \label{eq:nucleus generic}
\end{equation}
where $K$, $\nu$ and $g = |\partial_\sigma Y|$ are the curvature, tangential angle and local length element, respectively.
Keeping in mind the previous governing equations, it is possible to calculate the time derivative of the energy functional $E$, formally given by
\begin{align}
    \label{eq:derivative nucleus energy}
    \frac{d}{dt} E_n^{TOT} &= -k_b \int_{\Gamma_n} (\partial_\ell^2 K + \frac{1}{2}K^3) \beta + \Delta p_n \int_{\Gamma_n} \beta  + \mu_n \int_{\Gamma_n} (A_n -A_{n}^*) \beta  \\
    &+ \int_{\Gamma_n} \nabla W \cdot (\alpha T + \beta N) + \int_{\Gamma_n} (W K \beta - \alpha \partial_\ell W )  = \nonumber \\
  &= -k_b \int_{\Gamma_n} (\partial_\ell^2 K + \frac{1}{2}K^3) \beta + \Delta p_n \int_{\Gamma_n} \beta + \int_{\Gamma_n} (\nabla W \cdot N)\beta   + \int_{\Gamma_n} W K \beta  \,.
\end{align}
Then, the choice of $\beta$ is made to maximize the decrease of energy due to the normal component, whereas the tangential velocity $\alpha$, in principle, is a free parameter. In this paper, we refer to the asymptotically uniform tangential redistribution derived in \cite{mikula_computational_2004, benes_comparison_2009}, to impose the evolution of $\alpha$ in order to avoid nodes concentrating on points, leading to poor approximation and eventually inversion of ill-conditioned matrices.
Therefore ,we take
\begin{align}
   \partial_\ell \alpha &= -K\beta + \langle K \beta \rangle + (L_n/g-1) \zeta
    \label{eq:nucleus tangent velocity}
    \\
    \beta &= k_b \left(\partial_\ell^2 K + \frac{1}{2} K^3\right) - \Delta p_n  - \mu_n (A_n -A_{n}^*)  - \nabla W \cdot N - W K \,,    \label{eq:nucleus normal velocity}
\end{align}
where $\langle.\rangle = L_n^{-1} \int_{\Gamma_n} . \; \dd \ell$  is an averaging operator over the whole curve $\Gamma_n$, which makes \eqref{eq:nucleus tangent velocity} a non-local equation, whereas $\zeta >0$ is a given positive constant.

Recalling the identity:
\begin{align*}
    \partial_\ell^4 Y &= -(\partial_\ell^2 K) N - 2 \partial_\ell K \partial_\ell N - K \partial_\ell^2 N = -(\partial_\ell^2 K) N - 2 K \partial_\ell K T - K \partial_\ell (K T)
    \\
    &= -(\partial_\ell^2 K) N - 3 K \partial_\ell K T - K^2 \partial_\ell T = -(\partial_\ell^2 K) N - \frac{3}{2}\partial_\ell K^2 \partial_\ell Y - K^2 \partial_\ell^2 Y
\end{align*}
and $\partial_\ell^2 Y = -K N$, we have
\begin{align*}
    \left(\partial_\ell^2 K + \frac{1}{2}K^3\right) N &= -\partial_\ell^4 Y - \frac{3}{2} K^2 \partial_\ell^2 Y - \frac{3}{2} \partial_\ell K^2 \partial_\ell Y
    \\
    &= -\partial_\ell^4 Y - \frac{3}{2} \partial_\ell \left(K^2 \partial_\ell Y\right) \, ,
\end{align*}
so that, setting $g = \exp \eta$, we can rewrite \eqref{eq:nucleus generic} as
\begin{equation}
    \begin{cases}
    \partial_t K &= - k_b \left(\partial_\ell^4 K + \frac{1}{2} \partial_\ell^2(K^3)\right) + \partial_\ell (\alpha K) - K(K\beta + \partial_\ell \alpha) + \partial_\ell^2 (\nabla W \cdot N) + \partial_\ell^2 (W K)
    \\
     \partial_t \nu &= -k_b \partial_\ell^4 \nu - \frac{k_b}{2}\partial_\ell(\partial_\ell \nu)^3 + \partial_\ell(\nabla W\cdot N) + \partial_\ell (W \partial_\ell \nu ) + \alpha \partial_\ell \nu\
    \\
    \partial_t \eta &= K \beta + \partial_\ell \alpha
    \\
       \partial_t Y &= \alpha \partial_\ell Y - k_b \left(\partial_\ell^4 Y + \frac{3}{2} \partial_\ell(K^2 \partial_\ell Y)\right) + W \partial_\ell^2 Y - \Delta p_n \, N - \mu_n (A_n -A_{n}^*) N - \left( \nabla W \cdot N \right)  N
\label{eq:nucleus specialized}
\end{cases}
\end{equation}
We remark that the second equation in \eqref{eq:nucleus specialized} has been obtained remembering that $K = \partial_\ell \nu$.

\subsection{Quadrature formulae for the Mikula formulation} \label{App:Mikula_discretized}

We can now write the corresponding discretized equations of \eqref{eq:nucleus specialized}. To do so, we uniformly discretized the fixed parametrization interval $[0, 1]$
in $N$ subintervals, each of equal length $h=1/N_2$ and indexed by $i \in [0, N_2-1]$. Time is also discretized with
a time-step $\Delta t$. The discrete time is indexed by $j\in \mathbb{N}$, so that $t_j = j\, \Delta t$.
In what follows, indices corresponding to space (resp. time) will be in subscript (resp. superscript),
so that the point $Y(ih, j\Delta t)$ is written $Y_i^j$. The measure of the finite element $[Y_{i-1}^j, Y_i^j]$ at time $t_j$, is given by $r_i = |Y_{i-1} - Y_i|$. We also introduce the dual finite element $[\tilde{Y}_{i}^j, \tilde{Y}_{i+1}^j]$, containing the node $Y_{i}^j$, being $\tilde{Y}_{i}^j= \frac{Y_{i-1}^j+Y_i^j}{2}$. Le length of the dual finite element is denoted by $q_i^j = \frac{1}{2} \left(r_i^j + r_{i+1}^j\right)$.
Then, the system of equations  \eqref{eq:nucleus tangent velocity}- \eqref{eq:nucleus normal velocity}-\eqref{eq:nucleus specialized} is solved for the discrete quantities $\alpha_i^j$, $\beta_i^j$, $K_i^j$, $\nu_i^j$, $\eta_i^j$, $Y_i^j$.
In particular $\alpha_i^j$ denotes the tangential velocity of the node $Y_i^j$, whereas $\beta_i^j$, $K_i^j$, $\nu_i^j$, $\eta_i^j$ are piecewise constant approximations of the corresponding quantities on the finite element $[Y_{i-1}^j, Y_i^j]$.
The dependent variables $W(Y)$ and $\nabla W(Y)$ are consequently evaluated numerically at nodes, i.e. $W_i^j=W(Y_i^j)=W(Y(ih, j\Delta t))$, whereas the tangential and normal vectors are approximated on the finite element $[Y_{i-1}^j, Y_i^j]$, by the quantities $T_i^j = (\cos \nu_i^j, \sin \nu_i^j)^T$ and $N_i^j = (\sin \nu_i^j, - \cos \nu_i^j)^T$.
\\
We first integrate \eqref{eq:nucleus tangent velocity} on the finite element $[Y_{i-1}, Y_i]$
\begin{align*}
    \int_{Y_{i-1}}^{Y_i} \partial_\ell \alpha \dd \ell= \alpha_{i} - \alpha_{i-1} &= \int_{Y_{i-1}}^{Y_i} -K \beta + \langle K \beta \rangle + \left(L/g - 1\right)\omega \, \dd \ell
    \\
    & \simeq r_i \left( -K_i \beta_i + \langle K \beta\rangle\right) + \left(L/N_2 - r_i\right) \omega\,.
\end{align*}
and we discretize in time to get
\begin{equation*}
    \alpha_i^j = \alpha_{i-1}^j
    + r_i^{j-1} \left( -K_i^{j-1} \beta_i^{j-1} + B^{j-1}\right) + \left(L^{j-1}/ N_2 - r_i^{j-1}\right) \omega\,,
\end{equation*}
being $L^j = \sum_l r_l^j\,, \quad  B^j =  (L^j)^{-1} \sum_l r_l^j K_l^j \beta_l^j$ and $\beta_i^j$ obtained by the discretization of \eqref{eq:nucleus normal velocity}, i.e.,
\begin{gather*}
    \beta_i^j = \frac{k_b} {r_i^j} \left(\frac{K_{i+1}^j - K_i^j}{q_i^j} - \frac{K_i^j - K_{i-1}^j}{q_{i-1}^j}\right) + \frac{1}{2} k_b (K_i^j)^3
    - \Delta p_n - N_i^{j}\cdot \nabla \tilde{W}_i^{j} - \tilde{W}_i^{j} K_i^{j}\,.
\end{gather*}
This determines $\alpha_i^j$ for all $j$, provided a closure for $\alpha_0^j$.
It is possible to assume $\alpha_0^j=0$, i.e., the point $Y_0^j$ moves along the normal \cite{benes_comparison_2009}.
The local length of the finite element is updated in a similar way, using \eqref{eq:nucleus specialized}$_{c}$:
\begin{equation*}
    r_i^{j-1} \frac{\eta_i^j - \eta_i^{j-1}}{\Delta t} = r_i^{j-1} K_i^{j-1} \beta_i^{j-1} + \alpha_i^j - \alpha_{i-1}^j\,,
\end{equation*}
after which we set $r_i^j = \exp(\eta_i^j)$.
Integrating eq. \eqref{eq:nucleus specialized}$_{a}$ on $[Y_{i-1}, Y_i]$ we get
\begin{align}
    \int_{Y_{i-1}}^{Y_i} \partial_t K &= \int_{Y_{i-1}}^{Y_i} -k_b \partial_\ell^4 K - \frac{k_b}{2} \partial_\ell^2 (K^3) + \partial_\ell (\alpha K) - K(K\beta + \partial_\ell \alpha) + \partial_\ell^2 (\nabla W\cdot N) + \partial_\ell^2 (W K) \, \dd \ell\nonumber
    \\
    r_i^j \partial_t K_i &= -k_b [\partial_\ell^3 K]_{Y_{i-1}}^{Y_i} -\frac{k_b}{2} [\partial_\ell (K^3)]_{Y_{i-1}}^{Y_i}
    + [\alpha K]_{Y_{i-1}}^{Y_i} - K_i(K_i r_i \beta_i + (\alpha_i - \alpha_{i-1})) \nonumber  \\
    &+ [\partial_\ell ( \nabla W \cdot N)]_{Y_{i-1}}^{Y_i} + [\partial_\ell (W K)]_{Y_{i-1}}^{Y_i}= \nonumber \\
    & -k_b [\partial_\ell^3 K]_{Y_{i-1}}^{Y_i} -\frac{k_b}{2} [\partial_\ell (K^3)]_{Y_{i-1}}^{Y_i}
    + [\alpha K]_{Y_{i-1}}^{Y_i} - K_i(K_i r_i \beta_i + (\alpha_i - \alpha_{i-1}))  \nonumber \\
    &+ [\partial_\ell ( \nabla W \cdot N)]_{Y_{i-1}}^{Y_i} + [\partial_\ell (W K)]_{Y_{i-1}}^{Y_i}
    \label{eq:discretized k}
\end{align}
If we approximate the time derivative as $\Delta t^{-1} (K_i^j - K_i^{j-1})$, $K_i^j$ solves a linear system of the form
\begin{equation}
    a_i^j K_{i-2}^j + b_i^j K_{i-1}^j + c_i^j K_i^j + d_i^j K_{i+1}^j + e_i^j K_{i+2}^j = f_i^j\,.
    \label{eq:nucleus linear system k}
\end{equation}
Then, we have to solve the tangent angle equation \eqref{eq:nucleus specialized}$_{b}$, using the following approximation
\begin{align}
    \int_{Y_{i-1}}^{Y_i} \partial_t \nu &= \int_{Y_{i-1}}^{Y_i} - k_b \partial_\ell^4 \nu - \frac{k_b}{2} \partial_\ell (\partial_\ell \nu)^3 + \partial_\ell (\nabla W \cdot N)  + \partial_\ell (W \partial_\ell \nu) + \alpha \partial_\ell \nu \, \dd \ell
\nonumber    \\
    r_i^j \partial_t \nu_i &= -k_b [\partial_\ell^3 \nu]_{Y_{i-1}}^{Y_i} -\frac{k_b}{2} [(\partial_\ell \nu)^3]_{Y_{i-1}}^{Y_i}
    + [\nabla W \cdot N ]_{Y_{i-1}}^{Y_i} + [\alpha \nu]_{Y_{i-1}}^{Y_i}  + [W \partial_\ell \nu]_{Y_{i-1}}^{Y_i} - \nu_i (\alpha_i - \alpha_{i-1})
        \label{eq:discretized nu}
\end{align}
and approximating the time differential by finite differences, we obtain the following system for $\nu_i^j$
\begin{equation}
    \mathcal{A}_i^j \nu_{i-2}^j + \mathcal{B}_i^j \nu_{i-1}^j + \mathcal{C}_i^j \nu_i^j + \mathcal{D}_i^j \nu_{i+1}^j + \mathcal{E}_i^j \nu_{i+2}^j = \mathcal{F}_i^j\,.
    \label{eq:nucleus linear system nu}
\end{equation}

Finally, the discretization of the equation for $Y$ \eqref{eq:nucleus specialized}$_{d}$ is obtained by integrating on the dual element $[\tilde{Y}_{i+1}, \tilde{Y}_i]$, i.e.,
\begin{align*}
    \int_{\tilde{Y}_{i}}^{\tilde{Y}_{i+1}} \partial_t Y \dd \ell &=
    \int_{\tilde{Y}_{i}}^{\tilde{Y}_{i+1}} \alpha \partial_\ell Y - k_b \left(\partial_\ell^4 Y + \frac{3}{2} \partial_\ell (K^2 \partial_\ell Y)\right) + W \partial_\ell^2 Y
    - \Delta p_n \, N  \\
    & - \mu_n (A_n -A_{n}^*) N  - \left( \nabla W \cdot N \right)  N  \, \dd \ell \,
\end{align*}
So that the complete semi-discrete equation is:
\begin{multline}
    q_{i} \partial_t Y_i =  \alpha_i (\tilde{Y}_{i+1}-\tilde{Y}_{i}) -k_b [\partial_\ell^3 Y]_{\tilde{Y}_{i}}^{\tilde{Y}_{i+1}} - \frac{3}{2} k_b  [K^2 \partial_\ell Y]_{\tilde{Y}_{i}}^{\tilde{Y}_{i+1}}+ W_ i[ \partial_\ell Y]_{\tilde{Y}_i}^{\tilde{Y}_{i+1}}    \\
    - \frac{\Delta p_n  + \mu_n (A_n -A_{n}^*) }{2} \left( r_{i+1}  N_{i+1}  + r_{i} N_{i}  \right)
- \frac{1}{2}(r_{i+1}(N_{i+1} \cdot \nabla W_{i})N_{i+1} + r_i(n_i \cdot \nabla W_i)N_i) \,
    \label{eq:discretized Y}
\end{multline}
Again approximating the time differential by finite differences, we get a system of the following form for $Y_i^j$:
\begin{equation}
    A_i^j Y_{i-2}^j + B_i^j Y_{i-1}^j + C_i^j Y_i^j + D_i^j Y_{i+1}^j + E_i^j Y_{i+2}^j = F_i^j\,.
    \label{eq:nucleus linear system Y}
\end{equation}
The coefficients of the linear equations \eqref{eq:nucleus linear system k}, \eqref{eq:nucleus linear system nu} and \eqref{eq:nucleus linear system Y}  are reported in the following.

\subsubsection{Discretization of the equation for \texorpdfstring{$K$}{K}}
In the following we report how the different terms appearing in eq.  \eqref{eq:discretized k} have been discretized
\begin{align*}
[\partial_\ell^3 K]_{i-1}^i & \approx \left[ \frac{\partial_\ell^2 K (\tilde{Y}_{m+1})-\partial_\ell^2 K (\tilde{Y}_{m})}{q_m} \right]_{i-1}^i  \\
& \approx
\left[ \frac{1}{q_m} \left(  \frac{\partial_\ell K (Y_{m+1})-\partial_\ell K (Y_{m})}{r_{m+1}} \right) - \frac{1}{q_i} \left(  \frac{\partial_\ell K (Y_{m})-\partial_\ell K (Y_{m-1})}{r_{m}} \right) \right]_{i-1}^i \\
& \approx
\left[ \frac{1}{q_m r_{m+1}} \left(  \frac{K_{m+2}-K_{m+1}}{q_{m+1}} - \frac{K_{m+1}-K_{m}}{q_{m}} \right)
-\frac{1}{q_m r_{m}} \left(  \frac{K_{m+1}-K_{m}}{q_{m}} - \frac{K_m-K_{m-1}}{q_{m-1}} \right)\right]_{i-1}^i \\
[\partial_\ell K^3]_{i-1}^i & \approx\left[ \frac{K_{m+1}^3 - K_{m}^3}{q_m} \right]_{i-1}^i \\
[\alpha K]_{i-1}^i &\approx [\alpha_m \tilde K_m]_{i-1}^i = \left[ \alpha_m  \frac{K_{m+1}+K_{m}}{2}\right] _{i-1}^i \\
[\partial_\ell (\nabla W\cdot N)]_{i-1}^i & \approx [\frac{1}{q_m}(\nabla \tilde{W}_{m+1} \cdot N_{m+1} - \nabla \tilde{W}_{m} \cdot N_{m})]_{i-1}^i \\
&= [\frac{1}{2q_m}((\nabla W_{m+1} + \nabla W_m) \cdot N_{m+1} - (\nabla W_{m} + \nabla W_{m-1})\cdot N_{m})]_{i-1}^i \\
[\partial_\ell(KW)]_{i-1}^i & \approx \left[ \frac{\tilde{W}_{m+1}K_{m+1} - \tilde{W}_{m} K_{m}}{q_m} \right] _{i-1}^i =
\left[ \frac{1}{2 q_m} \left( \left(W_{m+1}+W_{m}\right) K_{m+1} - \left(W_{m}+W_{m-1} \right) K_{m} \right) \right] _{i-1}^i
\end{align*}

Therefore, the coefficient of  \eqref{eq:nucleus linear system k} are given as follows:
\begin{align*}
a_i^j &= \dfrac{k_b}{q_{i-2}^j r_{i-1}^j q_{i-1}^j} \, , \quad
e_i^j = \dfrac{k_b}{q_{i}^j r_{i+1}^j q_{i+1}^j} \\
b_i^j &= -k_b \left( \dfrac{1}{q_{i-1}^j r_{i}^j q_{i}^j} +  \dfrac{1}{ (q_{i-1}^j)^2 r_{i}^j }+  \dfrac{1}{ (q_{i-1}^j)^2 r_{i-1}^j} + \dfrac{1}{q_{i-1}^j r_{i-1}^j q_{i-2}^j}  \right) + \dfrac{\alpha_{i-1}^j}{2} - \dfrac{1}{2 q_{i-1}^j} (W_{i-1}^{j-1}+W_{i-2}^{j-1}) \\
d_i^j &= -k_b \left( \dfrac{1}{q_{i}^j r_{i+1}^j q_{i+1}^j} +  \dfrac{1}{ (q_{i}^j)^2 r_{i+1}^j }+  \dfrac{1}{ (q_{i}^j)^2 r_{i}^j} + \dfrac{1}{q_{i}^j r_{i}^j q_{i-1}^j}  \right) - \dfrac{\alpha_{i}^j}{2} - \dfrac{1}{2 q_{i}^j} (W_{i+1}^{j-1}+W_{i}^{j-1})  \\
c_i^j &= \dfrac{r_i^j}{\Delta t} + k_b \left( \dfrac{1}{ (q_{i}^j)^2 r_{i+1}^j } +  \dfrac{1}{ (q_{i}^j)^2 r_{i}^j} +  \dfrac{1}{ (q_{i-1}^j)^2 r_{i}^j} +  \dfrac{1}{ (q_{i-1}^j)^2 r_{i-1}^j} + \dfrac{2}{q_{i}^j r_{i}^j q_{i-1}^j}  \right) + \dfrac{\alpha_{i}^j}{2} - \dfrac{\alpha_{i-1}^j}{2} \\
&+ r_i^{j-1} k_i^{j-1} \beta_i^{j-1} + \dfrac{1}{2}\left( \dfrac{1}{q_i^j}+\dfrac{1}{q_{i-1}^j} \right) (W_i^{j-1}+W_{i-1}^{j-1}) \\
f_i^j &=
\dfrac{r_i^j}{\Delta t} k_{i}^{j-1}
    - k_b \dfrac{ (k_{i+1}^{j-1})^3 - (k_{i}^{j-1})^3 }{2  q_{i}^{ j-1}}
    + k_b \dfrac{ (k_i^{j-1})^3 - (k_{i-1}^{j-1})^3 }{2  q_{i-1}^{j-1}}
    +\dfrac{\nabla W_{i+1}^{j-1} + \nabla W_{i}^{j-1}}{2q_i^{ j-1}} \cdot N_{i+1}^{j-1} \\
      &-\left( \dfrac{1}{2q_i^{j-1}} +\dfrac{1}{2q_{i-1}^{j-1}}\right) (\nabla W_{i}^{j-1}+ \nabla W_{i-1}^{j-1}) \cdot N_i^{j-1}
      + \dfrac{\nabla W_{i-1}^{j-1} + \nabla W_{i-2}^{j-1}}{2q_{i-1}^{ j-1}} \cdot N_{i-1}^{j-1}
\end{align*}

\subsubsection{Discretization of the equation for \texorpdfstring{$\nu$}{nu}}
In the following we report how the different terms appearing in eq.  \eqref{eq:discretized nu} have been discretized
\begin{align*}
[(\partial_\ell \nu)^3]_{i-1}^i &\approx  [(\frac{1}{q_m}(\nu_{m+1} - \nu_{m}))^3]_{i-1}^i\\
[\nabla W \cdot N]_{i-1}^i &\approx \left[ \nabla W_m \cdot \widetilde{N_m} \right] _{i-1}^i = \left[ \frac{1}{2} \nabla W_m \cdot  \left( N_{m+1} + N_{m} \right) \right] _{i-1}^i \\
[\alpha \nu]_{i-1}^i &\approx \left[ \alpha_{m}\tilde\nu_{m} \right] _{i-1}^i = \left[ \frac{1}{2} \alpha_{m} (\nu_{m}+ \nu_{m+1} \right)] _{i-1}^i \\
[W \partial_\ell \nu]_{i-1}^i &\approx [\frac{ W_m}{q_m}(\nu_{m+1}-\nu_m)]_{i-1}^i
\end{align*}

Then, the coefficient of  \eqref{eq:nucleus linear system nu}  are
\begin{align*}
\mathcal{A}_i^j &= \dfrac{k_b}{q_{i-2}^j r_{i-1}^j q_{i-1}^j}  \, , \quad
\mathcal{E}_i^j = \dfrac{k_b}{q_{i}^j r_{i+1}^j q_{i+1}^j} \, ,\\
\quad \mathcal{B}_i^j &=-k_b \left( \dfrac{1}{q_{i-1}^j r_{i}^j q_{i}^j} +  \dfrac{1}{ (q_{i-1}^j)^2 r_{i}^j }+  \dfrac{1}{ (q_{i-1}^j)^2 r_{i-1}^j} + \dfrac{1}{q_{i-1}^j r_{i-1}^j q_{i-2}^j}  \right) + \dfrac{\alpha_{i-1}^j}{2} - \dfrac{W_{i-1}^{j-1}}{q_{i-1}^j}\, , \\
\mathcal{D}_i^j &= -k_b \left( \dfrac{1}{q_{i}^j r_{i+1}^j q_{i+1}^j} +  \dfrac{1}{ (q_{i}^j)^2 r_{i+1}^j }+  \dfrac{1}{ (q_{i}^j)^2 r_{i}^j} + \dfrac{1}{q_{i}^j r_{i}^j q_{i-1}^j}  \right) - \dfrac{\alpha_{i}^j}{2}  - \dfrac{W_i^{j-1}}{q_{i}^j} \, ,\\
\mathcal{C}_i^j &= \dfrac{r_i^j}{\Delta t} + k_b \left( \dfrac{1}{ (q_{i}^j)^2 r_{i+1}^j } +  \dfrac{1}{ (q_{i}^j)^2 r_{i}^j} +  \dfrac{1}{ (q_{i-1}^j)^2 r_{i}^j} +  \dfrac{1}{ (q_{i-1}^j)^2 r_{i-1}^j} + \dfrac{2}{q_{i}^j r_{i}^j q_{i-1}^j}  \right) \\
&+ \dfrac{\alpha_{i}^j}{2} - \dfrac{\alpha_{i-1}^j}{2} + \left( \dfrac{W_i^{j-1}}{q_{i}^j} + \dfrac{W_{i-1}^{j-1}}{q_{i-1}^j}\right) \\
\mathcal{F}_i^j &= \dfrac{r_i^j}{\Delta t} \nu_{i}^{j-1} - \dfrac{k_b}{2} \left(\dfrac{ \nu_{i+1}^{j-1} - \nu_{i}^{j-1}}{q_{i}^{ j-1}} \right)^3 + \dfrac{k_b}{2} \left( \dfrac{ \nu_i^{j-1} - \nu_{i-1}^{j-1} }{q_{i-1}^{j-1}} \right)^3 \\
&
+ \frac{1}{2}\nabla W_{i}^{j-1} \cdot (N_{i+1}^{j-1} + N_i^{j-1})
+ \frac{1}{2}\nabla W_{i-1}^{j-1} \cdot (N_{i}^{j-1} + N_{i-1}^{j-1})
\end{align*}

\subsubsection{Discretization of the equation for \texorpdfstring{$Y_i$}{Yi}}
In the following we report how the different terms appearing in eq.  \eqref{eq:discretized Y} have been discretized
$$ W_i [\partial_\ell Y]_{\widetilde{i}}^{\widetilde {i+1}} = W_i [\frac{Y_{m}-Y_{m-1}}{r_{m}}]_{i}^{i+1}$$

The coefficient of   \eqref{eq:nucleus linear system Y} are
\begin{align*}
A_i^j &= \dfrac{k_b}{r_{i-1}^j q_{i-1}^j r_{i}^j} \, , \quad
E_i^j = \dfrac{k_b}{r_{i+1}^j q_{i+1}^j r_{i+2}^j} \\
B_i^j &= -k_b \left( \dfrac{1}{r_{i-1}^j q_{i-1}^j r_{i}^j} +  \dfrac{1}{ (r_{i}^j)^2 q_{i-1}^j }+  \dfrac{1}{ (r_{i}^j)^2 q_{i}^j} + \dfrac{1}{r_{i}^j q_{i}^j r_{i+1}^j}  \right) + \dfrac{3}{2}k_b \dfrac{(k_i^j)^2}{r_i^j}+\dfrac{\alpha_{i}^j}{2} - \dfrac{W_i^{j-1}}{r_i^j}  \\
D_i^j &= -k_b \left( \dfrac{1}{r_{i}^j q_{i}^j r_{i+1}^j} +  \dfrac{1}{ (r_{i+1}^j)^2 q_{i}^j }+  \dfrac{1}{ (r_{i+1}^j)^2 q_{i+1}^j} + \dfrac{1}{r_{i+1}^j q_{i+1}^j r_{i+2}^j}  \right) + \dfrac{3}{2}k_b \dfrac{(k_{i+1}^j)^2}{r_{i+1}^j} - \dfrac{\alpha_{i}^j}{2}  - \dfrac{W_i^{j-1}}{r_{i+1}^j}  \\
C_i^j &= \dfrac{ q_i^j}{\Delta t} - \left( A_i^j + B_i^j +D_i^j + E_i^j\right) \\
F_i^j &= \dfrac{q_i^j}{\Delta t} y_{i}^{j-1}
- \frac{\Delta p_n   + \mu_n (A_n -A_{n}^*) }{2} \left( r_{i+1}^{j-1}  N_{i+1}^{j-1} + r_i^{j-1} N_i^{j-1} \right) \\
      &- \frac{1}{2}(r_{i+1}^{j-1}(N_{i+1}^{j-1} \cdot \nabla W_{i}^{j-1})N_{i+1}^{j-1} + r_i^{j-1}(N_i^{j-1} \cdot \nabla W_i^{j-1})N_i^{j-1}) \,
\end{align*}

We remark that linear dependency on $Y$ inside $W(Y)$ can be possibly treated implictly.

\bibliography{GJPS.bib}
\bibliographystyle{abbrv}

\end{document}